\documentclass{article}

\usepackage[letterpaper,top=2cm,bottom=2cm,left=3cm,right=3cm,marginparwidth=1.75cm]{geometry}
\usepackage{graphicx,  subcaption, upgreek, enumerate}
\usepackage{comment}

\usepackage{bm, xcolor}
\usepackage{comment}
\usepackage{amsfonts,amsmath,amssymb,amsthm}
\usepackage{hyperref,url}
\newtheorem{theorem}{Theorem}
\newtheorem{lemma}{Lemma}

\newtheorem{definition}{Definition}
\newtheorem{proposition}{Proposition}
\newtheorem{remark}{Remark}
\newtheorem{example}{Example}

\newcommand{\rg}[1]{{\color{blue} #1}}

\newcommand{\ys}[1]{{\color{brown} #1}}

\newcommand{\matr}[1]{\mathbf{#1}}
\newcommand{\barx}{\overline{x}}

\allowdisplaybreaks

\title{Normal approximation for the posterior in exponential families}
\author{Adrian Fischer\footnote{Adrian Fischer, University of Oxford, UK. E-mail: adrian.fischer@stats.ac.uk}, Robert
  E. Gaunt\footnote{Robert E. Gaunt, The University of Manchester,
    UK. E-mail: robert.gaunt@manchester.ac.uk}, Gesine
  Reinert\footnote{Gesine Reinert, University of Oxford, UK. E-mail:
    reinert@stats.ox.ac.uk}, and Yvik Swan\footnote{ Yvik Swan,
    Universit\'e libre de Bruxelles, Belgium. E-mail:
    yvik.swan@ulb.be}}

\begin{document}

\maketitle

\begin{abstract}
  In this paper we obtain quantitative, non-asymptotic, and data dependent,  {\it Bernstein-von Mises type} bounds on the normal approximation of the posterior distribution in exponential family models with arbitrary centring and scaling. Our bounds, {which are stated in the total variation and Wasserstein distances}, are valid for univariate and multivariate posteriors alike, and do not require a conjugate prior setting. 
  {They} are obtained through a {refined version of Stein's method of comparison of operators
  that allows for improved dimensional dependence in high-dimensional settings, and {that} may also be of interest in other problems}. Our approach is 
  rather flexible and in certain settings allows for the derivation of bounds with rates of convergence faster than the usual $O(n^{-1/2})$ rate (when $n$ is the sample size).  We illustrate our findings on a variety of exponential family distributions, including the 
  {Weibull}, multinomial and linear regression with unknown variance. The resulting bounds have an explicit dependence on the prior distribution and on sufficient statistics of the data from the sample, and thus provide insight into how these factors {may} affect the quality of the normal approximation. {Insights that can be gleaned from our examples include identification of conditions under which faster $O(n^{-1})$ convergence rates occur for Bernoulli data, 
  illustrations of how the quality of the normal approximation is influenced by the choice of standardisation, and dimensional dependence in high-dimensional settings.}
\end{abstract}

\


\
	

\section{Introduction}

The Bernstein-von Mises (BvM) Theorem is a cornerstone in Bayesian
statistics. Loosely put, this the- orem reconciles Bayesian and
frequentist large sample theory by guaranteeing that suitable scalings
of posterior distributions are asymptotically normal. In particular,
this implies that the contribution of the prior vanishes in the
asymptotic posterior. Beyond its 
philosophical implications, this
result is central to Bayesian inference, as it provides theoretical
justification for \emph{Laplace approximation} which consists in using
Gaussian approximations for (typically intractable) posterior
integrals. The BvM Theorem is named after the works \cite{bernnnn,
  VMRW}, but can 
  be traced back to work by Laplace
\cite{laplace1820theorie}; we refer to \cite[Section
8]{miller2021asymptotic} for an overview, including a historical
\emph{tour d'horizon}.

There exist many different versions of the BvM Theorem; one of the
most basic forms goes as follows.  Let
$\theta \in \mathbb{R}^d \sim \pi_0$ follow the prior distribution,
and conditional on $\theta$, let $X_1, X_2, \ldots$ be independent
random variables in $\mathbb{R}^r$ from  a probability measure $P_{\theta}$
 belonging to a model $(P_{\theta} : \theta \in \Theta)$.  For
$n \ge 1$ let $\mathbf{x} = (x_1, \ldots, x_n)$ be a realisation of $\mathbf{X} = (X_1, \ldots, X_n)$  and let
$P_{2}(\theta \, | \, \mathbf{x})$ denote the posterior distribution
of $\theta$ given  $\mathbf{X} = \mathbf{x}.$ 
Let
$\tilde \theta_n$ denote the posterior mode (we call it the MAP, for Maximum A Posteriori estimator) and $I(\tilde \theta_n)$
the model Fisher information matrix of $P_{\theta}$ evaluated at
$\tilde \theta_n$.  The BvM Theorem then assures us that, under
regularity conditions ({including} invertibility of the Fisher information {$I(\tilde \theta_n)$},
quadratic mean differentiability, identifiability
), 
\begin{equation}
  \label{eq:BVMTv1}
  d_{\mathrm{TV}}\big(P_2(\theta \, | \, \mathbf{X}) , \mathcal{N}(\tilde
  \theta_n , 
  I^{-1}(\tilde\theta_n){/n})\big) =
  o_{P_{\theta_0}^n}(1)
\end{equation}
where $d_{\mathrm{TV}}(P, Q)$ denotes the total variation distance
between distributions $P$ and $Q$, and the convergence holds in
probability with respect to the law
$P_{\theta_0}^n = P_{\theta_0}\otimes \dots \otimes P_{\theta_0}$ of
the data generating process when the true parameter is $\theta_0$.  A
similar conclusion holds for a centring around the maximum likelihood
estimate (MLE) and, more generally, around any efficient estimator of
the true parameter;  see \cite[Section 10.2]{van2000asymptotic} for
details and a proof.

Approximate normality of posterior distributions still holds, 
with a different limiting covariance matrix, in the misspecifed model setting that 
the true underlying distribution $P_0$ is no longer 
of the form
$P_{\theta_0}$.
An early reference on misspecified BvM Theorems is
\cite{kleijn2012bernstein} where the corresponding version of
\eqref{eq:BVMTv1} is obtained under a Local Asymptotic Normality (LAN)
condition; see also \cite{koers2023misspecified} for a BvM theorem in
the context of misspecified, non-i.i.d., hierarchical models, and 
\cite{bochkina2022bernstein,bochkina2023bernstein} for an overview
of BvMs under misspecification. Extensions to the semiparametric
setting have  been obtained e.g.\ in
\cite{bickel2012semiparametric} under a LAN condition, or in
\cite{castillo2015bernstein} for smooth functionals, both for i.i.d.\
and non-i.i.d.\ samples.  Versions of the BvM result also hold true in
nonparametric settings, see e.g.\ \cite{leahu2011bernstein,
  castillo2013nonparametric, castillo2014bernstein}. We refer to
\cite{rousseau2016frequentist} for an overview of 
semi- and
nonparametric BvM Theorems. The high-dimensional
setting 
is also well studied, see e.g.\ 
\cite{ghosal2000asymptotic,boucheron2009bernstein, johnstone2010high,
  bontemps2011bernstein, 
jin2024high, sarkar2023high}. Almost sure versions are obtained
in \cite{miller2021asymptotic} for 
various commonly encountered
situations. Finally, extensions of the BvM statement beyond normal
approximation have also been considered recently, see e.g.\
\cite{durante2023skewed} wherein skew normal approximations are
studied.

In practice, there is a need for non-asymptotic statements which,
arguably, are of a different nature than \eqref{eq:BVMTv1}. Here, the 
objective is to study,  for a fixed sample,  the proximity between the
standardised posterior distribution (conditioned on the observation)
and a standard normal distribution. An important reference on this
topic is \cite{panov2015finite} where a finite sample semiparametric
version of the BvM theorem is obtained for non-informative and flat
Gaussian priors; it is seen furthermore that, under certain conditions, the
normal approximation holds as soon as $d^3\ll n$ (which is an
improvement from \cite{ghosal2000asymptotic} where the condition
$d^4 \log d \ll n$ was identified). Similar spirited fixed sample
bounds are also obtained in \cite{spokoiny2019accuracy} for Gaussian
priors in high-dimensional or nonparametric setups, while
\cite{spokoiny2022dimension} considers Gaussian priors and log-concave
likelihoods; in 
\cite{dehaene2019deterministic} 
the
proximity is expressed in terms of the Kullback-Leibler
divergence. Recently, \cite{kasprzak2022good} proposed fixed sample bounds of the
form
\begin{equation}\label{eq:kasp}
  \mathcal{D} \big( P^{\star}(\theta
  \, | \, \mathbf{x}) , \mathcal{N}(0, I) \big) \le \mathcal{U}_n(\mathbf{x})
\end{equation}
with $\mathcal{D}(P, Q)$ standing for either the total variation distance $d_{\mathrm {TV}}(P, Q)$, 
the Wasserstein distance $d_{\mathrm {Wass}}(P, Q)$,  or the 
difference of covariances between laws $P$ and $Q$, with
$P^{\star}(\theta \, | \, \mathbf{x})$ the law of the standardised
posterior centred around the MAP or the MLE, and {with} 
$\mathcal{U}_n(\mathbf{x})$ a statistic.
The quantity $\mathcal{U}_n(\mathbf{x})$ is, in principle,
computable explicitly, 
although exact computation turns out to be challenging in any specific
example.  Concurrently, the preprints
\cite{katsevich2023improved,katsevich2023laplace} 
prove 
that the
condition $d^2 \ll n $ is necessary for total variation distance and obtain bounds 
with explicit dependence on the dimension for specific examples (multinomial, logistic regression). The papers  \cite{kasprzak2022good, katsevich2023improved,katsevich2023laplace}
also provide excellent, detailed, and comparative literature reviews
on BvM Theorems and Laplace approximation.


\medskip 
In this paper, we obtain bounds of the same spirit as \eqref{eq:kasp}, 
allowing to assess, at fixed sample size and sample, the quality of normal
approximation of the posterior under the assumption that the
observations are independent (but not necessarily identically
distributed) from an \emph{exponential family}.  We 
make no assumptions 
{(other than a mild differentiability assumption)} on the prior $\pi_0$, and consider general centring and
scaling. Specifically, the question we shall address can be stated as
follows.

\medskip
\noindent \textbf{Question\label{question}:} \label{question}\emph{
  Given $\theta_0\in \mathbb{R}^d$ and $\matr K$ a symmetric $d\times d$ matrix,
  let $\theta \sim P_2(\cdot \, | \, \mathbf{x})$ and consider
  \begin{equation}\label{eq:thetastargenintro} 
    \theta^\star =  \matr{K}\, ( \theta -\theta_0).
  \end{equation}
  Under what conditions on the prior, on the model, on the observation
  $\mathbf{x}$, and on the standardisation ($\theta_0$, $\matr{K}$) is
  the law of $\theta^{\star}$ close to 
  the standard normal
  distribution?  \medskip}

We will answer this question both in total variation distance and in
Wasserstein (a.k.a.\ Kantorovitch) distance.  In dimension 1,  our bounds are
sharp (sometimes even equalities); the following example, detailed in Section \ref{ex:weibull}, is
illustrative of the conclusions we obtain.
\begin{example}[Weibull data]
    \label{ex:weibullintro}
    Consider an observation sampled independently from Weibull data with density
    $ p_1(x \vert \ell,m) = ({m}/{\ell^m}) x^{m-1} \mathrm{exp} \{
    -(x/\ell)^m \} $ on the positive real line.  We suppose $m$ is
    known, and $\ell$ is the parameter of interest.  Fix, for
    $\tau_1>0$ and $\tau_2>-1$, a (conjugate) gamma prior on
    $\theta=\ell^{-m}$ of the form
    $\pi_0(\theta) \propto e^{-\eta(\theta)}$, with
    $\eta(\theta) = \tau_1\theta -\tau_2 \log \theta$ (for
    \(\tau_1=0\) and \(\tau_2=-1\) one obtains the improper Jeffreys
    prior for which our bounds also hold 
    when $n\geq 2$).  We consider the so-called MAP
    standardisation; 
    in the notation from \eqref{eq:thetastargenintro},
    we set $\theta_{\mathrm{MAP}}^{\star} = K(\theta - \theta_0)$ with
    $\theta_0 = \tilde \theta$ the posterior mode and
    $K = (\lambda''(\tilde \theta))^{1/2}$ where
    $ \lambda(\theta) = \tau_1\theta -(\tau_2 + n ) \log \theta$. Then our 
    bounds from Section \ref{ex:weibull} are 
    \begin{align*}  
     d_\mathrm{TV}( \theta^\star_{\mathrm{MAP}}, N)  
      \le  \sqrt{\frac{\pi}{8}}\frac{1}{  \sqrt{n+\tau_2} }, \quad d_\mathrm{Wass}( \theta^\star_{\mathrm{MAP}}, N)  
      =  \frac{1}{  \sqrt{n+\tau_2} }.
    \end{align*}
   This setting is also considered in \cite[Example
    5.2]{kasprzak2022good} where a Wasserstein {upper} bound is obtained;
    in contrast, we achieve equality.
    Our above bounds do not depend     on the sample nor on the parameter $\tau_1$; 
    {however}, letting $\tau_2$ be large artificially improves the quality of the Gaussian
    approximation.   In case we standardise around the MLE with a uniform prior, our bounds lead to a similar conclusion as above (in fact, the same bounds with $\tau_2$ set to 0), whereas taking a conjugate prior leads to (rather inelegant) 
    bounds on the total variation and Wasserstein distances which now   depend on all the parameters and the data $\mathbf{x}$  in an intricate way.
  \end{example}
  In dimension $d \ge 2$ our bounds remain explicitly computable.
  The different expressions  can also  be readily gauged in terms of the ``typical'' behavior of samples,  hereby  shedding new  light for instance on the dimensional dependence in the high-dimensional setting.   The next example is illustrative.

  \begin{example} [Categorical data]\label{ex:multiintro} In Section \ref{ex:multi}, we consider
    i.i.d.\ multinomial data with parameter 
    $p = (p_1,\ldots,p_k)^\intercal$ and individual likelihood
    $ p_1(x \vert p) = \prod_{j=1}^kp_j^{x_j}$ where
    $x = (x_1, \ldots,  x_k) \in \{0,1\}^k$ with $\sum_{j=1}^kx_j=1$, and $0<p_1,\ldots,p_k<1$ with $\sum_{j=1}^kp_j=1$.
    We fix a
    conjugate Dirichlet prior on $p$ of the form
$\pi_0(p) \propto (1 - \sum_{j=1}^{k-1} p_j)^{\tau_k}\prod_{j=1}^{k-1}
p_j^{\tau_j} $ for some parameters $\tau_1,\ldots,\tau_k>-1$. The natural parameter is
$\theta=(\theta_1,\ldots,\theta_{k-1})^\intercal$, where
$\theta_j = \log (p_j/p_k)$, $j=1,\ldots,k-1$, and we centre around the MAP as in Example \ref{ex:weibullintro}. 
Scaling with the {(unique symmetric)} square root of the Hessian
    evaluated at the MAP, we prove that there exist positive constants $\mathcal{C}_1$ and $\mathcal{C}_2$ such that
\begin{align}\label{eq:12mu}
 \quad d_{\mathrm{TV}}(\theta^{\star}_{\mathrm{MAP}}, N) \le   \frac{\mathcal{C}_1k}{\sqrt{\min \bm\pi}},  \quad d_{\mathrm{Wass}}(\theta^{\star}_{\mathrm{MAP}}, N) \le  \frac{\mathcal{C}_2k^{3/2}}{\sqrt{\min \bm\pi}}, 
\end{align}
where $\min\bm\pi=\min_{1\leq u\leq k}(n\bar{x}_u+\tau_u)$, and the Wasserstein distance bound is derived under the mild assumption that $\min\bm\pi > \max\{8,ck\}$ for some constant $c>0$. 
Moreover, we obtain the estimate $\mathcal{C}_1\leq 8.46$, and we also show that $\mathcal{C}_2\leq 7.40$ under the mild assumption that $\sqrt{\min\bm\pi}\geq7.40k/\sqrt{2}$.  
%
 Many bounds  in total
variation distance  are available for  BvMs with multinomial data, see e.g.\  \cite{boucheron2009bernstein, ghosal2000asymptotic, belloni2014posterior, ouimet2022multivariate, katsevich2023improved,katsevitchimprovedevenbetter}; as discussed in \cite{boucheron2009bernstein, katsevitchimprovedevenbetter}, bounds on different standardisations are not directly comparable.    We are not aware of a Wasserstein
bound (though \cite{katsevitchimprovedevenbetter} bounds the difference of means). 

In dimension 1 (i.e.\ when $k=2$, with Bernoulli data) our approach leads to  refined upper and lower  bounds 
which identify 
the role of the various sufficient statistics in the convergence, and also display 
the effect of the value of the ground truth success probability $p^\star$ on the quality of the approximations. For instance, we show the curious fact that  when $p^\star = (3 \pm \sqrt 3)/6$ then the upper and lower bounds on the Wasserstein distance coincide, leading to an explicit expression in this specific case.  Moreover, we also observe a phenomenon which, to the best of our knowledge, has not been reported previously:  in case of a ground truth success probability  $p^\star = 0.5$,  the choice of prior \emph{has an influence} on the asymptotic behavior of the posterior, in the sense that  the rate of convergence of the posterior towards the normal is no longer  $n^{-1/2}$ but rather $n^{-1}$, at least in case of a MAP centring with conjugate prior or an MLE centring with flat prior. When considering an MLE centring with conjugate prior, the rate of convergence only changes when $p^\star = 0.5$ for specific choices of the prior parameters. All our claims are supported by 
simulations, see Figures \ref{fig:binomial} and \ref{fig:binomial3};
Section \ref{example1} has 
more details.  
\end{example}

A complete statement of our bounds requires notations that would be
cumbersome to introduce already, and we refer the 
reader to the
forthcoming Theorem \ref{thm:abstract-results}. In a nutshell, our
bounds are exactly of the form \eqref{eq:kasp} with
$  \mathcal{U}_n(\mathbf{x})$ depending on several very natural
quantities: (i) a first order term on the log-likelihood measuring how
close $\theta_0$ is to a critical point (this term disappears when
$\theta_0$ is the MAP); (ii) a second order term
reflecting how close $\matr K$ is to the Hessian (this term also disappears when
$\theta_0$ is the MAP); (iii) a third order term measuring the
flatness of the log-likelihood around $\theta_0$; (iv)  remainder
terms (which are negligible in all our examples). In all these quantities the dependence on the data is through sufficient statistics only.  We could in principle work from here, as in most literature on
the topic, by identifying abstract conditions on the behavior of the log-likelihood which guarantee correct scaling with sample size and dimension; this would come at the cost of transparency and interpretability. We 
therefore choose rather to focus on specific concrete settings  for illustrating the power
of our approach. 

The method of proof rests on 
Stein's method of comparison of operators, see 
e.g.\ 
\cite{mrs21}. 
{T}his   method   provides bounds on integral probability metrics (see, e.g.,  \cite{zolotarev1984probability}),  
including the total variation and Wasserstein distances) between arbitrary probability measures under regularity conditions on their unnormalised densities.
Using 
results on the regularity of solutions to the standard normal Stein equation from  \cite{gms18}, we will show
via an application of Stein's method of comparison of operators that
that if $X$ has differentiable density $p^X$ and $N$ is standard normal then (see Lemma \ref{lemma:elem-steins-methWASSSBOUND})  
\begin{equation} 
  \label{eq:wassbounnnnintro}
  d_{\mathrm{Wass}}(X, N) \le \sup_{f \in \mathfrak{C}_1} \left| \mathbb{E} \left[
      \left\langle  \nabla \log p^X(X) +X, \, \nabla f(X) 
      \right\rangle \right]\right|
\end{equation}
with $\mathfrak{C}_1$ the collection of functions
    $f : \mathbb{R}^d \to \mathbb{R}$ which are twice differentiable
    with bounded second derivative and such that $\nabla f(x)$ 
    lies in the unit cube $\mathcal C_1$ in
    $\mathbb{R}^d$,     and  (see Lemma \ref{lemma:elem-steins-methTVBOUND})
\begin{equation}
  \label{eq:TVbounnnn1intro}
  d_{\mathrm{TV}}(X, N) \le C
  \sup_{f  \in \mathfrak{S}_1}  \left| \mathbb{E} \left[
      \left\langle  \nabla \log p^X(X) +X, \, \nabla f(X)
      \right\rangle \right]\right|
\end{equation}
 with $C$ an explicit constant and $\mathfrak{S}_1$  the collection of functions
  $f: \mathbb{R}^d \to \mathbb{R}$ which are twice differentiable with
  bounded first and second derivatives and such that
  $ \nabla f(x)$ 
  lies in the unit ball $\mathcal S_1$  in $\mathbb{R}^d$. 
  
  The Wasserstein distance bound (\ref{eq:wassbounnnnintro}), while not explicitly stated in the literature, follows directly from existing results. In contrast, the total variation distance bound (\ref{eq:TVbounnnn1intro}) is genuinely novel. It depends on new estimates for the solution of the standard normal Stein equation, which provide tighter control over \(\nabla f\) in the total variation framework (see Lemma \ref{sec:some-prop-solut-1}). These estimates allow the supremum in (\ref{eq:TVbounnnn1intro}) to be taken over a smaller set than what current results permit.
This shift from a unit-cube constraint for 
\(\mathfrak{C}_{1}\) in \eqref{eq:wassbounnnnintro} to a unit-ball constraint 
for \(\mathfrak{S}_{1}\) in \eqref{eq:TVbounnnn1intro} 
ensures an improved dimensional 
dependence for 
our bound on the total variation distance. For instance, the Wasserstein distance bound in~\eqref{eq:12mu} 
(which follows from \eqref{eq:wassbounnnnintro}) scales at a higher rate in the dimension than the total variation distance bound in~\eqref{eq:12mu}  (which follows from \eqref{eq:TVbounnnn1intro}). 
To our knowledge,
standard treatments of Stein's method 
skipped or simplified this aspect, thus missing the improved 
dependence on the dimension. This stronger control on \(\nabla f\) may also prove valuable for applications of Stein's method 
in other contexts, such as high-dimensional limit theorems or concentration-of-measure inequalities.

\medskip 
The paper is structured as follows. 
Section \ref{sec2} contains the core of
the paper: the detailed setting and assumptions in Section
\ref{sec:setup-1}, the Stein's method arguments and related general
bounds 
in Section~\ref{sec:elem-steins-meth}, and 
abstract bounds for  BvM-type settings 
in Section \ref{sec:back-bayes}. Section \ref{sec:applications}
contains 
various applications,
namely in the
univariate case: Bernoulli data (Section \ref{example1}), Poisson data
(Section \ref{poissonex}), normal with known mean data (Section
\ref{normalex}), Weibull data (Section \ref{ex:weibull}), and in the
multivariate case: multinomial data (Section \ref{ex:multi}), normal
data (Section \ref{ex:normal2}) and a linear regression setting
(Section
\ref{ex:linreg}). The Appendix 
contains further results
and the less illuminating technical proofs.

 \section{Main results} \label{sec2}

\subsection{Notations} 
\label{sec:notation}

{V}ectors $t \in \mathbb{R}^d$ are column vectors; the transpose
$t^{\intercal}$ is 
a row vector. Given two vectors $r, s$ in
$\mathbb{R}^d$ we write $r \otimes s$ as the matrix with components
$r_u s_v$, for $1 \le u, v \le d$; higher order tensor products are
defined similarly. For  two tensors $R, S$ of the same dimensions we
write
$\left\langle R, S \right\rangle = \sum_{i_1, \ldots, i_k}R_{i_1
  \ldots i_k} S_{i_1 \ldots i_k}$, the sum being taken over all
possible indices.  If $t \in \mathbb{R}^d$ we write
$\|t\|_2 = \sqrt{t^{\intercal}t} = \sqrt{\left\langle t, t
  \right\rangle}$ for its Euclidean norm.  If $\matr R$ is a symmetric {$d \times d$}
matrix we write
 $$
 \|\matr R\|_{\mathrm{diag}}= \sqrt{\max_{1 \le u \le d}|\matr
   R_{uu}|} \mbox{ and } \|\matr R\|_{\mathrm{sp}} = \sup \left\{ x'
   \matr R y \, | \, x, y \in \mathbb{R}^d
\mbox{ and }   \|x\|_2 = \|y\|_2 = 1 \right\}$$
 (the spectral norm). 
{T}he function  $\mathbb{I}_A(\cdot)$ denotes
the indicator function of the set $A  \in \mathcal{B}(\mathbb{R}^d)$,
where  $\mathcal{B}(\mathbb{R}^d)$ is  the
Borel $\sigma$-algebra on $\mathbb{R}^d$. 
Given a function $g : \mathbb{R}^d \to \mathbb{R}$ we write
$\|g\|_{\infty}$ for its essential supremum. If $g$ is 
$k$ times
differentiable at some point $t \in \mathbb{R}^d$ we denote by 
$\mathcal{D}^kg(t)$ the tensor with entries the $k$th partial
derivatives
$\partial^k_{i_1, \ldots, i_k} g(t) 
= \partial^k/(\partial t_{i_1}
\cdots \partial t_{i_k})g(t)
$ for $1 \le i_1, \ldots, i_k \le d$. We
let $\mathcal{D}^1 = \mathcal{D}$. If $g$ is differentiable,  its
gradient is the column vector
$\nabla g(t) = \mathcal D g(t) = (\partial_1 g(t), \ldots, \partial_d
g(t))^{\intercal}$. The Laplacian operator is denoted $\Delta$, and the Hessian of a twice differentiable function  
$g : \mathbb R^d \to \mathbb R$ 
is the symmetric matrix
$\matr{H}_g(t) = \mathcal{D}^2g(t) = \left( \partial_{u, v}^2g(t)
\right)_{1 \le u, v \le d}$. If $\matr{H}_g(t)$ 
is invertible, 
we
write $\matr{F}_g(t)$  for its inverse Hessian at $t$.  If in addition
$\matr{H}_g(t)$ is positive definite, it 
{has} 
a unique
invertible symmetric square root which we denote $\matr{K}_g(t)$; we
denote the inverse of this last matrix $\matr{J}_g(t)$. In other words,
\begin{equation}\label{eq:HKJ}
  \matr{H}_g(t) = \mathcal D^2 g(t), \, \matr{K}_g(t) =
  \left(  \matr{H}_g(t) \right)^{1/2}, \,
  \mathbf{F}_g(t) =  (\matr{H}_g(t))^{-1}   \mbox{ and }
  \matr{J}_g(t) = 
 \left(  \matr{H}_g(t) \right)^{-1/2}.
\end{equation}
We 
denote the entries of these matrices
$\matr H^g_{u, v}(t) ( = \partial^2_{u,v}g(t) )$,
$\matr K^g_{u,v}(t)$, $\matr F^g_{u,v}(t)$, and $\matr J^g_{u, v}(t)$
for $1 \le u, v \le d$, respectively.


We 
shall bound so-called 
\emph{Integral Probability Metrics
  ($\mathfrak{F}$-IPM)  which are of the form
$d_{\mathfrak{F}}(X, Y) = \sup_{\phi \in \mathfrak{F}} | \mathbb{E}
[\phi(X)] - \mathbb{E} [\phi(Y)]|$,  
where
$\mathfrak{F}$ is a class of test functions $\phi : \mathbb{R}^d \to \mathbb{R}$
such that $\phi \in L^1(X) \cap L^1(Y)$, see e.g.\ \cite{zolotarev1984probability}.}
Particular instances are 
${\mathfrak{F}}_{\mathrm{Wass}}=\{\phi:\mathbb{R}^d\rightarrow\mathbb{R}\,:\,
\phi \mbox{ is Lipschitz}, \|\phi\|_{\mathrm{Lip}}\leq1\}$, the class
of Lipschitz functions on $\mathbb{R}^d$ 
{with Lipschitz constant}
$\|\phi\|_{\mathrm{Lip}}:=\sup_{s\not=t}|\phi(s)-\phi(t)|/\|s-t\|_2
\le 1$,  
and
${\mathfrak{F}}_{\mathrm{TV}}=\{\mathbb{I}_B (\cdot) \,:\,B\in
\mathcal{B}(\mathbb{R}^d)\}$.
{Then}
\begin{align}
  d_{\mathrm{Wass}}(X, Y) = \sup_{\phi \in \mathfrak{F}_{\mathrm{Wass}}} | \mathbb{E} [\phi(X)]
  - \mathbb{E} [\phi(Y)]| \mbox{ and }    d_{\mathrm{TV}}(X, Y) = \sup_{\phi
  \in \mathfrak{F}_{\mathrm{TV}}} | \mathbb{E} [\phi(X)] 
  - \mathbb{E} [\phi(Y)]|,\label{eq:IPMss}
\end{align}
{denote} the (1-)Wasserstein  and total variation distances
between the laws of $X$ and $Y$.

{In this paper,} $N\sim N(0,I_d)$ {stands} for a standard normal random vector
in $\mathbb{R}^d$,  with density $\gamma_d(t)
$ on $\mathbb{R}^d$; 
$Y_1, Y_2, \ldots$
denote independent and identically distributed (i.i.d.)\ univariate uniform $U(0,1)$
random variables {which are independent of all other random elements.} 
For a random vector $X$ in $\mathbb{R}^d$, 
$L^p(X)$
{is} the collection of functions $\phi : \mathbb{R}^d \to \mathbb{R}$ such
that $|\phi^p(X)|$ has finite mean. 




\subsection{Setup}
\label{sec:setup-1}


{Let $n, r$ and $d$ be positive integers.}  Let
${\bf X} = (X_1, \ldots, X_n)$ be a vector of independent random
vectors,
each in $\mathbb{R}^r$, from regular $d$-parameter exponential
families with respective probability density functions
$ p_{1, i}(\cdot \, | \, \theta)$, $i = 1, \ldots, n$ (each with respect
to the same positive, $\sigma$-finite dominating measure $\mu$ defined
on the Borel sets of $\mathbb{R}^r$) given in the \emph{canonical
  form}
\begin{equation} \label{eq:likili}
   p_{1, i}(x \, | \, \theta) = \mathrm{exp}\big( \left\langle 
   h_i(x), \theta \right\rangle - \alpha_i(x) -  \beta_i(\theta)\big),
 \quad x \in \mathbb{R}^r,\quad  i = 1, \ldots, n,
 \end{equation}
 where $\theta=(\theta_1,\ldots,\theta_d)^\intercal$ is the common
 parameter of interest which we suppose to belong to some open set
 $ \Theta \subseteq \mathbb{R}^d$, and the functions
 $h_i : \mathbb R^r \to \mathbb R^d$ are chosen so that the
 normalising constants $\mathrm{exp}(-\beta_i(\theta))$ are finite for
 $i = 1, \ldots, n$. 
 For i.i.d.\ data, we shall write $\beta(\theta)=\beta_i(\theta)$ for all $i=1,\ldots,n$.
 Exponential family distributions in canonical form (\ref{eq:likili}) 
 are  analytically convenient and come without  
 loss 
 of generality as any exponential family
 distribution can be reduced to the canonical form via a
 re-parametrisation, see, e.g., 
 \cite{diaconis1979conjugate}; 
 a review of the results for  this paper is
 in Appendix \ref{sec:moments-expon-famil}. 

 Now let ${\bf x} = (x_1, \ldots, x_n) \in \mathbb{R}^{r \times n}$ be
 a fixed observation from model (\ref{eq:likili}) and choose a
 (possibly improper) prior $\pi_0$ on 
 ${\mathcal{B}}(\mathbb R^d)$, which, for notational convenience, we write as
 $ \pi_0(\theta) = \mathrm{exp}(-\eta(\theta)) $ for
 $\eta : \mathbb{R}^d \to \mathbb{R}$ a scalar function on which
 appropriate assumptions will be imposed below. The posterior
 distribution given $\mathbf{x}$ is then
\begin{align} 
  p_2 (\theta \,|\, \mathbf{x}) 
  & = \kappa(\mathbf{x})  \exp \big( n  \langle   
   \overline{h}(\mathbf{x}), \theta \rangle   - n \overline{\beta}(\theta)  -
    \eta(\theta) 
    \big),\label{eq:post} 
\end{align} 
with $\overline{h}(\mathbf{x}) = n^{-1} \sum_{i=1}^n h_{i}(x_i)$,
$\overline{\beta}(t) = n^{-1} \sum_{i=1}^n \beta_i(t) \in \mathbb{R}$
and $\kappa(\mathbf{x})$ the corresponding normalising constant. A key function in this paper is 
\begin{equation}
  \label{eq:lmbdaa}
  \lambda(\theta) = n \overline{\beta}(\theta)  +
\eta(\theta). 
\end{equation}

\medskip

Finally, let $\theta_0 =(\theta^0_1, \ldots, \theta^0_d) \in \Theta$
and $\matr{K}$ be a symmetric, invertible and positive definite
$d \times d$ matrix.  In line with the conventions from equation
\eqref{eq:HKJ} in Section~\ref{sec:notation}, we set
$\matr{J} = \matr{K}^{-1}$, $\matr{F} = \matr{J}^2$ and
$\matr{H} = \matr{K}^2$. Given
$\theta \sim p_2( \cdot \, | \, \mathbf{x}) $ we define
\begin{equation}\label{eq:thetastargen} 
    \theta^\star =  \matr{K}\, ( \theta -\theta_0).
  \end{equation}
 Then
  $\theta^{\star} ( = \theta^\star(\mathbf{x}, \matr K, \theta_0))$
  has density on 
$
\Theta^{\star} = \left\{ \matr{K}(\theta- \theta_0) \, | \,
  \theta \in \Theta\right\}$ given by
\begin{align}
  p^{\star}(t \, | \, \mathbf{x}) & = \mathrm{det}(\matr{J})
  p_2(\theta_0 +
  \matr{J}\, t \,
  | \, \mathbf{x})
   = \kappa(\mathbf{x}) \mathrm{det}(\matr{J})
     \exp \big( n\big\langle 
  \overline{h}(\mathbf{x}), \theta_0 + \matr{J}\, t \big\rangle 
         - \lambda \big(\theta_0 +
    \matr{J}\, t\big)    \big)  \label{eq:pstarlambda}.
\end{align}

\medskip

We can now reformulate the question from the introduction in a more
precise manner.

\medskip

\noindent \textbf{Question\label{question}:} \label{question}\emph{ Given
$\theta \sim p_2(\cdot \, | \, \mathbf{x})$ from model (\ref{eq:post})
and a standardisation as (\ref{eq:thetastargen}), under what
conditions on the prior $\eta$, on the model
${\beta_1, \ldots, \beta_n}$, on $\mathbf{x}$, and on the
standardisation ($\theta_0$, $\matr{K}$) is the law of
$\theta^{\star}$ close to $N(0,I_d)$
?
\medskip}

Our standing assumptions 
on
the model (\ref{eq:likili}) and the corresponding likelihood
(\ref{eq:post}) are as follows.

\medskip

\noindent {\label{assumtionA}{\bf  Assumption A}.}

\begin{enumerate}
  \setcounter{enumi}{-1}
\item\label{item:A0} 
{T}he maps
  $x  \mapsto p_{1, i} (x \,| \, \theta), 1 \le i \le  n,$ 
  have the same
  support $\mathcal{S}\subseteq\mathbb{R}^r$ which does not depend on
  $\theta$.
\item \label{item:A1}For all $1 \le i \le n$ the parameter $\theta$ is identifiable, in the sense that
  the maps $\theta \mapsto p_{1, i} (x \,| \, \theta)$ are 
  one-to-one for all $x \in \mathcal{S}$. The parameter space $\Theta$ is {open} and convex.

\item \label{item:A2}The posterior density
  $\theta \mapsto p_2( \theta \,|\, {\bf{x}})$ is {positive} and
  differentiable throughout $\Theta$, and for all
$f : \mathbb{R}^d \to \mathbb{R}$ with uniformly bounded first and second derivatives, 
\begin{equation}
  \label{eq:6}
  \int_{\Theta}
 \mathrm{div}( \nabla f( \theta )      p_2(\theta \, | \,
  \mathbf{x}) ) \,  d \theta = 0. 
\end{equation}
 
       \item \label{item:A3}${\lambda}$ is twice
        differentiable on $\Theta$ and the Hessian
         $\matr{H}_{\lambda}(\theta)$
  is symmetric
positive definite for all
   $\theta \in \Theta$.

\end{enumerate}

\begin{remark}
  Assumptions A\ref{item:A0} and A\ref{item:A1} are generic. Assumption A\ref{item:A2}, which
  depends on $\mathbf{x}$, ensures that all forthcoming integration by
  parts can be carried out without boundary terms. We could work under
  a weaker formulation, at the cost of additional boundary terms; see
  the forthcoming proof of Lemma~\ref{lma:steinopissteinop}. Finally,
  Assumption A\ref{item:A3} is very natural but could be weakened, see Remark
  \ref{rk:aboutheassum}.
\end{remark}

Although we will keep the statements of the results as general as
possible in terms of $\theta_0$ and $\matr K$, in all the examples
worked out in Section \ref{sec:applications} we will consider one of
the following two situations.

\begin{example}[Standardisation around the MAP]
  \label{sec:bounds-mapINIT}
 
  In (\ref{eq:thetastargen}) we fix
  $ \theta_0 = \tilde \theta = \tilde
  \theta_{\mathrm{MAP}}(\mathbf{x})$ the posterior mode which we
  suppose to be uniquely defined through the equation
  $
    \nabla \lambda(\tilde \theta) = n \bar h(\mathbf{x}). 
$
  We also take
  $\widetilde{\matr{K}}= \matr{K}_{\lambda}(\tilde \theta)$ as the
  unique symmetric square root of $\matr{H}_{\lambda}(\tilde \theta)$, and
  consider the normal approximation of
\begin{equation}\label{eq:thetastar} 
  \theta^\star_{\mathrm{MAP}} = \widetilde{\matr{K}}\big( \theta -\tilde \theta\big),
\end{equation} 
which is distributed according to the standardised posterior
distribution centred at the posterior mode.  
\end{example}

\begin{example}[Standardisation around the MLE]
  \label{sec:bounds-mleINIT}

  In (\ref{eq:thetastargen}) we fix
  $\theta_0 =\hat \theta  
  = \hat
  \theta_{\mathrm{MLE}}(\mathbf x) 
  $ the maximum likelihood estimate
  which is assumed 
  to be uniquely defined through the equation
  $
    \nabla \overline{\beta} (\hat \theta) =  \bar h(\mathbf{x}). 
 $
  We 
  take $\widehat{\matr{K}} = \sqrt n \matr{K}_{\overline{\beta}}(\hat \theta)$
  where $\matr{K}_{\overline{\beta}}(\hat \theta)$ is the unique symmetric
  square root of $\matr{H}_{\overline{\beta}}(\hat \theta)$, the Hessian of
  $\overline{\beta}$ evaluated at $\hat \theta$.  Thus we consider the
  normal approximation of
\begin{equation}
\label{eq:thetastar2} 
    \theta^\star_{\mathrm{MLE}}  = \widehat{\matr{K}}\big( \theta -\hat
    \theta\big), 
  \end{equation}
  which is distributed according to the standardised posterior
  distribution centred at the MLE.
\end{example}

Our approach to tackling our \hyperref[question]{Question} rests on two basic ingredients,
namely (i) general bounds on the total variation and Wasserstein distances
that are obtained from Stein's method of comparison of operators (see
Section \ref{sec:elem-steins-meth}) and (ii) Taylor expansions of the
log-likelihood ratio
$ \log(p^{\star} (\cdot \, | \, \mathbf{x})/\gamma_d(\cdot))$ (see
Lemma~\ref{lma:3rdordertayl} and  Appendix~\ref{sec4.2}).


\subsection{Stein's method of comparison of operators}
\label{sec:elem-steins-meth}

To sketch how we use Stein's method for Gaussian approximation via comparison of operators, first for a random vector $X$ on
$\mathbb{R}^d$ with differentiable density $p^X$ with support
$\mathcal{S}_X$, we introduce the Stein density operator (or \emph{Langevin-Stein operator}, see \cite{anastasiou2021stein}), 
for
$x \in \mathcal{S}_X$ as
\begin{equation} \label{steindensity} T_{X}f(x) =
  \frac{\mathrm{div}(p^X(x) \nabla f(x))}{p^X(x)} = \Delta f(x) +
  \langle \nabla \log p^X(x) , \nabla f(x)\rangle. 
 \end{equation} 
 Then 
 $ \mathbb{E}[ T_X f (X)] 
 = \int_{\mathcal{S}_X} \mathrm{div}(p^X(x)
 \nabla f(x)) dx
 = 0 $ for a large class of functions
 $f: \mathbb{R}^d \to \mathbb{R}$.
 We extend $T_X$ to act on $\mathbb{R}^d$ by
 setting it to 0 outside of $\mathcal{S}_X$.  
 The class of functions such that $\mathbb{E}[ T_X f (X)] =0$ is called
 the \emph{Stein class} for $p^X$ (or, equivalently, for $X$); we
 denote it $\mathcal{F}(X)$.  
 For $N\sim N(0,I_d)$ 
 its density $\gamma_d$ 
 satisfies $\nabla \log \gamma_d(x) = - x$ for all $x \in \mathbb{R}^d$, and $N$ hence has  density 
 Stein operator \eqref{steindensity} 
\begin{equation}\label{mvnop}T_{N}f(x) = \Delta f(x) - \langle  x ,
  \nabla f(x)\rangle \mbox{ on } \mathbb{R}^d;
\end{equation}
the corresponding Stein class $\mathcal{F}(N)$ 
contains for instance all twice differentiable functions with bounded
gradient.  This operator was first used in \cite{barbour90,gotze} and 
t
has been 
well studied in the literature. Of
particular importance are the 
PDEs of the form 
\begin{equation}\label{mvnsteineqn}
 T_Nf_{\phi}(x) 
 =    \Delta f_{\phi}(x) - \langle x, \nabla
    f_{\phi}(x)\rangle
    =\phi(x)-\mathbb{E}[\phi(N)];
  \end{equation}
with $\phi: \mathbb{R}^d \to \mathbb{R}$ belonging to $L^1(N)$. These
equations are   called the ($\phi$-)\emph{standard normal Stein
  equations} and allow to derive the  following (essentially
well-known) lemma.

  \begin{lemma}[Bound on differences of
    expectations]\label{lemma:elem-steins-methdiffexp}
    Let $N \sim \gamma_d$ and $X \sim p^X$ be as above.  Let
    $\phi : \mathbb{R}^d \to \mathbb{R}$ be such that there exists a
    twice differentiable solution $f_{\phi}$ to the Stein equation
    \eqref{mvnsteineqn} which, moreover, satisfies
    $f_{\phi} \in \mathcal{F}(X)$. Then
  \begin{align}
    \left| \mathbb{E}[\phi(X)] - \mathbb{E}[\phi(N)] \right|
    & \le    \left| \mathbb{E} \left[
      \left\langle  \nabla \log p^{X}(X) +X, \,  \nabla f_{\phi}(X) 
      \right\rangle \right]\right|.   \label{eq:diffexpboun} 
  \end{align}
\end{lemma}

\begin{proof}
  By construction,
  $(T_N - T_X )f(x) = - \langle \nabla \log p^X(x) +x, \, \nabla f(x)
  \rangle$ for all $x$ at which $f$ is differentiable. The stated
  assumptions are tailored to ensure \eqref{eq:diffexpboun} follows
  from \eqref{steindensity}, \eqref{mvnop} and \eqref{mvnsteineqn}
  along with the fact that $\mathbb{E} [T_X f_{\phi} (X) ] = 0$
  because $f_{\phi} \in \mathcal{F}(X)$.
\end{proof}

Recalling the integral probability metrics introduced at the end of Section
\ref{sec:notation},
one sees how the right-hand-side of \eqref{eq:diffexpboun}
provides bounds on IPMs, at least formally.  
To apply
Lemma \ref{lemma:elem-steins-methdiffexp} in practice it is crucial to
have good information on the behavior of the gradients of solutions to
(\ref{mvnsteineqn}). This is a well-studied problem. In dimension
$d = 1$ we refer e.g.\ to \cite[Example 2.31]{es22} where we can read the
following (here we correct  a typo in the last line of that example).

\begin{lemma}[Bounds in dimension $d=1$] \label{lma:one2d} For
  $\phi \in L^1( \gamma_1)$ let $f_{\phi} : \mathbb{R} \to \mathbb{R}$
  be a differentiable function with first derivative
  $f_{\phi}'(x) = e^{x^2/2} \int_{-\infty}^x (\phi(u) -
  \mathbb{E}[\phi(N)]) e^{-u^2/2}\, du$. Then $f_{\phi}$ is a twice
  differentiable solution to the Stein equation
  $f''(x) - x f'(x) = \phi(x) - \mathbb{E}[\phi(N)]$. Moreover,
  \begin{itemize}
  \item if $\phi \in \mathcal{F}_{\mathrm{Wass}}$ (Lipschitz with Lipschitz constant 1) then
    $\|f_{\phi}'\|_{\infty} \le 1$ and     $\|f_{\phi}''\|_{\infty} \le \sqrt{2/\pi}$; 
  \item if $\phi \in \mathcal{F}_{\mathrm{TV}}$ (indicators of  Borel sets) then
    $\|f_{\phi}'\|_{\infty} \le  \sqrt{\pi/8}$ and     $\|f_{\phi}''\|_{\infty} \le  2$.
  \end{itemize}
  \end{lemma}
  The case $d \ge 2$ is more complicated. We start with a regularity
  result from \cite[Proposition 2.1]{gms18} which will lead to
  Wasserstein bounds.
\begin{lemma}[Lipschitz test functions] \label{lema:reguwass} For
  $0<t<1$ and $x \in \mathbb{R}^d$ we set
  $ N_{t, x} = \sqrt t x + \sqrt{1-t} N.$ If $\phi$ is differentiable
  on $\mathbb{R}^d$ then the function $f_{\phi}$ defined by
  \begin{equation}
    \label{eq:solsteineq}
    f_{\phi}(x) = - \int_0^1 \frac{1}{2t} (\mathbb{E} [\phi
( N_{t, x} ) ] - \mathbb{E}[\phi(N)] )\, dt
\end{equation}
is a twice differentiable solution to (\ref{mvnsteineqn}). Let
$\mathcal{C}_1$ be the unit cube in $\mathbb{R}^d$.  If furthermore
$\phi \in \mathfrak{F}_{\mathrm{Wass}}$ 
then
  \begin{equation}\label{fhbound}
    \nabla f_{\phi} (x) \in \mathcal{C}_1  \mbox{ for all }  x \in
    \mathbb{R}^d \mbox{  and } \sup_{x \in \mathbb{R}^d} \| \nabla^2 f_\phi(x) \|_{\mathrm{H.S.}}
    \le \sqrt d. 
  \end{equation}
  \end{lemma}
  \begin{proof}
    The first statement comes from \cite[Proposition 2.1]{gms18}
    where it is shown to hold under the weaker condition that $\phi$
    is $\alpha$-H\"older. If $\phi$ is 
    Lipschitz
    then 
   the first claim in  (\ref{fhbound}) follows from  
    \cite{gr96} and \cite[Proposition
    2.1]{gms18}.  
    The second statement in (\ref{fhbound}) follows from \cite[Lemma
    3]{chatterjee2007multivariate}.
  \end{proof}

Using Lemma \ref{lema:reguwass} (resp., the first statement in Lemma
\ref{lma:one2d})  in the bound (\ref{eq:diffexpboun}) from Lemma
\ref{lemma:elem-steins-methdiffexp} we immediately reap a
  Wasserstein-1 bound on multivariate (resp., univariate) normal approximation.

  \begin{lemma}[Wasserstein Stein
    bound]\label{lemma:elem-steins-methWASSSBOUND}
    Let $\mathfrak{C}_1$ 
    be as in \eqref{eq:wassbounnnnintro} 
    Let $X$ be a random variable on $\mathbb{R}^d$ with
    density $p^X$ such that $\mathfrak C_1 \subseteq
    \mathcal{F}(X)$. Then
\begin{equation}
  \label{eq:wassbounnnn}
  d_{\mathrm{Wass}}(X, N) \le \sup_{f \in \mathfrak{C}_1} \left| \mathbb{E} \left[
      \left\langle  \nabla \log p^X(X) +X, \, \nabla f(X) 
      \right\rangle \right]\right|.
\end{equation}

\end{lemma}

Our next two results concern total variation distance. The one-dimensional case is covered by Lemma \ref{lma:one2d}. In dimension
$d \ge 2$, we note that solutions to (\ref{mvnsteineqn}) being only
valid for points $x\in\mathbb{R}^d$ where $\phi(x)$ is locally
H\"older continuous (see \cite[p.\ 726]{gotze}), we cannot directly
work with the Stein equation (\ref{mvnsteineqn}) with indicator test
functions $\phi \in \mathfrak{F}_{\mathrm TV}$ of the form
$\phi(\cdot) = \mathbb I_B (\cdot)$ with
$B\subset\mathcal{B}(\mathbb{R}^d)$.  
Instead we apply 
the smoothing approach used in \cite[Lemma 2.1]{gotze}.
\begin{lemma}
  \label{sec:some-prop-solut-1}
For  $0<t<1$ and  $\phi \in\mathfrak{F}_{\mathrm{TV}}$, let
$\phi_t(x) = \mathbb{E} [ \phi(\sqrt{t} N +
  \sqrt{1-t}x)].$
Then $\phi_t$ is differentiable
  on $\mathbb{R}^d$ and $ \|\partial_j \phi_t\|_{\infty} < \infty$ for all
  $j = 1, \ldots, d$.  
  The functions
\begin{equation}
\label{solnpsi}
    \tilde f_{\phi, t}(x):=-\int_t^1\frac{1}{2(1-s)}(\phi_s(x)-\mathbb{E}[\phi(N)])\,ds
\end{equation}
are well-defined {twice differentiable}  solutions of the
$\phi_t$-Stein equation (\ref{mvnsteineqn}) {with bounded
  first and second derivatives. Moreover, letting $\mathcal{B}(\delta)$ be
  the centred ball in $\mathbb{R}^d$} with radius $\delta$, we have 
\begin{equation}
\label{psibound}
 \nabla  \tilde f_{\phi, t}(x)  \in \mathcal{B}({{\sqrt{\pi/2}}})  \mbox{ for all } x \in \mathbb{R}^d
 \mbox{ and all } t \in (0, 1).
\end{equation} 
\end{lemma}
\begin{proof}
A change of variables gives 
\begin{align}
  \phi_t(x)  = \int_{\mathbb{R}^d} \phi(t^\frac12 
  y + ({1-t})^\frac12 x)
  \gamma_d( y)\, d y  = t^{-d/2}\int_{\mathbb{R}^d} \phi(u)
  \gamma_d((u -({1-t})^\frac12 x)/ t^\frac12 )\, du.   \label{eq:change}
\end{align}
Since $\phi$ is bounded and $\gamma_d$ is integrable, we can exchange
derivatives with respect to $x_j$ and integrals, recall  that  $\partial_j \gamma_d(x) = -x_j \gamma_d(x)$, and change variables back  to get
\begin{align*}
  \partial_j \phi_t(x)
  & = \left( \frac{1-t}{t}\right)^\frac12 \mathbb{E} \left[N_j \phi(t^\frac12 N +
    ({1-t})^\frac12 x)   \right] .
\end{align*}
  Since
$\phi$ is bounded 
it follows that $| \partial_j \phi_t(x)|<\infty$ for all $x\in\mathbb{R}^d$.
This proves the first part of the claim. 

The fact that
$ \tilde f_{\phi, t}$ is a well-defined twice differentiable solution
of the $\phi_t$-Stein equation \eqref{mvnsteineqn} follows immediately
from \cite[Lemma 2.1]{gotze}.  
For (\ref{psibound}), 
a change of variables as in  \eqref{eq:change} yields
\begin{align*}
    \tilde f_{\phi, t}(w)&
    =-\frac{1}{2}\int_t^1\!\int_{\mathbb{R}^d}\frac{1}{s^{d/2}(1-s)}\left(\phi(u)-\mathbb{E}[\phi(N)]\right)\gamma_d\left(\frac{u-(1-s)^{1/2}w}{s^{1/2}}\right)\,du\,ds.  
\end{align*} 
Therefore, by dominated convergence, since $\phi$ is bounded, we
obtain, for $x\in\mathbb{R}^d$ 
and 
$\xi \in
\mathbb{R}^d$: 
\begin{align*}
&   | \langle \xi, \nabla  \tilde f_{\phi, t}(x)\rangle|   =  \bigg|   \sum_{j=1}^d
                 \xi_j\partial_{j} \tilde f_{\phi, t}(x) \bigg|\\
  &=\frac{1}{2}\bigg|\sum_{j=1}^d
                 \xi_j\int_t^1\!\int_{\mathbb{R}^d}\frac{s^{-(d+1)/2}}{(1-s)^{1/2}}(u-(1-s)^{1/2}x)_j\left(\phi(u)-\mathbb{E}[\phi(N)]\right)\gamma_{d}\left(\frac{u-(1-s)^{1/2}x}{s^{1/2}}\right)\,du\,ds\bigg|\\
    &=\frac{1}{2}\bigg|\int_t^1\!\int_{\mathbb{R}^d}\frac{1}{({s(1-s)}^{\frac12}}\sum_{j=1}^d
                 \xi_jy_j\left(\phi(s^{1/2}y+(1-s)^{1/2}x)-\mathbb{E}[\phi(N)]\right)\gamma_{d}(y)\,dy\,ds\bigg|\\
  &\leq\frac{1}{2}\|\phi-\mathbb{E}[\phi(N)]\|_\infty\int_0^1\!\int_{\mathbb{R}^d}\frac{1}{\sqrt{s(1-s)}}
    \bigg|\sum_{j=1}^d 
                 \xi_jy_j\bigg|\gamma_{d}(y)\,dy\,ds. 
\end{align*}
If $\phi \in \mathfrak{F}_{\mathrm{TV}}$ then
$\|\phi-\mathbb{E}[\phi(N)]\|_\infty \le 1$.  Direct evaluation gives
$\int_0^1(s(1-s))^{-1/2}\,ds=\pi$. Also, letting $N_1, \ldots, N_d$
denote i.i.d.\ standard normal random variables, we have
$$\int_{\mathbb{R}^d} 
\bigg|\sum_{j=1}^d \xi_jy_j\bigg|\gamma_{d}(y)\,dy = \mathbb{E}
\bigg[\bigg|\sum_{j=1}^d \xi_{j}N_j\bigg|\bigg] = 
\sqrt{\frac{2}{\pi}}\|\xi\|_2,
$$
because $\sum_{j=1}^d \xi_{j}N_j \sim \mathcal{N}(0, \sum_{j=1}^d
\xi_j^2)$. Hence $| \langle \xi, \nabla  \tilde f_{\phi, t}(x)\rangle| \le\sqrt{
  {\pi}/{2}} \| \xi\|_2$  for all $t$, 
  $x$ and 
  $\xi$.
\end{proof}
\begin{remark}
  The statement \eqref{psibound} is independent of $t$; this
  property 
  is not enjoyed by higher order derivatives
  of $\tilde f_{\phi, t}$ (see \cite[Lemma 2.1]{gotze}). 
While the regularity of the solution of multivariate normal Stein equation is well studied, 
the fact that the gradient of the solution of the $\phi_t$-Stein equation belongs to a centred ball in $\mathbb{R}^d$, rather than just a centred cube in $\mathbb{R}^d$, appears to have been overlooked. 
This subtle difference 
will be crucial for obtaining 
bounds with an improved dimensional dependence in the total variation distance.   
\end{remark}

Similarly as in the Wasserstein case, 
Lemma
\ref{sec:some-prop-solut-1} (resp., the second statement in Lemma
\ref{lma:one2d}) in \eqref{eq:diffexpboun} leads to the following
total variation distance bound on 
normal
approximation.
\begin{lemma}[Global total variation
  bound]\label{lemma:elem-steins-methTVBOUND}
  Let $\mathfrak{S}_1$ be as in \eqref{eq:TVbounnnn1intro}. 
  Let $X \in \mathbb{R}^d$ be a random
  variable with density $p^X$ such that
  $\mathfrak S_1 \subseteq \mathcal{F}(X)$. Then
  \begin{equation}
  \label{eq:TVbounnnn1}
  d_{\mathrm{TV}}(X, N) \le C_d
  \sup_{f  \in \mathfrak{S}_1}  \left| \mathbb{E} \left[
      \left\langle  \nabla \log p^X(X) +X, \, \nabla f(X)
      \right\rangle \right]\right|
\end{equation}
where $C_1 = \sqrt{\pi/8}$ and $C_d = 
\sqrt{{\pi}/{2}}$ for all $d \ge 2$. 
\end{lemma}

\begin{proof}
  The claim in dimension $d=1$ is known (see e.g.\ \cite{es22}). For $d\geq2$, using
  \cite[(4.6) and (4.9)]{bhattacharya2010exposition} we get that for
  all $0<t<1$ there exists absolute constants $A_d\geq1$ and $B_d\geq0$ such that
  $ d_{\mathrm{TV}}(X,N)\leq A_d\sup_{\phi \in \mathfrak
    F_{\mathrm{TV}}}|\mathbb{E}[\phi_t(X)]-\mathbb{E}[\phi_t(N)]|+
  {B_d\sqrt{t}}/({1-t})$ with $\phi_t$ as in Lemma
  \ref{sec:some-prop-solut-1}. We then apply the same approach as in
  Lemma \ref{lemma:elem-steins-methWASSSBOUND} on the supremum, with
  the added fact that since the bound on the difference of
  expectations does not ultimately depend on $t$ we can take $t=0$ (which on examining the argument of \cite[Section 4]{bhattacharya2010exposition} allows us to also let $A_d=1$) to
  get the conclusion.
\end{proof}

Sharper bounds than in Lemma \ref{lemma:elem-steins-methTVBOUND} can often be obtained via the following refinement.

\begin{lemma}[Local total variation bound]\label{lemma:elem-steins-methTVBOUNDcondi}
  Let $\mathfrak S_1$  and $C_d$, $d \ge 1$, be as in Lemma
  \ref{lemma:elem-steins-methTVBOUND}.  Let
  $\mathcal{X} \in \mathcal{B}(\mathbb{R}^d)$ be compact with
  piecewise smooth boundary $\delta(\mathcal{X})$.  
  Let $X$
  have 
  density $p^X$ with respect to the Lebesgue measure on
  $\mathbb{R}^d$ which is positive and differentiable on
  $\mathcal{X}$. Finally, suppose that
  $\mathfrak S_1 \subseteq \mathcal{F}(X)$.  Then
\begin{align}
  \label{eq:TVbounnnn2}
  d_{\mathrm{TV}}(X, N) & \le C_d
                          \sup_{f\in \mathfrak{S}_1} \left| \mathbb{E}
                          \left[
                          \left\langle  \nabla \log p^X(X) +X, \,
                          \nabla f(X) 
                          \right\rangle \mathbb{I}_{\mathcal{X}}(X)\right]\right|+
                                                  \mathbb{P}[X \notin \mathcal{X}] +
                          r(\mathcal{X}),  
\end{align}
with
\begin{equation}\label{eq:errortermtvmz}
  r(\mathcal{X})  =  C_d \mathrm{A}_{d}(\delta(\mathcal{X})) \sup \{
    p^X(s ) \, | \, s \in
\delta(\mathcal{X})\}
  , 
                    \end{equation}
                    where 
                    $ \mathrm{A}_{1} (\delta(\mathcal{X})) = 1$ and $
                    \mathrm{A}_{d}(\delta(\mathcal{X}))$ 
                    is the surface area of  $\delta(\mathcal{X}) \subset \mathbb{R}^d$ for all
                    $d \ge 2$.
\end{lemma}

\begin{proof}
  First let $B \in \mathcal{B}(\mathbb{R}^d)$.  By using the law of
  total probability, we have
\begin{align*}
  |\mathbb{P}[X\in
  B]-\mathbb{P}[N\in B]| 
  &=|\mathbb{P}[X\in B\,|\,  X \in \mathcal{X}]
    \, \mathbb{P}[X
    \in\mathcal{X}]+\mathbb{P}[X\in B \,|\ X
    \not\in\mathcal{X}] \, \mathbb{P}[X
    \not\in\mathcal{X}]\\ 
  &\quad-\mathbb{P}[N\in B\,|\, X \in \mathcal{X}] \,
    \mathbb{P}[X \in \mathcal{X}]-\mathbb{P}[N\in
    B\,|\,X \not\in\mathcal{X}] \, \mathbb{P}[X
    \not\in \mathcal{X}]|\\ 
  &\leq|
    |\mathbb{P}[X_{\mathcal{X}}\in B]-\mathbb{P}[N\in B]|\mathbb{P}[X
    \in\mathcal{X}]+\mathbb{P}[X \not\in \mathcal{X}]
\end{align*}
with  $X_{\mathcal{X}}$ the random variable with
distribution $\mathcal{L}(X \, | \, X \in
\mathcal{X})$ and $N$ independent of $X$. Hence, 
 $ d_{\mathrm{TV}}(X, N)  \le   d_{\mathrm{TV}}(
  X_{\mathcal{X}}, N) \mathbb{P}[X \in\mathcal{X}]+ 
  \mathbb{P}[X \not\in \mathcal{X}]. $ 
Again we use the inequality   $ d_{\mathrm{TV}}(X_{\mathcal{X}},N)\leq\sup_{\phi \in \mathfrak
    F_{\mathrm{TV}}}|\mathbb{E}[\phi_t(X_{\mathcal{X}})]-\mathbb{E}[\phi_t(N)]|$.
  Some easy algebraic manipulations
  give, with $T_N$ as in \eqref{mvnop}, 
\begin{align*}
  &     \mathbb{E} \left[ \phi_t(X_{\mathcal{X}}) - \mathbb{E} [\phi_t(N)]
    \right]  \mathbb{P}[X \in \mathcal{X}]  = \int_{\mathcal{X}}  (
    T_N \tilde f_{\phi, t}(x) )  p^X(x) \,dx \\
  & = - \mathbb{E} \left[ \left\langle \nabla  \log p^X(X) + X, \nabla
    \tilde f_{\phi, t}(X) \right\rangle 
    \mathbb{I}_{\mathcal{X}}(X) \right] +  \int_{\mathcal{X}}
    \mathrm{div}( p^X(x) \nabla \tilde f_{\phi, t}(x))\, dx.
\end{align*}
  The first term above gives the first term in
\eqref{eq:TVbounnnn2}. By the divergence theorem we deduce
  \begin{align*}
  \int_{\mathcal{X}} \mathrm{div}(p^X(u )  \nabla
{\tilde f_{\phi, t}}(u))\, du  
    & =  \oint_{\delta(\mathcal{X})} p^X (s )  \left\langle  \nabla
      { \tilde f_{\phi, t}}(s),  n(s)\right\rangle  d \sigma(s) 
  \end{align*}
  with $\sigma$ the surface measure and $n$ the outward pointing unit
  normal vector.  Using the Cauchy-Schwarz inequality with the inequality $\|
  \nabla  \tilde f_{\phi, t} \|_2 \le \sqrt{\pi/2}$ in (\ref{psibound}) (and in dimension $d=1$ we have $\|f_\phi'\|_\infty\leq\sqrt{\pi/8}$, from Lemma \ref{lma:one2d}) we get
  \begin{align*}
\int_{\mathcal{X}} \mathrm{div}(p^X(u )  \nabla
{ \tilde f_{\phi, t}}(u))\, du  
 & \le 
 C_d
 \oint_{\delta(\mathcal{X})}
                      p^X (s)  \, d \sigma(s)
                      \le 
                      C_d\sup \left\{
                      p^X (s ) \, | \, s \in
\delta(\mathcal{X}) \right\}
                      \mathrm{A}_{d}({\delta(\mathcal{X})}).
  \end{align*}
  Adding the various pieces gives the claim.
\end{proof}

\begin{remark}\label{rk:XiSTH}
  At least formally, \eqref{eq:TVbounnnn1} follows from
  \eqref{eq:TVbounnnn2} by taking $\mathcal{X} = \Theta$ and employing
  the convention $r(\Theta) = 0$. This will allow us to simultaneously state both the
  local and global total variation bounds.
\end{remark} 



\subsection{Back to the Bayes} \label{sec:back-bayes}

Recall  the exponential family setup and  all the corresponding notations  from Section \ref{sec:setup-1}. 
{
We now  address our  \hyperref[question]{Question} by applying the Stein's method developed in  Section \ref{sec:elem-steins-meth} to  the specific choice 
$X =\theta^{\star} = \matr{K}(\theta - \theta_0)$ for   some $\theta$ distributed according to the posterior distribution in an exponential family as in \eqref{eq:post},  $\theta_0 \in \Theta$ some arbitrary point,  and $\matr{K}$ some symmetric, invertible and positive definite
$d \times d$ matrix. 
We begin by  noting the following.
\begin{lemma}\label{lma:steinopissteinop}
  Let  $p^{\star}$ be the  density  of $\theta^\star$ on  $\Theta^{\star}$ as given in
\eqref{eq:pstarlambda}.  Under  Assumption \hyperref[assumtionA]{A}, the Stein operator $T_{\theta^{\star}}$ for 
  $\theta^{\star}$, given for $t^{\star} \in \Theta^{\star}$ by
  $$T_{\theta^{\star}}f(t^{\star}) = \Delta f(t^{\star}) +\left\langle \nabla \log
    p^{\star}(t^{\star} \, | \, \mathbf{x}), \nabla f(t^{\star})
  \right\rangle$$ 
  and set to 0 outside of $\Theta^{\star}$, 
  satisfies
  $\mathbb E[ T_{\theta^\star} f(\theta^\star)] = 0$ for all twice
  differentiable functions $f : \mathbb{R}^d \to \mathbb{R}$ with
  bounded first and second derivatives.
 \end{lemma}

                            \begin{proof}
                              Let $f : \mathbb{R}^d \to \mathbb{R}$ be
                              a differentiable function.  We note 
                             (see e.g.\ \cite{mrs21}) that for
                              all $t^{\star} \in \Theta^{\star}$ the Stein operator
                              satisfies
$                                p^{\star}(t^{\star} \, | \, \mathbf{x})
   T_{\theta^{\star}}f(t^{\star}) = \mathrm{div}(
   \nabla f(t)  p^{\star}(t^{\star} \, | \,
  \mathbf{x})). 
$
  Hence
  \begin{equation*}
   \mathbb{E} \left[
     T_{\theta^{\star}}f(\theta^{\star})
   \right] = \int_{\Theta^{\star}} \mathrm{div}(
   \nabla f(t^{\star})  p^{\star}(t^{\star} \, | \,
    \mathbf{x}))\, dt^{\star}. 
   \end{equation*}
  The conclusion then follows directly
  after a change of variables in
    Assumption (A\ref{item:A2}).
   \end{proof}
 
   Moreover, the  following can be seen to hold (see 
    Appendix \ref{sec4.2}).
 \begin{lemma}[Third order Taylor
   expansion]\label{lma:3rdordertayl}
    Let $\Theta_0\subseteq \Theta$ be a
   measurable neighborhood of $\theta_0$
    and set 
   $\Theta_0^{\star} = \matr K (\Theta_0
     -\theta_0)$.    Let $\matr J = \matr K^{-1}$, 
     set 
   $\bar h$ as in \eqref{eq:post}   and  $\lambda$ as in  \eqref{eq:lmbdaa}.   If $\lambda$ is three
    times differentiable on $\Theta_0$, then
                              for all $t^{\star} \in \Theta_0^{\star}$
                              and all $\xi \in \mathbb{R}^d$ we have
\begin{align}
  \left\langle  \nabla \log p^{\star}(t^{\star} \, | \, \mathbf{x}) +t^{\star}, \, \xi
  \right\rangle 
  & = \left\langle \nabla \lambda(\theta_0) -  n \bar
    h(\mathbf{x}) , \,                                        
    \matr{J}   \xi \right\rangle \nonumber +                                    
    \left\langle                                                                        
    \matr{H}_{\lambda}(
    \theta_0) 
    - \matr{K}^2, \,   \left( \matr{J} \xi \right) 
    \otimes
\left(     \matr{J} t^{\star} \right)
    \right\rangle\\
  &    \quad     +  \left\langle \bm \Upsilon^3(\lambda, \theta_0; t^{\star}), \left( \matr{J} \xi
    \right)
    \otimes
    (\matr{J} t^{\star})^{\otimes
    2} 
    \right\rangle
                                                                     \label{eq:3dorder} 
\end{align}
where for $g: \mathbb{R}^d \to \mathbb{R}$ three times
  differentiable on $\Theta_0$  and $t^{\star} \in \Theta_0^{\star}$ we
  set
$$ \bm\Upsilon^3 (g, \theta_0; t^{\star}) = \mathbb{E} \left[ Y_1 \mathcal{D}^3
  g(\theta_0 + Y_1 Y_2 (\matr{J} t^{\star}) )\right]$$ with $Y_1, Y_2$
i.i.d.\ uniform random variables on $(0,1)$ and
$\matr J = \matr K^{-1}$.

\end{lemma}


In order to efficiently state our bounds we will need some final
notations.
\begin{definition}[Bespoke derivatives]\label{def:bespoke}
With all notations as above we set  
\begin{align}
  \label{eq:1}
  \mathfrak{D}_1(g, \theta_0, \matr K \, | \,  \Theta_0^{\star})  & =  
  \sum_{u=1}^d                     \left|                  \partial_ug(
    \theta_0
    ) - n \bar
    h_u(\mathbf{x})
  \right| \mathbb{P}[\theta^{\star} \in \Theta_0^{\star}], \\
  \label{eq:2}
  \mathfrak{D}_{2}(g,
  \theta_0, 
  \matr{K}\, | \,  \Theta_0^{\star})  &  = 
  \sum_{u=1}^d
    \sum_{v=1}^d
        \left|
   \partial^2_{u,
       v}g(\theta_0)
      - (\matr{K}^{2})_{u,
      v} \right|\mathbb{E}
     \left[ \left|
     ( \matr J
       \theta^{\star})_v \right| \mathbb{I}_{\Theta_0^{\star}}(\theta^{\star})\right], \\
  \label{eq:4}
  \mathfrak{D}_{3}(g, \theta_0, \matr K \, | \, \Theta_0^{\star})          &
      = \sum_{u=1}^d
    \sum_{v=1}^d \sum_{w=1}^d  \mathbb{E}\left[ \left|
                 \bm  \Upsilon^3_{u, v, w}(g, \theta_0;
      \theta^{\star})(\matr
      J \theta^{\star})_v
       (\matr J 
     \theta^{\star})_w\right|\mathbb{I}_{\Theta_0^{\star}}(\theta^{\star})\right]. 
\end{align}
Finally, we set 
  \begin{align}\label{eq:5}
    \Delta_{3}(g,
    \theta_0, \matr K \, | \,   \Theta_0^{\star})
    & = \sum_{j=1}^3  \mathfrak{D}_j(g, \theta_0,\matr K \, | \,   \Theta_0^{\star}).  \end{align}
 \end{definition} 
  \begin{remark} \label{rem:aee}
    In line with the comment from Remark \ref{rk:XiSTH} we will drop
    the notation ``$\, | \, \Theta_0^{\star}$'' in any of the above
    when $\Theta_0= \Theta$ (or, equivalently, when
    $\Theta_0^{\star} = \Theta^{\star}$).
  \end{remark}

With these quantities defined, we are ready to state our bounds.

\begin{theorem}[Third order bounds] \label{thm:abstract-results} Let
  $\theta$ be distributed according to (\ref{eq:post}) and suppose
  Assumption \hyperref[assumtionA]{A} is satisfied. Let
  $\theta^{\star}$ be defined as in (\ref{eq:thetastargen}) for
  $\matr K$ a square symmetric definite positive matrix and
  $\theta_0 \in \Theta$; define $\Theta^{\star}$ accordingly. Let
  $\matr J = \matr K^{-1}$ and $\matr F = \matr J^2$.  Set
  $C_1 = \sqrt{\pi/8}$ and $C_d =  \sqrt{\pi/2}$ for $d\ge 2$ (as
  in Lemma~\ref{lemma:elem-steins-methTVBOUND}).
  %
  Then
  the following statements  hold true:
\begin{itemize}
\item  If $\lambda$ is three times differentiable
  on $\Theta^{\star}$ then
\begin{align}
  & d_{\mathrm{Wass}}(\theta^{\star}, N) \le   \sqrt{d} \| \matr
  F\|_{\mathrm{diag}}  \Delta_{3}(\lambda,
  \theta_0, \matr K),   \label{eq:boundwass} \\
&     d_{\mathrm{TV}}(\theta^{\star}, N) \le   C_d \| \matr
  F\|_{\mathrm{diag}}  \Delta_{3}(\lambda,
    \theta_0, \matr K).   \label{eq:bounTV}
\end{align}
\item Let
  $\Theta_0^{\star}\subset \Theta^{\star} \subsetneq \mathbb{R}^d$ be
  a compact set with piecewise smooth boundary.  If $\lambda$ is three
  times differentiable on $\Theta_0^{\star}$ then
\begin{align}
  \label{eq:TVbounnnnthetastarcond}
  d_{\mathrm{TV}}(\theta^{\star}, N) & \le C_d\| \matr F\|_{\mathrm{diag}}
                                       \,  \Delta_{3}(\lambda, 
                                       \theta_0,
                                       \matr K \, | \,
                                      \Theta_0^{\star} )  +  \mathbb{P}[ 
                                       \theta^{\star}  \notin
                                      \Theta_0^{\star}] + 
                                       r_d^{\star}(\Theta_0^{\star})
\end{align} 
where 
$ r_d^{\star}(\Theta_0^{\star}) = C_d  A_d^{\star} p^{\star}_{\infty}$ with $p^{\star}_{\infty} = \sup \left\{ p^{\star} (\theta^{\star} \, | \,
    \mathbf{x}) \, | \, \theta^{\star} \in \delta(\Theta_0^{\star}\right\}$, 
 $A_1^{\star} = 1$ and 
  $A_d^{\star}$ is the surface area of $\delta(\Theta_0^{\star})$ if $d \ge
  2. $
\end{itemize}
\end{theorem}

 \begin{proof}
   Reading directly from Lemma \ref{lma:3rdordertayl} we see that 
   it suffices to tackle the quantities
   \begin{align*}& 
\left|      \left\langle  \nabla
     \lambda(\theta_0) - n \bar h(\mathbf{x}), 
                   \matr J\mathbb{E} \left[ \nabla f(\theta^{\star}) \right] 
     \right\rangle  \right|  +  \left|\left\langle   \matr
     H_{\lambda}(\theta_0) - \matr K^2,  \mathbb{E} \left[  (\matr J \nabla f(\theta^{\star})) \otimes
 (\matr J  \theta^{\star}) \right] 
     \right\rangle \right| \\
     & +   \left|\mathbb{E} \left[ \left\langle  \bm \Upsilon^3(\lambda, \theta_0; \theta) ,  (\matr J \nabla f(\theta^{\star})) \otimes
      (\matr J  \theta^{\star})^{\otimes 2}  
     \right\rangle  \right]  \right| 
   \end{align*}
   for $f$ belonging either to $\mathfrak{C}_1$ in order to apply
   Lemma \ref{lemma:elem-steins-methWASSSBOUND} and obtain the Wasserstein bound
   \eqref{eq:boundwass}, or to $\mathfrak S_1$ in order to apply
   Lemmas \ref{lemma:elem-steins-methTVBOUND} and
   \ref{lemma:elem-steins-methTVBOUNDcondi} on
   $\mathcal{X} = \mathcal{B}(\delta)$ and obtain
   the total variation bound \eqref{eq:TVbounnnnthetastarcond}.

   For 
   \eqref{eq:boundwass}, note that
   for all $\xi \in \mathcal{C}_1$ (the unit cube in $\mathbb{R}^d$)
   and all $1 \le u \le d$ we have
   \begin{align*}
    | (\matr J \xi)_u | \le  \left|\sum_{v=1}^d\matr J_{uv} \xi_v\right|  \le \sqrt{\sum_{v=1}^d\matr
     J_{uv}^2 \sum_{w=1}^d \xi_w^2}  \le \sqrt d \sqrt{
      \matr F_{uu}}  \le \sqrt d \| \matr F \|_{\mathrm{diag}}.
   \end{align*}
Bound \eqref{eq:boundwass} follows  from Lemma \ref{lemma:elem-steins-methWASSSBOUND}. 
For 
the bound 
\eqref{eq:bounTV} we apply Lemma
\ref{lemma:elem-steins-methTVBOUND}    along
with the fact that for all $\xi \in \mathcal{B}_1$ (the unit disk in
$\mathbb{R}^d$) and all $u = 1, \ldots, d$ we have
    \begin{align*}
    | (\matr J \xi)_u | \le       \left|  \sum_{v=1}^d \matr J_{uv} \xi_v  \right|  \le \sqrt{\sum_{v=1}^d \matr
      J_{vu}^2} =  \sqrt{
      \matr F_{uu}}  \le  \| \matr F \|_{\mathrm{diag}}.
    \end{align*} 
Finally, the bound 
 \eqref{eq:TVbounnnnthetastarcond} follows from the same arguments,
 through  Lemma 
\ref{lemma:elem-steins-methTVBOUNDcondi} on $\Theta_0^{\star}$. 
 \end{proof}

 \begin{remark}\label{rk:aboutrem} A simple
   choice of $\Theta_0^{\star}$ is the ball $\mathcal{B}(\delta)$ with
   radius $\delta$ whose surface area is given by 
   $ A_d^{\star} = {\delta^{d-1} {d}\pi^{d/2}}/{\Gamma(d/2+1)} \sim
\pi^{-1/2}(2e\pi)^{d/2}\delta^{d-1}/d^{(d-1)/2} $, see e.g.\
   \cite[Section 3]{tee2004surface}. 
   This quantity is negligible in $d$ as long as $\delta$
   does not grow too fast with the
   dimension. 
 \end{remark}
 
 \begin{remark}[About the assumptions]\label{rk:aboutheassum}
   Lemma \ref{lma:3rdordertayl} is a particular case of a general
   order Taylor expansion provided in Lemma \ref{lem:diffexp} in
   Appendix \ref{sec4.2}. Although in Theorem
   \ref{thm:abstract-results} we focus on third order expansions, we
   could 
   work out higher order bounds (under higher
   order differentiability assumptions) from \eqref{eq:kdorder} in 
   Appendix \ref{sec4.2}.  We could also dispense with the third order
   differentiability requirement and only take derivatives up to order
   2, by using \eqref{eq:2dorder}. Also, 
 we could
 easily dissociate the assumptions on $\overline{\beta}$ and
   $\eta$ (for instance imposing higher order regularity on
   $\overline{\beta}$ but only very weak differentiability on the
   prior $\eta$) and expand to different orders in $\beta$ and $\eta$. 
   We do not detail all the possible variants, for conciseness, 
   but provide the following univariate fourth order expansion, which we apply in Example \ref{example1} to obtain a faster convergence rate with respect to the sample size $n$ when the model takes a particular form.
   Let
\begin{align}
 \Delta_{4}(g,
    \theta_0, \matr K \, | \,   \Theta_0^{\star})
    & = \sum_{j=1}^2  \mathfrak{D}_j(g, \theta_0,\matr K \, | \,   \Theta_0^{\star})+\frac{1}{2}\big|g^{(3)}(\theta_0)\big|\mathbb{E}\left[(J\theta*)^2\mathbb{I}_{\Theta_0^{\star}}(\theta^{\star})\right]\nonumber\\
    &\quad+\mathbb{E}\left[\left|\mathbb{E}\left[Y_1Y_2^2g^{(4)}(\theta_0+Y_1Y_2Y_3J\theta^*)\right](J\theta^*)^3\right|\mathbb{I}_{\Theta_0^{\star}}(\theta^{\star})\right],\label{delta4}
\end{align}   
where $Y_1,Y_2,Y_3$ are i.i.d.\,
$U(0,1)$. 
Then, if $\lambda:\mathbb{R}\rightarrow\mathbb{R}$ is four times differentiable
  on $\Theta^{\star}$,
\begin{align}
  & d_{\mathrm{Wass}}(\theta^{\star}, N) \le   \|\matr{F} \|_{\mathrm{diag}} \Delta_{4}(\lambda,
  \theta_0, \matr K),   
  \quad    d_{\mathrm{TV}}(\theta^{\star}, N) \le  \sqrt{\frac{\pi}{8}} \|\matr{F} \|_{\mathrm{diag}} \Delta_{4}(\lambda,
    \theta_0, \matr K).   \label{eq:bounTV2}
\end{align}
 \end{remark}

 The 
 bounds depend on
 $\|\matr{F} \|_{\mathrm{diag}}$ 
 and
 $\Delta_3(\lambda, \theta_0, \matr K \, | \, \Theta_0^{\star}) = \sum_{j=1}^3 \mathfrak
 D_j(\lambda, \theta_0, \matr K \, | \, \Theta_0^{\star})$.  The first
 two summands in this last term simply reflect the influence of the
 centring and scaling parameters from (\ref{eq:thetastargen}) on the
 Gaussian proximity: $\mathfrak D_1(\lambda, \theta_0, \matr K \, | \, \Theta_0^{\star})$ indicates how
 close $\theta_0$ is to a critical point of
 $p_2(\cdot \, | \, \mathbf{x})$ and
 $\mathfrak D_2(\lambda, \theta_0, \matr K \, | \, \Theta_0^{\star})$ indicates how close
 the scaling $\matr K$ is to the Fisher information. We thus have the
 following simplication under MAP standardisation.
 
\begin{example}[MAP standardisation, continued]
  \label{sec:bounds-map}
Using the notations from Example \ref{sec:bounds-mapINIT}, it follows that
\begin{align}\label{eq:thingscancel}
  \mathfrak{D}_1 (\lambda, \tilde \theta ,  \widetilde{\matr{K}} \, |
  \,  \Theta_0^{\star}) = \mathfrak{D}_{2} (\lambda,
  \tilde
  \theta
  , \widetilde{\matr{K}} \, |\,   \Theta_0^{\star}) = 0
\end{align}
for all $\Theta_0^{\star}$ (hence also for the unconditional
version). 
\end{example}
 
For general standardisation, 
as $\lambda = n \overline \beta + \eta$,  
linearity and the triangle inequality yield that
\begin{align*}
  \mathfrak D_1 (\lambda,  \theta_0, \matr K \, | \, \Theta_0^{\star}) & \le n
   \mathfrak D_1 (\overline{\beta},  \theta_0, \matr K \, | \,
                                                                        \Theta_0^{\star} ) +   \sum_{u=1}^d \left|
                                                        \partial_u
                                                        \eta(\theta_0)
                                                        \right| \mathbb{P} [\theta^{\star} \in
      \mathcal{B}(\delta)],  \\
   \mathfrak D_2 ( \lambda, \theta_0, \matr K \, | \, \Theta_0^{\star}) & \le n
   \mathfrak D_2 (\overline{\beta}, \theta_0,  \matr K \, | \,
                                                                          \Theta_0^{\star})
                                                                          +
                                                                          \sum_{u=1}^{d}
                                                                          \sum_{v=1}^d
                                                                          \left|
                                                                          \matr
                                                                          H^{\eta}_{u,
                                                                          v}(\theta_0)
                                                                        \right|\mathbb{E}
                                                         \left[ \left|
                                                                            ( \matr J
                                                         \theta^{\star})_u \right| \mathbb{I}_{\Theta_0^{\star}}(\theta^{\star})\right]
\end{align*}
and
\begin{align}
   \mathfrak{D}_{3}(\lambda, \theta_0 , \matr K\, | \,
 \Theta_0^{\star})  \le  n \mathfrak{D}_{3}(\bar \beta, \theta_0 , \matr K\, | \,
 \Theta_0^{\star}) + \mathfrak{D}_{3}(\eta, \theta_0 , \matr K\, | \,
 \Theta_0^{\star}).  \label{d3ineq}
\end{align}
The prior dependent terms vanish in all of the above under a flat
(uninformative) prior; under a Gaussian prior we have
                     \begin{equation*}
                        \sum_{w=1}^d
      \left|  \partial_w\eta( 
       \theta_0
      )  \right| = \sum_{w=1}^d | \theta^0_w| \mbox{ and }  \matr
     H_{\eta}(\theta_0) = I_d. 
 \end{equation*} 
 In the case of MLE standardisation some further simplifications occur,
                     as follows.

\begin{example}[MLE standardisation, continued]
  \label{sec:bounds-mle}
  Using the notations from Example \ref{sec:bounds-mleINIT}, it
  follows that
   $\| \widehat{\matr{F}}\|_{\mathrm{diag}} = n^{-1/2} \| \matr
    F_{\overline \beta}(\hat \theta)\|_{\mathrm{diag}}; $
  moreover,
  $ \mathfrak{D}_1 (\bar \beta, \hat \theta, \widehat{\matr{K}} \, | \,
 \Theta_0^{\star}) = \mathfrak{D}_2 (\bar \beta, \hat \theta, \widehat{\matr{K}} \, | \,
 \Theta_0^{\star}) = 0$ and 
      \begin{align*}
         \mathfrak D_1 (\lambda,  \hat \theta, \widehat{\matr{K}} \, | \,
        \Theta_0^{\star}) & \le   \sum_{u=1}^d \left|
                                                        \partial_u
                                                        \eta(\hat\theta)
                                                        \right| \mathbb{P} [\theta^{\star} \in
 \Theta_0^{\star}], \\
 \mathfrak {D}_{2} (\lambda,
     \hat\theta,
 \widehat{\matr{K}} \, | \,
                                                        \Theta_0^{\star})      
      & \le   \frac{1}{\sqrt n} \sum_{u=1}^{d}
                                                                          \sum_{v=1}^d
                                                                          \left|
                                                                          \matr
                                                                          H^{\eta}_{u,
                                                                          v}(\hat\theta)
                                                                        \right|\mathbb{E}
                                                         \left[ \left|
                                                                            (
    \matr J_{\overline \beta}(\hat \theta)
                                                         \theta^{\star})_u \right| \mathbb{I}_{\Theta_0^{\star}}(\theta^{\star})\right] 
    .  \end{align*}
  
\end{example}
 
\begin{example}[MAP standardisation, continued] 
As $\lambda=n\bar\beta+\eta$, if 
the various quantities in inequality (\ref{d3ineq}) exist, we have that $\mathfrak{D}_{3}(\lambda, \theta_0 , \matr K\, | \,
 \Theta_0^{\star})  \le  n \mathfrak{D}_{3}(\bar \beta, \theta_0 , \matr K\, | \,
 \Theta_0^{\star})\{1+r_n^1\}$, where $r_{n}^1$ is a remainder term of order $O(n^{-1})$. 
 As the posterior mode $\tilde\theta$ solves $\nabla\bar\beta(\theta)+n^{-1}\nabla\eta(\theta)=\bar{h}(\mathbf{x})$ whilst the MLE $\hat\theta$ solves $\nabla\bar\beta(\theta)=\bar{h}(\mathbf{x})$, 
 provided the functions $\bar\beta$ and $\eta$ possess suitability regularity, we have that $\tilde\theta=\hat\theta\{1+r_n^2\}$, where $r_{n}^2$ is a remainder term of order $o(1)$ as $n\rightarrow\infty$ (in all our examples of Section \ref{sec:applications} we have that $r_n^2=O(n^{-1})$). We therefore obtain from (\ref{eq:boundwass}) and (\ref{eq:bounTV}) that, if $\lambda$ is three times differentiable
  on $\Theta^{\star}$,
\begin{align*}
d_{\mathrm{Wass}}(\theta_{\mathrm{MAP}}^{\star}, N) \le \frac{\sqrt{d} }{\sqrt n}U_n(\mathbf{x})\{1+r_n(\mathbf{x},\eta)\}, \quad     d_{\mathrm{TV}}(\theta_{\mathrm{MAP}}^{\star}, N) \le    \frac{C_d}{\sqrt n}U_n(\mathbf{x})\{1+r_n(\mathbf{x},\eta)\},
\end{align*}  
where $U_n(\mathbf{x})=\| \matr
    F_{\overline \beta}(\hat \theta)\|_{\mathrm{diag}}  \mathfrak D_3(n\bar\beta,
  \hat\theta, \widehat{\matr K})$ depends solely on the sample size $n$ and sample data $\mathbf{x}$ and is independent of the prior distribution, and $r_n(\mathbf{x},\eta)$ is a remainder term (which is in general dependent on the prior distribution $\eta$) of order $o(1)$ as $n\rightarrow\infty$. 

Although we will tend to not spell this out for sake of conciseness, in all our examples in Section \ref{sec:applications}, our upper bounds are of the type 
the decomposition $\mathcal{U}_1(\mathbf{x})n^{-1/2}+\mathcal{U}_2(\mathbf{x},\bm\tau)n^{-3/2}$, where $\mathcal{U}_1(\mathbf{x})$ is independent of the prior hyperparameter $\bm\tau$, whilst $\mathcal{U}_2(\mathbf{x},\bm\tau)$ does in general depend on the prior hyperparameter. In contrast, our upper bounds on $d_{\mathrm{Wass}}(\theta_{\mathrm{MLE}}^{\star}, N)$ and $d_{\mathrm{TV}}(\theta_{\mathrm{MLE}}^{\star}, N)$ do not enjoy such a decomposition, being of the form $\mathcal{U}_3(\mathbf{x},\bm\tau)n^{-1/2}$; we also obtain lower bounds of the form $\mathcal{L}(\mathbf{x},\bm\tau)n^{-1/2}$ which suggest that in general under MLE standardisation we cannot expect to obtain upper bounds of the form $\mathcal{U}_4(\mathbf{x})n^{-1/2}+\mathcal{U}_5(\mathbf{x},\bm\tau)n^{-3/2}$. 
\end{example}

It remains to deal with the third order term in our bounds. 





\begin{lemma}\label{sec:main-resultsv2}
  Suppose that $g$ is three times differentiable on
   $ \Theta_0^{\star}$, and let $S(t)$ be a function such that $\sum_{u=1}^d|\partial^3_{u, v, w} g(t)|\leq S(t)$ for all $1 \le v, w \le d$ with $t\in\Theta_0^{\star}$.
  Then
   \begin{align*}
     \mathfrak{D}_{3}(g, \theta_0, \matr K \, | \,  \Theta_0^{\star})  & \le 
                                                   d  \mathbb{E} \left[
                                                   Y_1 S(R_2) 
                                                   \|
                                                                         \theta^{\star}\|^2_{\matr F} \,
                                                                         \mathbb{I}_{\Theta_0^{\star}}(\theta^{\star})
                                                   \right] 
   \end{align*}
   with $R_2 = \theta_0 + Y_1 Y_2 \matr J \theta^{\star}$ and $\|
                                                                         \theta^{\star}\|^2_{\matr
                                                                           F}
                                                                         =
                                                                         (\theta^{\star})^{\intercal}
                                                                         \matr
                                                                         F
                                                                         \theta^{\star}$.
   \end{lemma}

\begin{proof}
  Let $\Theta_0 = \matr J \Theta_0^{\star} + \theta_0$.  Then, by triangle
  inequality,
  \begin{align*} \mathfrak{D}_{3}(g, \theta_0, \matr K   \, | \,  \Theta_0^{\star}) & \le  
    \sum_{u=1}^d \sum_{v=1}^d \sum_{w=1}^d \mathbb{E} \left[ Y_1
      |\partial^3_{u, v, w} g( R_2)| \, |(\matr J\theta^{\star})_v|\,
      |(\matr J \theta^{\star})_w| \, \mathbb{I}_{\Theta_0^{\star}}(\theta^{\star})
    \right]
    \\
                                                            & \le 
                                                              \sum_{v=1}^d
                                                              \sum_{w=1}^d
                                                              \mathbb{E}
                                                              \left[
                                                              Y_1
                                                              S(R_2)
                                                              \,
                                                              |\theta_v
                                                              -
                                                              \theta_v^0|\,
                                                              |\theta_{w}
                                                              -
                                                              \theta_w^0|
                                                              \,
                                                             \mathbb{I}_{\Theta_0}(\theta)
                                                              \right] \\
                                                              & =                                                         \mathbb{E}
                                                              \bigg[
                                                              Y_1
                                                              S(R_2)
                                                              \,
                                                               \bigg( \sum_{v=1}^d
                                                               |\theta_v
                                                              -
                                                              \theta_v^0| \bigg)^2
                                                              \,
                                                             \mathbb{I}_{\Theta_0}(\theta)
                                                              \bigg]
     \le d     \mathbb{E}
                                                              \bigg[
                                                              Y_1
                                                              S(R_2)
                                                              \,
                                                                                                                               \sum_{v=1}^d(\theta_v
                                                              -
                                                              \theta^0_v)^{2}
                                                              \,
                                                             \mathbb{I}_{\Theta_0}(\theta)
                                                              \bigg],
\end{align*}
applying the Cauchy-Schwarz 
inequality on the sum. Finally, $  \sum_{v=1}^d(\theta_v
         -
          \theta^0_v)^{2}
          = (\matr
                              J
        \theta^{\star})^{\intercal}
   (\matr J \theta^{\star})
      $
           which
   gives the claim. 
   \end{proof}

\begin{remark}Sharper bounds can be obtained if certain third order partial derivatives of $g$ vanish. For example, suppose that
$\partial_{u,v,w}^3g(t)=0$ if $u,v,w$ all belong to the set $\{1,\ldots, d-1\}$ (
as in Example \ref{ex:linreg}). Then, 
from the proof of Lemma \ref{sec:main-resultsv2}, it is easily seen that
\begin{equation}\mathfrak{D}_{3}(g, \theta_0, \matr K \, | \,  \Theta_0^{\star})\leq\sum_{v=1}^d\mathbb{E}\Big[Y_1S(R_2)|(\theta_v-\theta_v^0)(\theta_d-\theta_d^0)|\mathbb{I}_{\Theta_0^{\star}}(\theta^{\star})\Big].
    \label{lembd}
\end{equation}
\end{remark}

  Lemma \ref{sec:main-resultsv2} is 
  simple to apply when $S$ is a bounded function.
  This assumption
   is 
   natural but the coarseness of the approach leads to a
   non-optimal dependence on the dimension. We can improve the
   dependence 
   considerably 
   under the next assumption.

\vspace{2mm}

\noindent {\bf Assumption ${\textrm P}$1. \label{assumptionP2}} There exist functions
$p_u(t), 1 \le u \le d$, and constants
$C_1, C_2, C_3$ such that, for all $t \in \mathcal{X}$,
$|\partial^3_{u, u, u} g(t)| \le C_1 p_u(t)$, 
$ |\partial^3_{u, u, v} g(t)| \le C_2 p_u(t) p_v(t) $  and
$|\partial^3_{u, v, w} g(t)| \le C_3 p_u(t) p_v(t) p_w(t)$ for all
$1 \le u \neq v \neq w \le d$. 



\begin{lemma}\label{sec:main-results}
  Suppose that $g$ satisfies Assumption \hyperref[assumptionP2]{P1} on $\Theta_0^{\star}$. 
  Let
  $C_{\infty}(t) := C_1+ 2 C_2 \sum_{u=1}^{d} p_u(t)+ C_3
  \big(\sum_{u=1}^{d} p_u(t)\big)^2.  $ Then
   \begin{align}\label{eq:d3gth}
     \mathfrak{D}_{3}(g, \theta_0, \matr K   \, | \,  \Theta_0^{\star})  & \le 
                                                   \mathbb{E} \big[
                                                   Y_1C_{\infty}(R_2)
                                                   \|
                                                                           \theta^{\star}\|^2_{\matr
                                                                           P(R_2)\matr
                                                                           F} \, \mathbb{I}_{\Theta_0^{\star}}(\theta)
                                                   \big] 
   \end{align}
   with $R_2 = \theta_0 + Y_1 Y_2 \matr J \theta^{\star}$ and 
   $\matr P(t) = \mathrm{diag}(p(t))$ is the diagonal matrix with
   entries $p_u(t)$, $1 \le u \le d$.
   
\end{lemma}

\begin{proof}
Let $\Theta_0 = \matr J \Theta_0^{\star} + \theta_0$.  Under Assumption \hyperref[assumptionP2]{P1} we can write
  \begin{align*} \mathfrak{D}_{3}(g, \theta_0, \matr K  \, | \,
    \Theta_0^{\star}) &
                                                                \le    \sum_{u=1}^d
                                                                                 \sum_{v=1}^d
                                                                                 \sum_{w=1}^d
                                                                                 \mathbb{E}
                                                                                 \left[
                                                                                 Y_1
                                                                                 |\partial^3_{u, v, w} g(
                                                                                 R_2)|
                                                                      \, 
                                                                                 |(\matr
                        J \theta^{\star})_{v}
                                                                                 |(\matr
                        J \theta^{\star})_{w}| \,
                                                                                 \mathbb{I}_{\Theta_0^{\star}}(\theta^{\star})
                                                                \right]
                                                            \le
                                                            \sum_{i=1}^4 I_i(g)
\end{align*}
with \begin{align*} I_1(g) & = C_1\sum_{u=1}^{d} \mathbb{E} \left[Y_1
    p_u(R_2) (\theta_u - \theta^0_u)^2 \,
    \mathbb{I}_{\Theta_0}(\theta) \right]
     \end{align*}
     and
     \begin{align*}
       I_2(g) & = C_2 \sum_{u=1}^{d} \sum_{v \neq u} \mathbb{E}
                \left[Y_1 p_u(R_2)p_v(R_2) (\theta_v - \theta^0_v)^2
                \, \mathbb{I}_{\Theta_0}(\theta)\right] \\
              & \le C_2 
                 \mathbb{E}
                \bigg[Y_1\Big(\sum_{v=1}^d p_v(R_2)\Big) \sum_{u=1}^{d}
                p_u(R_2) (\theta_u - \theta^0_u)^2 \,
                \mathbb{I}_{\Theta_0}(\theta)
                \bigg]
     \end{align*}
     (we switch notation in the summation indices in the last step for
     clarity),  
     \begin{align*}
       I_3(g) & = C_2\sum_{u=1}^{d} \sum_{w \neq u} \mathbb{E}
                \left[Y_1 p_u(R_2)p_w(R_2) |\theta_u - \theta^0_u|
                |\theta_w - \theta^{0}_w| \,\mathbb{I}_{\Theta_0}(\theta)\right]  \\
                           & \le C_2 \mathbb{E} \bigg[Y_1 \left(
                             \sum_{u=1}^{d} p_u(R_2) |\theta_u -
                             \theta^0_u| \right)^2 \,
                             \mathbb{I}_{\Theta_0}(\theta)\bigg]\\
&                             \le C_2 
                             \mathbb{E} \bigg[Y_1  
                       \Big(   \sum_{u=1}^{d}    p_u(R_2)\Big)  \sum_{u=1}^{d}    p_u(R_2) (\theta_u -
                             \theta^0_u)^2  \,
                             \mathbb{I}_{\Theta_0^{\star}}(\theta)\bigg],
                                                                             \end{align*}
                                                                             where
                                                                             the
                                                                             last
                                                                             inequality
                                                                             follows
                                                                             from
                                                                             applying
                                                                             the
                                                                             Cauchy-Schwarz
                                                                             inequality
                                                                             on
                                                                             the
                                                                             sum,
                                                                             and
                                                                             finally
     \begin{align*}
       I_4(g)& = C_3  \sum_{u=1}^{d} \sum_{v \neq u} \sum_{w \neq u, v}
            \mathbb{E} \left[ Y_1 p_u(R_2) p_v(R_2) p_w(R_2) |\theta_v
            - \theta^0_v| |\theta_w - \theta^0_w| \, \mathbb{I}_{\Theta_0}(\theta) \right]\\
       & \le C_3  
            \mathbb{E} \bigg[ Y_1  \Big(\sum_{u=1}^{d}  p_u(R_2) \Big)  \bigg( \sum_{v=1}^{d}  p_v(R_2)  |\theta_v
            - \theta^0_v|  \bigg)^2 \,\mathbb{I}_{\Theta_0}(\theta) \bigg]\\
       &  \le C_3     \mathbb{E} \bigg[Y_1  
                             \Big( \sum_{u=1}^{d} p_u(R_2)  \Big)^2 \sum_{u=1}^{d} p_u(R_2) (\theta_u -
                             \theta^0_u)^2  \,
                            \mathbb{I}_{\Theta_0}(\theta)\bigg]
     \end{align*}
     with the same justification.   
   \end{proof}

\begin{remark}\label{rk:main-results-1}
  Assumption \hyperref[assumptionP2]{P1} may seem ad hoc, but it is actually quite natural in
  our context. Indeed, taking $g=\beta$ as in model (\ref{eq:likili})
  then standard properties of exponential families in canonical form
  give (see Appendix \ref{sec:moments-expon-famil}), for all $1 \le u,
  v, w \le d$,
  \begin{equation*}
    \partial^3_{u,v,w}\beta(t) = \mathbb{E} \left[
      (h_u(X(t))-\mathbb{E}[h_u(X(t))])
      (h_v(X(t))-\mathbb{E}[h_v(X(t))])
      (h_w(X(t))-\mathbb{E}[h_w(X(t))]) \right]
  \end{equation*}
  with $X(t) \sim p_1(\cdot \, | \, t)$. If the $h_u(X(t))$ are
  {negatively correlated} (as is the case for multinomial data
  generating process, see Example \ref{ex:multi}) then Assumption \hyperref[assumptionP2]{P1} is satisfied.
\end{remark}

   \begin{remark}
     Note that for $x \in \mathbb{R}^d$ and $\matr X$ a $d\times d$
     symmetric matrix we have
$\| x\|_{\matr X}   = x^{\intercal} \matr X x \le \| \matr
    X \|_{\mathrm{sp}} \| x\|_2^2 $ 
with $\| \matr X \|_{\mathrm{sp}}$ the spectral norm defined in
Section \ref{sec:notation}.
   \end{remark}

The univariate case allows 
the following observation in
Wasserstein distance.
 As the identity function
$\mathrm{Id}:\mathbb{R}\rightarrow\mathbb{R} : x \mapsto
x$ is Lipschitz-1, we have the following lower bound: 
\begin{equation}
  \label{wasslow}
  d_{\mathrm{Wass}}(\theta^\star,N)\geq|\mathbb{E}[\theta^\star]|= {K}
  \big|\mathbb{E}[\theta]-
  \theta_0\big|. 
\end{equation}
  It is of interest to know when
the choice $\mathrm{Id}$ achieves the supremum in the Wasserstein
distance. In this context, the notion of stochastic ordering comes into
the picture.  We recall that two real-valued  random variables $X, Y$ are
stochastically ordered if the function
$x \mapsto \mathbb{P}[X \le x] - \mathbb{P}[Y \le x]$ does not change
sign over $x \in \mathbb{R}$.  The following result can be found e.g.\
in \cite[Theorem 1.A.11]{shaked2007stochastic}.
\begin{proposition}\label{prop:stochoooo}
 Let $X, Y$ have finite mean. Then $X, Y$ are stochastically ordered  if and only if
  $d_{\mathrm{Wass}}(X, Y) = |\mathbb{E}[X] - \mathbb{E}[Y]|. $
\end{proposition}

The next proposition ensues (see also \cite{yu2009stochastic} for related material).
 
\begin{proposition}\label{prop:stochord}
  If the function $ w \mapsto \lambda'(w)  - n \bar h(\mathbf{x}) - {K}^2 (w - \theta_0)$  does not change sign on the support of
$p^\star(\cdot | \mathbf{x})$  
then $\theta^\star$ and $N$ are
  stochastically ordered, and 
  $d_{\mathrm{Wass}}(\theta^\star, N) = {K} | \mathbb{E}[\theta]- {\theta_0}|.$
\end{proposition}

\begin{proof}
  
  A well-known sufficient condition for stochastic ordering of two
  absolutely continuous distributions is when their log-likelihood is
  monotone (see e.g.\ \cite{lrs17}). Straightforward computations show
  that this happens as soon as the stated condition holds.
\end{proof}


\section{Applications}
\label{sec:applications}
\subsection{Bernoulli data}
\label{example1}

Consider an observation $\mathbf{x}=(x_1, \ldots, x_n)$ sampled
independently from a Bernoulli distribution with parameter
$p \in (0, 1)$, a member of 
an
exponential family with natural parameter $\theta$ 
with 
\begin{equation*}
\theta=\log(p/(1-p)), \quad \Theta =(-\infty, \infty), \quad h_i(x) =
x, \quad \beta(\theta)=\log(1+e^\theta)
\end{equation*}
(recall \eqref{eq:likili}); $\beta$ is infinitely differentiable on $\Theta$.
With prior $\pi_0(\theta) = \exp({-\eta(\theta)})$
      the function  from
 \eqref{eq:lmbdaa} is 
    \begin{align*}
        \lambda(\theta) & =: n \bar\beta(\theta) + \eta(\theta)=n \beta(\theta) + \eta(\theta) =  n \log(1+e^{\theta}) + \eta(\theta).
        \end{align*}
We first consider   {normal
approximation of the posterior distribution centred at its mode} as in
Example~\ref{sec:bounds-mapINIT}, i.e.\ 
        \begin{align}\label{eq:postmode}
          \tilde \theta =
          (\lambda')^{-1}(n \bar x),\: \tilde{K} = \sqrt{\lambda''(\tilde
          \theta)}  \mbox{ and } \theta^{\star}_{\mathrm{MAP}} =
          \sqrt{\lambda''(\tilde \theta)} (\theta - \tilde \theta).
        \end{align}
We fix a
       {conjugate prior} on $\theta$, i.e.
        $\eta(\theta) = -\tau_1 \theta + \tau_2 \log(1+ e^\theta)$ for
        some $\tau_1>0$, $\tau_2-\tau_1>0$.  Then
        \begin{align*}
&           \lambda(\theta)  =  -\tau_1
        \theta + (n +\tau_2) \log(1+ e^\theta),  \quad        \lambda'(\theta)   = -\tau_1 +  (n +\tau_2) \frac{e^\theta}{1+ e^\theta}, \\
  &  \lambda''(\theta)   =  (n +\tau_2) \frac{e^\theta}{(1+ e^\theta)^2}, \mbox{ and } 
    \lambda'''(\theta)  =  (n +\tau_2) \frac{ e^\theta-e^{2\theta}}{(1+ e^\theta)^3}.
\end{align*}
For ease of notation (and also to connect with the forthcoming Example
\ref{ex:multi}) we introduce the quantities $\pi_1 = n\bar{x}+ \tau_1$
and $\pi_2 = n(1-\bar{x})+ \tau_2 - \tau_1$; note 
that
$\pi_1 + \pi_2 = n + \tau_2$ 
does not depend on the sample. 
Simple computations yield
 \begin{equation*} 
     \tilde\theta  = \log \left(\frac{\pi_1}{\pi_2}\right), 
\mbox{ and }
 \lambda''(\tilde \theta) = \frac{\pi_1\pi_2}{\pi_1 + \pi_2}.
   \end{equation*}
   We have $\tilde{K} 
  =   \sqrt{{\pi_1\pi_2}/({\pi_1 + \pi_2})}=   \sqrt{{\pi_1\pi_2}/({n+\tau_2})}$ and thus $\|\tilde F
   \|_{\mathrm{diag}} = \tilde J = \sqrt{1/\pi_1 + 1/\pi_2}$. 
   Since $0 \le ({ e^{x}-e^{2x}})/{(1+ e^{x})^3} \le {1}/{(6\sqrt3)}$
   for all $x \in \mathbb{R}$, 
  we can
   apply Lemma~\ref{sec:main-resultsv2} 
with $S(t)=1/(6\sqrt{3})$  to get 
   \begin{equation*}
 \Delta_3 (\lambda, \tilde \theta, \tilde K )  \le \frac{1}{12
       \sqrt 3} 
       \frac{(n+ \tau_2)^2}{\pi_1 \pi_2}
     \mathbb{E} \left[  (\theta^{\star}_{\mathrm{MAP}})^2
     \right] .
     \end{equation*}
It follows that the main term in our upper bounds in Theorem \ref{thm:abstract-results} satisfies  
     \begin{equation*}
          \tilde J \Delta_3 (\lambda, \tilde \theta, \tilde K 
 ) \le  \frac{1}{12
       \sqrt 3} 
       \frac{1}{\sqrt{n}} \frac{(1 + \tau_2/n)^{5/2} }{\big((\bar x + \tau_2/n)(1 - \bar x + (\tau_2 - \tau_1)/n) \big)^{3/2}}
           \mathbb{E} \left[  (\theta^{\star}_{\mathrm{MAP}})^2
\right].
     \end{equation*}

The conjugate prior setting allows us to compute the posterior moments
explicitly. 
With $\Gamma$ the 
gamma function, $\psi(t) = \Gamma'(t) / \Gamma(t)$ the digamma
function and $\rho(t) = \log(t) - \psi(t)$ we have 
\begin{align*}
\mathbb{E} \left[ \theta \right] - \tilde \theta  & =\rho(\pi_1)  -  \rho
                                                         (\pi_2),
                                                    \mbox{ and }
\mathbb{E} \big[ (\theta - \tilde \theta)^2 \big]  = \psi'(\pi_1) + \psi' (\pi_2)  
+ (\rho(\pi_1) - \rho(\pi_2))^2;
\end{align*}
see Appendices \ref{sec:moments-expon-famil} and \ref{appendixd1} for details.
Combining these formulas with  \eqref{eq:boundwass}
and \eqref{wasslow} we obtain 
\begin{align}
d_{\mathrm{TV}}(\theta^{\star}, N) & \le 
                                       \sqrt{\frac{\pi}{8}}
                                      U_n(\mathbf{x}, \bm \tau), \quad L_n(\mathbf{x}, \bm \tau)  \le
                 d_{\mathrm{Wass}}(\theta^{\star}_{\mathrm{MAP}}, N)  \le
                 U_n(\mathbf{x}, \bm \tau)  \label{simpwas}
\end{align}
with
\begin{align}
  L_n(\mathbf{x}, \bm \tau) & := \sqrt{\frac{\pi_1 \pi_2}{n+ \tau_2}}\left|\rho(\pi_1) - \rho(\pi_2) \right|,\label{eq:lowww} \\
  U_n(\mathbf{x}, \bm \tau) & :=  \sqrt{\frac{n + \tau_2}{\pi_1\pi_2} }\frac{1}{12
     \sqrt 3} \left(  \psi'(\pi_1) + \psi' (\pi_2)  
+ (\rho(\pi_1) - \rho(\pi_2))^2 \right).  \label{eq:boundwd3}
\end{align}

Supposing that the data follows a Bernoulli distribution with ground truth parameter $p^\star$, direct computations (using 
$\lim_{x \to \infty} x \rho(x) = 1/2$ and
$\lim_{x \to \infty} x \psi'(x) = 1$) 
give 
\begin{align}\label{wassliminf}    \liminf_{n \to \infty} \sqrt n \,
L_n( \mathbf{x}, \bm \tau) &\ge \frac{1}{2}
\frac{|1-2 p^\star|}{\sqrt{p^\star(1-p^\star)}}, \\ 
%
      \limsup_{n \to \infty} \sqrt n \,
U_n( \mathbf{x}, \bm \tau)  &\le \frac{1}{12
       \sqrt 3} 
        \frac{1}{(p^\star(1 - p^\star) )^{3/2}}. \label{wasslimsup} 
\end{align}

The lower bound in (\ref{wassliminf}) is equal to zero if $p^\star=1/2$, suggesting that the rate of convergence  may be faster than $O(n^{-1/2})$ in this case. We now use the fourth order expansions given in (\ref{eq:bounTV2}) to determine the rate of convergence of $d_{\mathrm{Wass}}(\theta^{\star}, N)$ and $d_{\mathrm{TV}}(\theta^{\star}, N)$ in the case $\bar{x}=1/2+O(n^{-1/2})$ (which would be the case for $p^*=0.5$, by the central limit theorem). 
We begin by bounding the quantity $\Delta_{4}(\lambda, \tilde \theta, \tilde K )$, defined in equation (\ref{delta4}). Simple calculations give that
$
\lambda'''(\tilde\theta)=\pi_1\pi_2(\pi_2-\pi_1)/(\pi_1+\pi_2)^2=\pi_1\pi_2(n(1-2\bar{x})+\tau_2-2\tau_1)/(\pi_1+\pi_2)^2   
$
and
$
|\lambda^{(4)}(x)|\leq(n+\tau_2)/8   
$ for all $x\in\mathbb{R}$, so that
\begin{align}
\tilde{J} \Delta_{4}(\lambda, \tilde \theta, \tilde K )\leq U_{2, n}(\mathbf{x}, \bm \tau)&:= \frac{1}{2(n+\tau_2)}\bigg(\frac{1}{\pi_1}+\frac{1}{\pi_2}\bigg)^{1/2}|n(1-2\bar{x})+\tau_2-2\tau_1\big| \mathbb{E} \left[  (\theta^{\star}_{\mathrm{MAP}})^2\right]\nonumber \\
    &\quad+\frac{n+\tau_2}{48}\bigg(\frac{1}{\pi_1}+\frac{1}{\pi_2}\bigg)^{2}\mathbb{E} \left[  \label{eq:upuD3v111} |\theta^{\star}_{\mathrm{MAP}}|^3\right]\\
    &=\frac{|1-2\bar{x}|}{2\sqrt{n\bar{x}(1-\bar{x})}}+\frac{1}{24}\sqrt{\frac{2}{\pi}}\frac{1}{(\bar{x}(1-\bar{x}))^2n}+r_n,  \label{eq:upuD3v112} 
\end{align}
where $r_n$ is a remainder term of order $O(n^{-3/2})$. In obtaining the above equality, we used that $ \mathbb{E} \left[  (\theta^{\star}_{\mathrm{MAP}})^2\right]=1+O(n^{-1/2})$ and $ \mathbb{E} \left[  |\theta^{\star}_{\mathrm{MAP}}|^3\right]=\mathbb{E}[|N|^3]+O(n^{-1/2})=2\sqrt{2/\pi}+O(n^{-1/2})$. From the inequalities in (\ref{eq:bounTV2}) we conclude that
\begin{align}
d_{\mathrm{TV}}(\theta^{\star}, N)  \leq\sqrt{\frac{\pi}{8}}                          U_{2, n}(\mathbf{x}, \bm \tau), \quad  d_{\mathrm{Wass}}(\theta^{\star}, N)  \leq                          U_{2, n}(\mathbf{x}, \bm \tau), \label{newtvwas}
\end{align}
and, since $1-2x\leq 1/(6\sqrt{3} x(1-x))$ for $0<x<1$, we see that \eqref{wasslimsup}  can be improved to 
\begin{align}\label{eq:upuD3v2}
  \limsup_{n \to \infty} \sqrt{n}  U_{2, n}(\mathbf{x}, \bm{\tau}) 
\leq\frac{1}{2}\frac{|1-2p^\star|}{\sqrt{p^\star(1-p^\star)}}; 
\end{align}
moreover, from \eqref{eq:upuD3v112},  if $\bar{x}=1/2+O(n^{-1/2})$ then the rate of convergence is of order $O(n^{-1})$. 

\begin{figure}[!]
    \centering
  \begin{subfigure}{.49\textwidth}
\captionsetup{width=.75\textwidth}
  \centering
  \includegraphics[width=6cm]{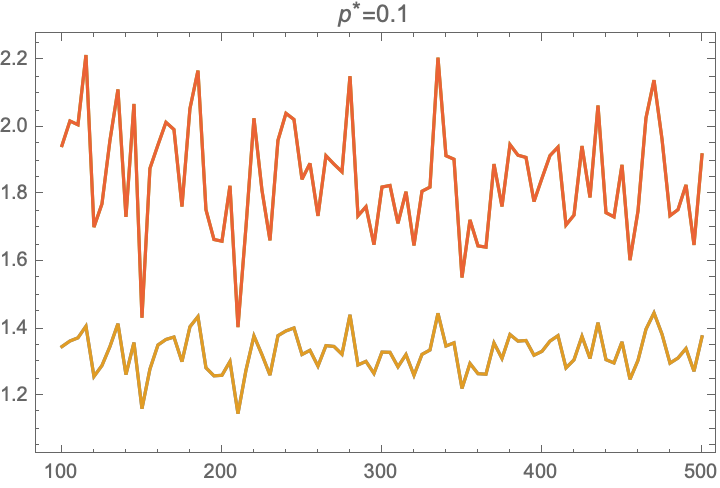}
    \caption{}\label{fig:sfig1berntv}
\end{subfigure}%
\begin{subfigure}{.49\textwidth}
\captionsetup{width=.75\textwidth}
  \centering
  \includegraphics[width=6cm]{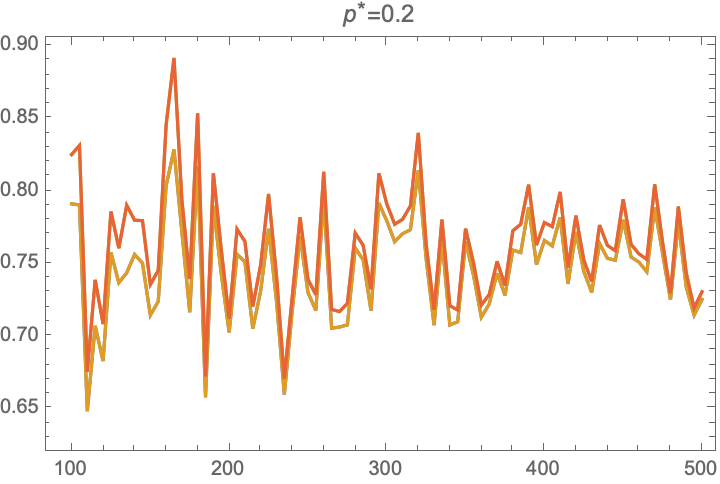}
  \caption{\label{fig:sfig1berntv2}}
\end{subfigure}%
\\
\begin{subfigure}{.49\textwidth}
\captionsetup{width=.75\textwidth}
  \centering
  \includegraphics[width=6cm]{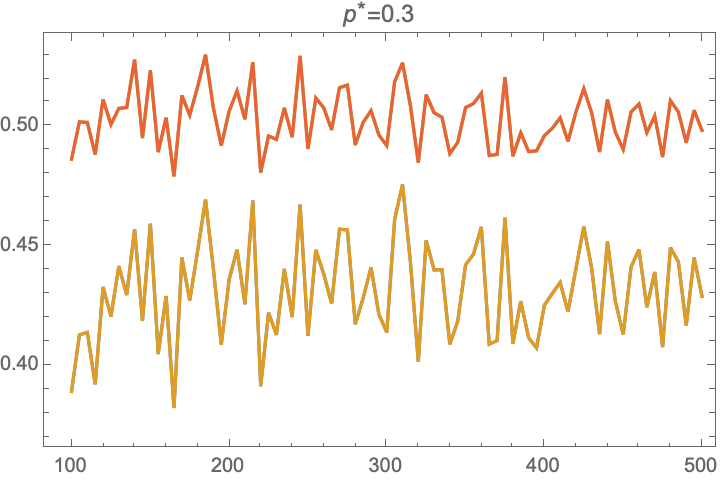}
  \caption{\label{fig:sfig1berntv3}}
\end{subfigure}%
\begin{subfigure}{.49\textwidth}
\captionsetup{width=.75\textwidth}
  \centering
  \includegraphics[width=6cm]{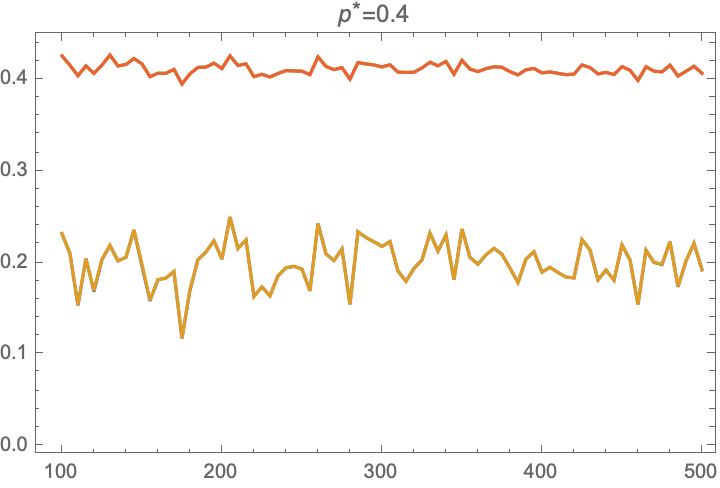}
  \caption{\label{fig:sfig1berntv3}}
\end{subfigure}%
\\
\begin{subfigure}{.49\textwidth}
\captionsetup{width=.75\textwidth}
  \centering
  \includegraphics[width=6cm]{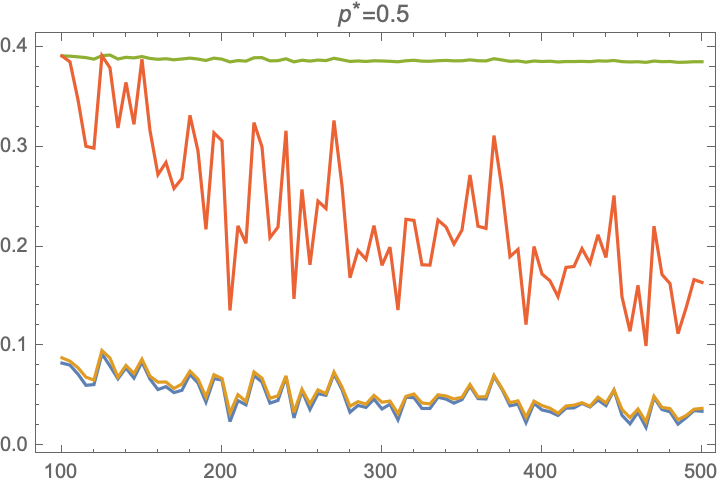}
  \caption{\label{fig:sfig1berntv5}}
\end{subfigure}
\begin{subfigure}{.49\textwidth}
\captionsetup{width=.75\textwidth}
  \centering
  \includegraphics[width=6cm]{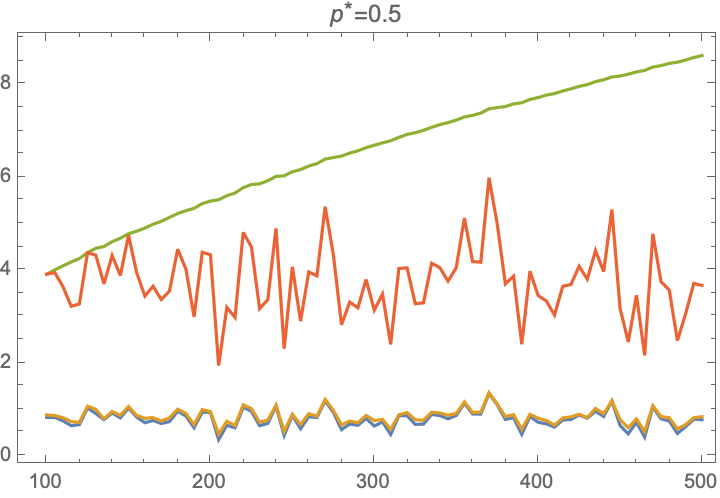}
  \caption{\label{fig:sfig1berntv5}}
\end{subfigure}%
\vspace{0.3cm}
\caption{\label{fig:binomial} \it 
This figure illustrates the results from Example \ref{example1} (Bernoulli data) using MAP centring and scaling with a conjugate prior. 
 \textbf{Orange curves:} True values of the Wasserstein distance, obtained numerically. 
 \textbf{Blue curves:} Lower bound of the Wasserstein distance, as 
 given by 
 equation \eqref{eq:lowww}. 
     \textbf{Red curves:} Upper bound of the Wasserstein distance, as given by equation \eqref{eq:upuD3v111}. 
     \textbf{Green curves:} Additional bound, computed using equation \eqref{eq:boundwd3}.
In panels {(a) to (e)}, the values are scaled by multiplying them by the square root of the sample sizes for the specified values of \( p^\star \). In panel {(f)}, the same data as in panel (e) is presented, but the values are instead multiplied by the sample size. 
Each curve represents the average of {10 simulations} performed for sample sizes \( n \in \{100, 120, \dots, 500\} \), all derived from the same Bernoulli dataset. The prior parameters are arbitrarily fixed at \( \tau_1 = 0.84 \) and \( \tau_2 = \sqrt{2} \). 
In many panels, the {blue curves overlap with the orange curves}, and the {green curves overlap with the red curves}, making the blue and green curves largely invisible.
}
\end{figure}

For the sake of illustration we computed the Wasserstein and total variation distance bounds in \eqref{simpwas},  
\eqref{newtvwas} for Bernoulli samples with
success parameters $p^\star=i/10$, $i=1, \dots,  5$, for sample sizes ranging
from $n=100$ to $n=500$; we also compute a numerical  approximation of the \emph{true} Wasserstein
and total variation distances since, in this simple setting, it is
easy to evaluate the corresponding integrals  numerically at least for moderate sample sizes 
up
to $n=500$.   We illustrate the results in Figure \ref{fig:binomial}, only in the case of the Wasserstein distance since the conclusions in the total variation distance are similar.
Several lessons can be learned from here.  First, it
appears that \eqref{eq:lowww} is in fact the exact value of the
Wasserstein distance in this case (in view of Proposition
\ref{prop:stochoooo} we 
conjecture the two
distributions are stochastically ordered
) and also, up to a multiplicative constant, the
exact value of the total variation distance.  Second, 
when  $p^\star=0.5$ the true  rate of convergence is of order $n^{-1}$; the third order upper bound from \eqref{simpwas} misses this information 
while the fourth order upper bound in  \eqref{newtvwas} 
 captures it.  
We illustrate this in subfigures (e) and (f), by changing the scaling factor on the same data.   We believe that the fact that the blue curve is sometimes slightly below the orange curve is  explained by numerical instability in the integrals.  Finally, we note that the proximity between the lower and upper bounds visible in panel (b) is 
the result of an intriguing situation:   the lower bound and the two upper bounds yield the same value when $p^\star = (3 \pm \sqrt{3})/6$ (i.e. $p^\star = 0.211... $ and $p^\star = 0.788...$);  hence whenever $\bar x$ is close to this value, the bounds are nearly identical.

Next, we study {normal approximation of the posterior
  standardised around the MLE}. Here $\theta$ is distributed according
to $p_2( \cdot \, | \, \mathbf{x})$ as in \eqref{eq:post} with some prior
$\eta: \mathbb{R} \to \mathbb{R}$, which we assume to have bounded third and fourth order derivatives on $\Theta$  and
\begin{align*}
 \theta^{\star}_{\mathrm{MLE}} =
          \sqrt{n\beta''(\hat \theta)} (\theta - \hat \theta) \mbox{
   with }   \hat  \theta =
          (\beta')^{-1}(\bar x). 
\end{align*}
Since  $\beta(\theta) =
\log(1+e^{\theta})$,
solving $\beta'(\theta) = \bar x$, we obtain 
$
   \hat\theta =\log (\bar{x}/(1-\bar{x})).
$
We also have
$\widehat K = (n \beta''(\hat \theta))^{1/2} = \sqrt{n \bar x (1 - \bar
  x)}$. Then
$\|\widehat F \|_{\mathrm{diag}} = \widehat J = 1/\sqrt{n \bar x(1 -
  \bar x)}$. 
  Also, again on using Lemma \ref{sec:main-resultsv2} to bound $\mathcal{D}_3(\lambda, \hat \theta, \widehat K  )$, we obtain
\begin{align*}
 &  \mathcal{D}_1(\lambda, \hat \theta, \widehat K  ) \le  
 |  \eta' ( \hat \theta) |, \quad \mathcal{D}_2(\lambda, \hat
  \theta,\widehat K )  =
|    \eta''(  \hat \theta  )|
\mathbb{E}[|\theta- \hat \theta|] ,\\
&\mathcal{D}_3(\lambda, \hat \theta, \widehat K  )\leq\frac{1}{\bar x(1-\bar x)}\bigg(\frac{1}{12\sqrt{3}}+\frac{\|\eta^{(3)}\|_\infty}{n}\bigg)\mathbb{E}\left[(\theta_{\mathrm{MLE}}^*)^2\right],\\
&\lambda^{(3)}(\hat\theta)= \bar{x}(1-\bar x)(1-2\bar x)+\eta^{(3)}(\hat\theta), \quad \|\lambda^{(4)}\|_\infty\leq n/8+\|\eta^{(4)}\|_\infty.
 \end{align*}
It follows that  the main term in the  bounds from Theorem \ref{thm:abstract-results} satisfies 
\begin{align}
  &  \widehat J    \Delta_3 (\lambda, \hat
    \theta, \widehat  K )\leq
    \frac{ |\eta'(\hat \theta)| }{\sqrt{n \bar x(1 -
    \bar x)}}   +  \frac{ |\eta''(\hat \theta)| }{n \bar x(1 -
    \bar x)}   \mathbb{E}[|\theta_{\mathrm{MLE}}^*| ] 
    +    \frac{  \frac{1}{12 \sqrt 3} +
    \frac{ \|\eta^{(3)}\|_\infty
    }{2n}}{\sqrt{n}(\bar x(1 -
    \bar x))^{3/2}}  \mathbb{E}\left[(\theta_{\mathrm{MLE}}^*)^2\right]
                                               \label{eq:boundmlegen}
\end{align}
which gives bounds on the total variation and Wasserstein distances in
this case as well.
Also, bounding the quantity $\Delta_{4}(\lambda, \hat \theta, \widehat K )$, defined in equation (\ref{delta4}) gives that
\begin{align*}
\widehat J  \Delta_{4}(\lambda, \hat \theta, \widehat K )
&\leq  \frac{ |\eta'(\hat \theta)| }{\sqrt{n \bar x(1 -
    \bar x)}}   +  \frac{ |\eta''(\hat \theta)| }{n \bar x(1 -
    \bar x)}   \mathbb{E}[|\theta_{\mathrm{MLE}}^*| ]+\frac{|1-2\bar x|}{2\sqrt{n\bar x(1-\bar x)}}\bigg(1+\frac{|\eta^{(3)}(\hat\theta)|}{n\bar x(1-\bar x)}\bigg)\mathbb{E}\left[(\theta_{\mathrm{MLE}}^*)^2\right] \\
&\quad+\frac{1}{48n(\bar x(1-\bar x)^2}\bigg(1+\frac{8\|\eta^{(4)}\|_\infty}{n}\bigg) \mathbb{E}\left[|\theta_{\mathrm{MLE}}^*|^3\right].    
\end{align*}
In the case of a
conjugate prior with exponent
$\eta(\theta)=-\tau_1\theta+\tau_2\log(1+e^\theta)$ we have that
\begin{align}
\widehat J \Delta_{4}(\lambda, \hat\theta, \widehat K )
&\leq\frac{|\tau_1-\tau_2\bar{x}|}{\sqrt{n\bar{x}(1-\bar{x})}}+\frac{\tau_2}{n} \mathbb{E} \left[  |\theta^{\star}_{\mathrm{MLE}}|\right] +\frac{(n+\tau_2)|1-2\bar{x}|}{2n^{3/2}\sqrt{\bar{x}(1-\bar{x})}}\mathbb{E}  \left[  (\theta^{\star}_{\mathrm{MLE}})^2\right]\nonumber \\
    &\quad+\frac{n+\tau_2}{48(n\bar{x}(1-\bar{x}))^2}\mathbb{E} \left[  |\theta^{\star}_{\mathrm{MLE}}|^3\right]\label{eq:newbound4}\\
    &=\label{eq:newbound4v2} \frac{\tau_2|1-2\bar{x}|+2|\tau_1-\tau_2\bar{x}|}{2\sqrt{n\bar{x}(1-\bar{x})}}+\sqrt{\frac{\pi}{2}}\frac{\tau_2}{n}+\frac{1}{24}\sqrt{\frac{2}{\pi}}\frac{1}{(\bar{x}(1-\bar{x}))^2n}+r_n, 
\end{align}
where $r_n$ is a remainder term of order $O(n^{-3/2})$. Observe that the bound (\ref{eq:newbound4v2}) is of order $n^{-1}$ if $\bar{x}=1/2+O(n^{-1})$ and $\tau_2=2\tau_1$.
In case of a flat (non-informative)
prior
we also get that
\begin{align*}
\widehat J \Delta_{4}(\lambda, \hat\theta, \widehat K )\leq \frac{|1-2\bar x|}{2\sqrt{n\bar x(1-\bar x)}}\mathbb{E}\left[(\theta_{\mathrm{MLE}}^*)^2\right] +\frac{1}{48n(\bar x(1-\bar x)^2} \mathbb{E}\left[|\theta_{\mathrm{MLE}}^*|^3\right].   
\end{align*}

As in the MAP centric case, we performed some numerical explorations over the same configurations. These lead to similar conclusions as those reported in Figure \ref{fig:binomial}, with no further  surprise in the case of a flat prior. In the case of a conjugate prior, however, we read from \eqref{eq:newbound4v2} that the rate of convergence will no longer drop to $n^{-1}$ when $p^\star = 0.5$ \emph{unless} $\tau_1 = 2 \tau_2$. This is 
also indicated in our simulations, see Figure \eqref{fig:binomial3}.

\begin{figure}[h]
    \centering
\vspace{.2cm}
 \begin{subfigure}{.49\textwidth}
\captionsetup{width=.75\textwidth}
  \centering
  \includegraphics[width=6cm]{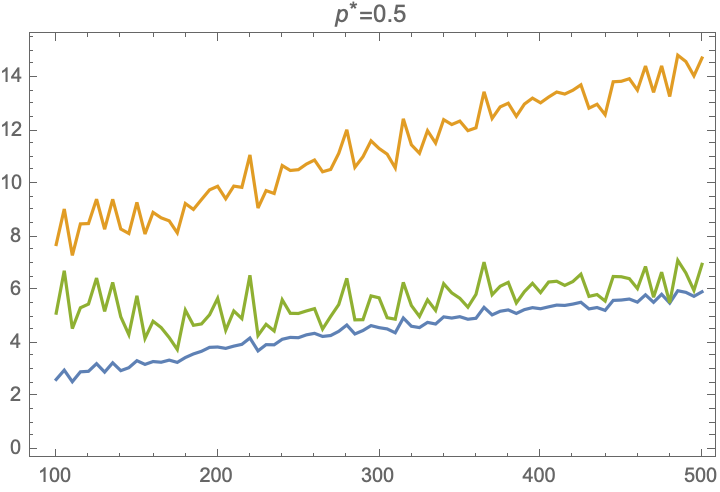}
  \caption{}
\end{subfigure}%
 \begin{subfigure}{.49\textwidth}
\captionsetup{width=.75\textwidth}
  \centering
  \includegraphics[width=6cm]{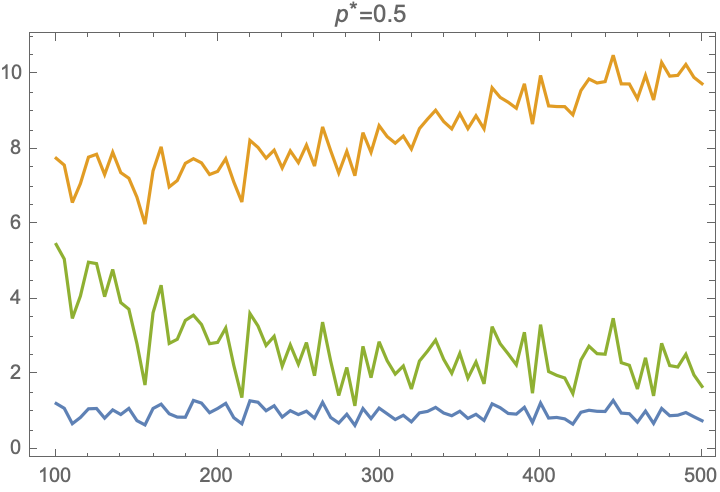}
  \caption{}
\end{subfigure}%

\vspace{0.3cm}
\caption{\label{fig:binomial3} \it Exactly the same configuration as in Figure \ref{fig:binomial}, but with MLE centring and scaling using a conjugate prior. The values are multiplied by the sample size. 
     \textbf{Blue curves:} True values of the Wasserstein distance.
     \textbf{Orange curves:} Bound as defined in equation \eqref{eq:boundmlegen}.
     \textbf{Green curves:} Bound as given in equation \eqref{eq:newbound4}.
All curves are multiplied by the sample size.  
In the left panel ({panel (a)}), the parameters are set to  \( \tau_1 = 0.84 \) and \( \tau_2 = \sqrt{2} \). In the right panel ({panel (b)}), the parameters are set to \( \tau_1 = 1 \) and \( \tau_2 = 2 \). In both panels, \( p^\star = 0.5 \).
}
\end{figure}

\subsection{Poisson data}
\label{poissonex}

Consider i.i.d.\ Poisson data with parameter $\mu>0$.  
This distribution is an exponential family with natural parameter $\theta$ with
\begin{equation*}
    \theta=\log \mu, \quad \Theta = (-\infty, \infty), \quad h_i(x) =
    x, \quad \beta(\theta)=e^\theta. 
  \end{equation*}
 We focus on {normal approximation of the posterior
      distribution centred at its mode} so that \eqref{eq:postmode} is
    the variable of interest.  We consider a {conjugate gamma
      prior} on $\mu$ of the form
    $\pi_0(\theta)\propto \exp({-\eta(\theta)})$ with
    $\eta(\theta) = -\tau_1 \theta + \tau_2 e^\theta$ for
    $\tau_1,\tau_2>0$. Then
$
    \lambda(\theta) = - \tau_1 \theta + (n + \tau_2) e^\theta,
$
so that 
$
       \lambda'(\theta)  = -\tau_1 +  (n +\tau_2) e^\theta$ and
       $\lambda^{(k)}(\theta) = (n +\tau_2) e^\theta$ for all  $k \ge 2$. 
Solving $\lambda'(
\theta) = n \bar x$ we obtain 
$$\tilde\theta=\log \left( \frac{\tau_1+ n \bar x }{n+\tau_2} \right)
\mbox{ and } \lambda^{(k)}(\tilde \theta) = \tau_1+ n \bar x, \quad
k\geq2.$$
In particular, $\tilde K = \sqrt{\tau_1 + n \bar x}$ and $\tilde J =
1/\sqrt{\tau_1+ n \bar x}$. 

A {Wasserstein distance bound} is easy to obtain. Indeed, since
$\lambda'''(\theta)$ is always positive, Proposition~\ref{prop:stochord} applies and we immediately obtain
 \begin{align}
\label{eq:wasspoi}
         d_\mathrm{Wass}( \theta^\star_{\mathrm{MAP}}, N)  
 = \sqrt{n\bar x + \tau_1}  \left(\log(n\bar x + \tau_1) - \psi(n\bar
   x + \tau_1)\right)  
   ,
 \end{align}
and thus, if the ground truth parameter is $\lambda^\star$, we get  
 $$\lim_{n\to \infty} \sqrt{n} d_\mathrm{Wass}( \theta^\star_{\mathrm{MAP}}, N) = 1/(2
 \sqrt{\lambda^\star})$$ for all admissible values $\tau_1, \tau_2$. 

 For a {total variation distance bound}, Lemma \ref{sec:main-resultsv2}
 (or Lemma \ref{sec:main-results}) can be applied but $\lambda'''$ is
 unbounded so this does not bring a true simplification. We use
 \eqref{eq:TVbounnnnthetastarcond}. We first compute
 \begin{align*}\mathbb{E}\left[Y_1  \lambda'''(\tilde \theta + Y_1Y_2
                                     w ) \right]    =(\tau_1 + n \bar x) \mathbb{E} \left[Y_1\exp \left(
                                                                    Y_1Y_2w  \right) \right]
                                                                    =\frac{
                                                                    \tau_1+n
                                                                    \bar
                                                                    x
                                                                    }{w^2}
                                                                    \left(
                                                                    e^{w
                                            }-1
                                                                    -
                                                                    w\right).
 \end{align*}
Next we note that for $|w| \le {\epsilon}$ it  holds that 
$
  |e^w-1 - w| \le  (w^2/2) e^{ \epsilon}
$
so that with
$\Theta_0^{\star} = (- \tilde K \epsilon, \tilde K \epsilon)$ we get 
that the various terms in \eqref{eq:TVbounnnnthetastarcond} satisfy
\begin{align*}
 &  \mathfrak{D}_3(\lambda, \tilde \theta, \tilde K \, | \, (- \tilde K
  \epsilon, \tilde K \epsilon)) \le \frac{e^{\epsilon}}{2}
   \mathbb{E} \left[ (\theta^{\star}_{\mathrm{MAP}})^2 \right],  \\
  & \mathbb{P} \left[ \theta^{\star}_{\mathrm{MAP}} \notin (- \tilde K
    \epsilon, \tilde K \epsilon)  \right] \le \frac{\mathbb{E} \left[
    (\theta^{\star}_{\mathrm{MAP}})^2  \right]}{\tilde
    K^2\epsilon^2},\quad r^{\star}_1( \tilde K \epsilon) =\sqrt{\frac{\pi}{8}} \max(p^{\star}(\tilde K
      \epsilon),p^{\star}(-\tilde K
      \epsilon))
\end{align*}
(for the first inequality we bound the indicator function by 1, for
the second we use Markov's inequality).  Recalling \eqref{eq:pstarlambda} we obtain after some simplifications 
\begin{align*}
  p^{\star}(\pm \tilde K \epsilon  )  &=   \frac{(n \bar x+ \tau_1)^{n \bar x  + \tau_1 - \frac12}}{\Gamma(n \bar x + \tau_1)} e^{-(n\bar x + \tau_1) (e^{\pm \epsilon} \mp \epsilon)} < \frac{1}{\sqrt{2\pi}}\exp\big(-(n\bar x + \tau_1) (e^{\pm \epsilon} \mp \epsilon-1)\big)\\
  &\leq \frac{1}{\sqrt{2\pi}}\exp\big(-(n\bar x + \tau_1) (e^{- \epsilon} + \epsilon-1)\big)
  =:r_n(\epsilon),
\end{align*}
where we used Stirling's inequality to obtain the first inequality, and the basic inequality $e^y-y\geq e^{-y}+y$, for      $y>0$, to obtain the second inequality. Since $e^{-\epsilon} + \epsilon-1 > 0 $ for all $\epsilon>0$, the term $r^{\star}_1( \tilde K \epsilon)< r_n(\epsilon)$ is exponentially negligible in $n$ for  $\epsilon$ not too small. 
We now have for all $\epsilon>0$, 
\begin{align}\label{produclog}
    & d_\mathrm{TV}( \theta^\star_{\mathrm{MAP}}, N) 
 \le  \frac{1}{\sqrt{\tau_1 + n \bar{x}}}\left( \sqrt{\frac{\pi}{8}}\frac{e^{\epsilon}}{2}
    + \frac{1}{\epsilon^2} \frac{1}{\sqrt{\tau_1 + n \bar x}} \right)  {\mathbb{E}[(\theta^{\star}_{\mathrm{MAP}})^2]}+\sqrt{\frac{\pi}{8}}r_n(\epsilon).
\end{align}
 Since for all $n, \bar x$ such that
$n \bar x + \tau_1 \ge 27$ there exists $\epsilon$ such that
$\sqrt{\frac{\pi}{8}}\frac{e^{\epsilon}}{2} + \frac{1}{\epsilon^2}
\frac{1}{\sqrt{\tau_1    + n \bar x}} \le 1$ (here we can take $\epsilon = 0.85489$),  
we can conclude that
\begin{align}
    d_\mathrm{TV}( \theta^\star_{\mathrm{MAP}}, N)  \le \frac{1}{\sqrt{\tau_1 + n \bar{x}}}{\mathbb{E}[(\theta^{\star}_{\mathrm{MAP}})^2]} 
\end{align}
which gives a bound on the rate of convergence of order $n^{-1/2}$ since $\mathbb{E}[ (\theta^{\star}_{\mathrm{MAP}})^2]=1+o(1)$ as $n\rightarrow\infty$.  
We can, as in Example
\ref{example1}, use the conjugate prior setting to compute directly 
\begin{align*}
    \mathbb{E}[(\theta - \tilde \theta)^2]  &=
                                              \psi'(n\bar{x}+\tau_1)+
                                              (\psi(n\bar{x}+\tau_1) -
                                              \log
                                              (n\bar{x}+\tau_1))^2. 
\end{align*}
This yields
\begin{align}\label{eq:TVBOUNDPOISS}
  d_\mathrm{TV}( \theta^\star_{\mathrm{MAP}}, N)  \le
  \sqrt{\tau_1+n \bar
  x}  
  \left( \psi'(n\bar{x}+\tau_1)+
  (\psi(n\bar{x}+\tau_1) -
  \log
  (n\bar{x}+\tau_1))^2
  \right)  \lesssim \frac{1}{\sqrt{\tau_1 + n \bar x}}.
\end{align}
Here $\lesssim$ indicates an inequality up to an unspecified absolute constant.
 Numerical explorations along the lines
of Figure \ref{fig:binomial} indicate that the bound
\eqref{eq:TVBOUNDPOISS} should be an equality up to a (nearly) constant factor
of roughly 8, see Figure \ref{fig:2poi}. Finally, as in the previous
example, the total variation and Wasserstein distances appear to be
(nearly) proportional to each other.

\begin{figure}[!]
    \centering
\vspace{.2cm}
   \begin{subfigure}{.49\textwidth}
\captionsetup{width=.75\textwidth}
  \centering
  \includegraphics[width=6cm]{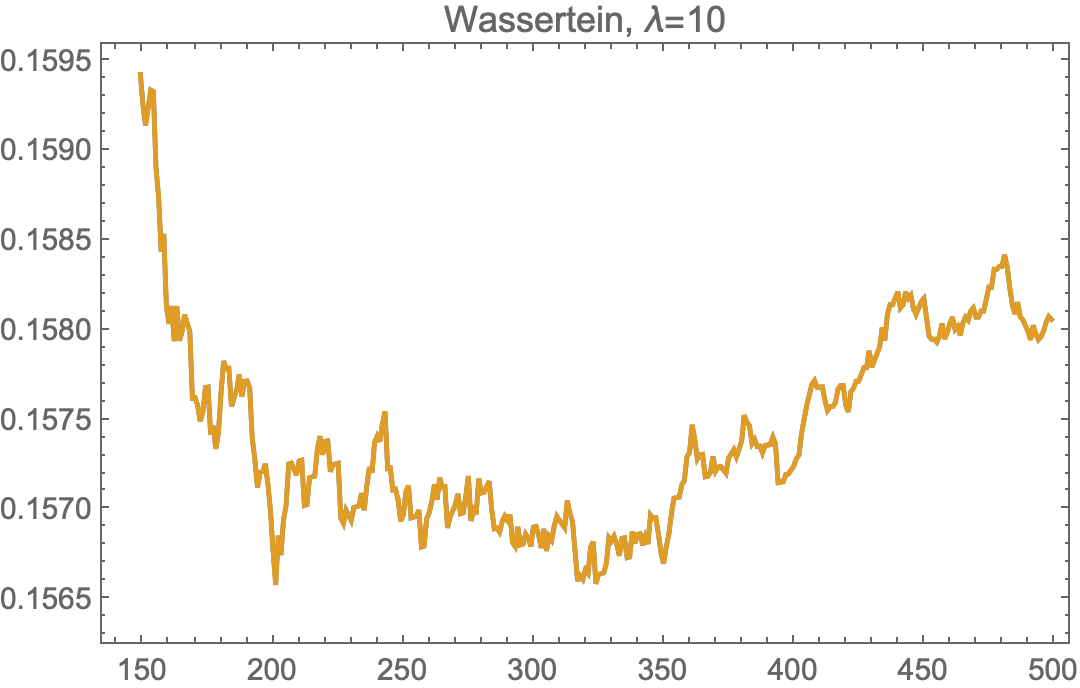}
  \label{fig:WassPOI}
\end{subfigure}%
\begin{subfigure}{.49\textwidth}
\captionsetup{width=.75\textwidth}
  \centering
  \includegraphics[width=6cm]{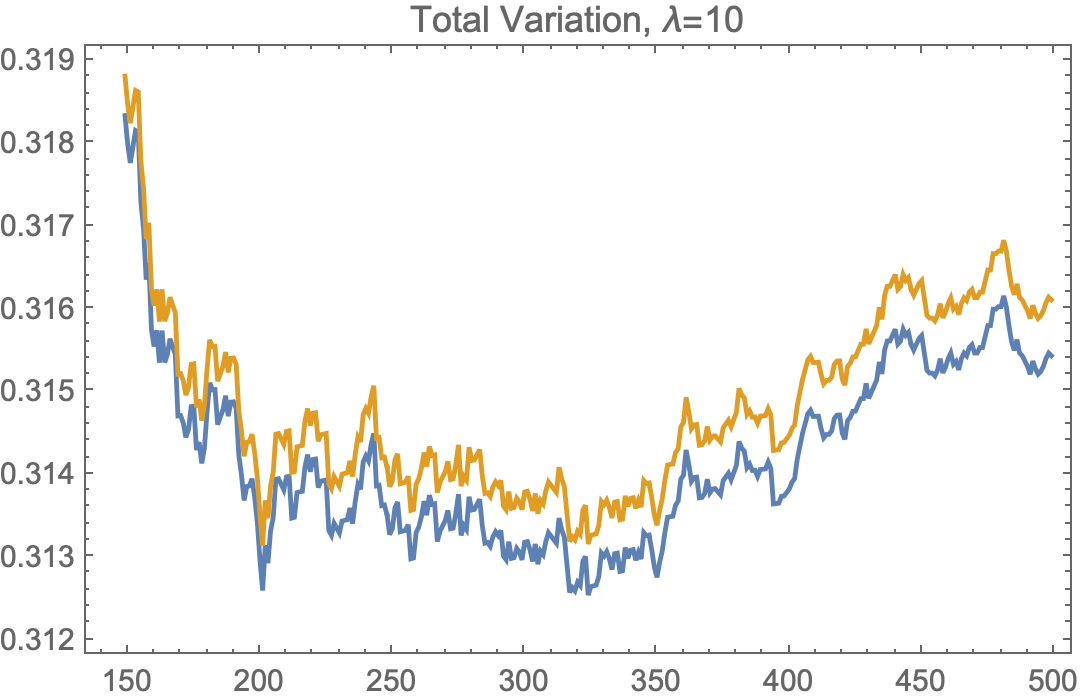}
    \label{fig:TVpPOI}
\end{subfigure}%

\vspace{0.3cm}
\caption{\label{fig:2poi} \it The \textbf{blue curves} represent the true values of the Wasserstein distances (left column) and total variation distances (right column) multiplied by the square root of the sample sizes; the total variation distance curve is further multiplied  by a factor of 7.5. The \textbf{orange curves} give the upper bounds, computed via equation \eqref{eq:wasspoi} (left column) and equation \eqref{eq:TVBOUNDPOISS} (right column), also multiplied by the square root of the sample size. 
The orange and blue curves are, obviously,  indistinguishable in the left panel. All curves are computed for the same sample of Poisson data with parameter 
\( \lambda = 10 \), for sample sizes \( n \in \{100, 120, \dots, 500\} \). The prior parameters are set to \( \tau_1 = 1 \) and \( \tau_2 = 3 \). 
}
\end{figure}

\subsection{Normal  with known mean data}
  \label{normalex}

    Consider i.i.d.\ normal data with known mean  $m$ and unknown variance $\sigma^2$. This distribution is written as an exponential family with natural parameter $\theta$ by taking 
\begin{equation*}
    \theta = \sigma^{-2}, \quad \Theta = (0, \infty), \quad h(x)=-(x-m)^2/2 \mbox{ and }\beta(\theta)=-(1/2)\log\theta.
\end{equation*}
 We denote the prior $\pi_0(\theta) \propto  \exp({-\eta(\theta)})$
    so that   
$
        \lambda(\theta) =  - (n/2) \log \theta + \eta(\theta).  
$
    Note that the function $\beta$ is infinitely differentiable.

We consider {normal approximation  of the posterior
  standardised around the posterior mode}, with  a  {conjugate gamma prior} on $\theta$  of the form $\pi_0(\theta)  \propto \exp({-\eta(\theta)})$, where 
$\eta(\theta) =  \tau_1\theta  
 -\tau_2 \log \theta$ for $\tau_1>0$ and $\tau_2>-1$.  The posterior distribution $\theta$ has a gamma distribution with shape parameter $n/2+\tau_2+1$, and so Assumption A\ref{item:A2} is satisfied if $n+2\tau_2>0$. 
Then 
$  \lambda(\theta) = \tau_1 \theta 
 -(\tau_2 + n/2) \log \theta,$
so that, for $n\geq2$, 
\begin{align*}
       \lambda'(\theta)   = \tau_1 -   \frac{n/2+\tau_2}{\theta}, \quad 
    \lambda''(\theta)  = \frac{n/2+\tau_2}{\theta^2} \mbox{ and } 
        \lambda'''(\theta)  = -\frac{n+2\tau_2}{\theta^3}.
\end{align*}
Let  $n{s_n}^2 = \sum_{i=1}^n(x_i-m)^2$. Solving $\lambda'(w) = -n s_n^2/2$ we obtain 
\begin{align*}
      \tilde\theta  =  \frac{n+2\tau_2}{n s_n^2+2\tau_1} \mbox{ so that } \lambda''(\tilde \theta) = \frac{(n s_n^2+2\tau_1)^2}{2(n+2\tau_2)}.
\end{align*}
In particular, $\tilde K = (n s_n^2+2\tau_1)/\sqrt{2(n+2\tau_2)}$ and
  $\tilde J = \sqrt{2(n+2\tau_2)}/(n s_n^2+2\tau_1)$.
 
We have  \begin{align*}
\left| \mathbb{E} \left[ Y_1  \lambda'''(\tilde \theta + Y_1Y_2w) \right]  \right| & = \left| \mathbb{E} \left[ Y_1 \frac{n+2\tau_2}{  \big( Y_1
  Y_2 w +\frac{n+2\tau_2}{n
  s_n^2+2\tau_1}  \big)^3}  \right]  \right| =  \frac{1}{2}  \frac{ (n
  s_n^2 + 2 \tau_1)^2}{(n + 2 \tau_2) \big|w + \frac{(n+2 \tau_2 )}{(n
  s_n^2 + 2 \tau_1) } \big|}.
\end{align*}
We get, after some straightforward simplifications,
\begin{align*}
  \mathfrak{D}_3 (\lambda, \tilde \theta, \tilde K) &
\leq                                                       
                                                               \frac{(n 
                                                               s_n^2+2\tau_1)^2}{2({n+2\tau_2})}
      \mathbb{E} \left[    
                                                               \frac{(\theta
                                                               -
                                                               \tilde
                                                               \theta)^2}{|\theta|} \right] 
                                                                 =
                                                                 \frac{n
                                                                 s_{n}^2+
                                                                 \tau_1
                                                                 }{n +
                                                                 2 \tau_2}
 \end{align*}
 (the equality follows since $\theta$ is gamma distributed with known
 parameters).  Hence from \eqref{eq:bounTV} we get
\begin{align} \label{bound_normal_precision_TV}
     d_\mathrm{TV}( \theta^\star_{\mathrm{MAP}}, N) &\leq \sqrt{\frac{\pi}{4}} \frac{1}{\sqrt{n+2\tau_2}}.
\end{align}
Similarly, the Wasserstein inequality \eqref{eq:boundwass}  yields 
$d_\mathrm{Wass}( \theta^\star_{\mathrm{MAP}}, N) \leq
\sqrt{2}/\sqrt{n+2\tau_2}$. Also, we can use inequality
(\ref{wasslow}) to obtain a lower bound for the Wasserstein distance
$d_{\mathrm{Wass}}(\theta^\star_{\mathrm{MAP}},N)\geq
{\sqrt{2}}/{\sqrt{n+2\tau_2}},$ whence
\begin{align*} d_{\mathrm{Wass}}(\theta^\star_{\mathrm{MAP}},N)= \frac{\sqrt{2}}{\sqrt{n+2\tau_2}}.  
\end{align*}
Numerical explorations indicate that the total
variation bounds are off from the simulated total variation distance by a constant factor and the total variation distance
appears, as for the previous examples and for some parameter
constellations, to be a direct multiple of the Wasserstein
distance. 
The specific values taken by the sample have no
impact on the Wasserstein distance and appear to have little to no impact on the
total variation distance. When $\tau_2$ is very large the distance is
small.

Next, we study {normal approximation of the posterior
  standardised around the MLE}, focusing on Wasserstein
distance. Here $\beta(\theta) = -(1/2)\log {\theta}$ and $\beta^{(\ell)}(\theta)=(-1)^\ell(\ell-1)!/(2\theta^\ell)$ for $\ell\geq1$.
Solving $\beta'(\theta) = - s_n^2/2$ we obtain
\begin{equation*}
  \hat \theta = 1/s_n^2 
\end{equation*}
so that ${\beta}''(\hat
  \theta) =  s_n^4/2$. Thus, $\widehat K = \sqrt{n/2}s_n^2$ and $\widehat J = \sqrt{2}/ (\sqrt n s_n^2)$.
We take note of the identity
\begin{align*}
 \left|  \mathbb{E} \left[
 {Y_1} \beta'''(\hat \theta + Y_1Y_2w)  
  \right]  \right| =   \bigg|  \mathbb{E} \bigg[
 \frac{Y_1}{(Y_1Y_2 w  + \hat \theta)^3} 
  \bigg]  \bigg| = \frac{1}{2 (\hat \theta)^2} \frac{1}{|w + \hat \theta|} 
\end{align*}
which, after simplifications, gives 
\begin{align*}
  \mathcal{D}_3(\lambda, \hat \theta, \widehat K)
  \le    \frac{(s_n^2)^2 }{2}  \mathbb{E} \bigg[ 
    \frac{(
    \theta-\hat \theta )^2}{|\theta|} \bigg]    + 
  \mathbb{E} \left[ Y_1\eta'''(R_2) (
    \theta-\hat \theta )^2\right]  
\end{align*}
with $R_2 = \hat \theta + Y_1 Y_2 (\theta - \hat \theta)$ (we express
everything in terms of $\theta$ rather than $\theta^{\star}$ because
this will simplify the computations slightly).  From
\eqref{eq:boundwass} (also recall Example \ref{sec:bounds-mle}) we get
\begin{align*}
     d_{\mathrm{Wass}}(\theta^{\star}_{\mathrm{MLE}}
    , N) &  \le \sqrt{\frac{2}{n}}\frac{|\eta'(s_n^{-2})|}{s_n^2}
           + \sqrt{\frac{2}{n}} \frac{ \left| \eta''(s_n^{-2})
  \right|}{s_n^2} \mathbb{E}[|\theta - \hat \theta|] 
  \\
  & \qquad + \frac{s_n^2}{\sqrt{2 n}} \mathbb{E} \bigg[
    \frac{(
    \theta-\hat \theta )^2}{|\theta|} \bigg]                                                                                            
  +  \sqrt{\frac{2}{n}} \frac{1}{s_n^2} 
 \mathbb{E} \left[ Y_1\eta'''(R_2) (
    \theta-\hat \theta )^2\right].  
\end{align*}
The moments can be obtained explicitly, and all terms are easy to
obtain once the properties of the prior are known.  For instance, if
the prior is $\eta(\theta) = \tau_1\theta -\tau_2 \log \theta$ for
$\tau_1>0$ and $\tau_2>-1$ we have 
$\mathbb{E}[|\theta^{\star}_{\mathrm{MLE}}|^2] \le {n \left(s_n^4
    (n+2 (\tau_2+1) (\tau_2+2))-4 s^2_n \tau_1
    (\tau_2+1)+2\tau_1^2\right)}/{(n s^2_n+2 \tau_1)^2}$) whence 
\begin{align*}
   &  d_{\mathrm{Wass}}(\theta^{\star}_{\mathrm{MLE}}
    , N)  \le 
     \frac{|\tau_1 - \tau_2s_n^2|}{\sqrt{n/2} s_n^2} + \frac{
      2 |\tau_2|}{n} \mathbb{E}[|\theta^{\star}_{\mathrm{MLE}}|]   + \sqrt{\frac{2}{n}}{\left|n - 2 \tau_2\right|}\frac{(n s_n^4 + 2 (s_n^4
  \tau_2 + (\tau_1 - s_n^2 \tau_2)^2))}{ s_n^2 (n s_n^2 + 2 \tau_1)(n+\tau_2)}.
\end{align*}
The leading terms in the above give the bound
\begin{align*}
  d_{\mathrm{Wass}}(\theta^{\star}_{\mathrm{MLE}}
    , N)  \le \sqrt{\frac{2}{n}}   \frac{|
  \tau_2s_n^2 - \tau_1| + s_n^2}{ s_n^2}
+ \mathrm{e}(n,
  \tau_1, \tau_2, s_n^2)
\end{align*}
with $\mathrm{e}(n, \tau_1, \tau_2, s_n^2)$ 
an (explicit) error term of lower order $O(n^{-1})$.
A lower bound can be calculated as well from \eqref{wasslow}, leading to
\begin{align*}
d_{\mathrm{Wass}}(\theta^{\star}_{\mathrm{MLE}},N) \ge    \sqrt{\frac{2}{n}}\frac{
  \vert \tau_2 s_n^2  -  \tau_1 + s_n^2 \vert }{
 s_n^2} \frac{1}{ 1 + 2 \tau_1/(ns_n^2)},
\end{align*}
which is very close to the upper bound.  Bounds on the total variation
distance follow from exactly the same argument, only with different
constants.  
The MLE depends on
\emph{all} the model parameters as well as the sample in this case,
contrarily to the posterior mode centring which led to distances
which depend only on $\tau_2$ and not on the sample. As before, all
these findings are confirmed by numerical explorations; {in particular these results
  indicate that there are parameter constellations leading to equality
  in the Wasserstein distance with the lower bound}.

 \subsection{Weibull data}
    \label{ex:weibull}
Consider i.i.d.\  Weibull data with density on the positive real line
\begin{align*}
    p_1(x \vert  \ell,m) &=  \frac{m}{\ell^m}x^{m-1} \mathrm{exp} \left\{ -(x/\ell)^m \right\}  = \mathrm{exp} \left\{ -x^m{\ell^{-m}} + m \log x - \log x +  \log m + \log (\ell^{-m}) \right\}. 
\end{align*}
We suppose $m$ is known, and $\ell$ is the parameter of interest.
Fix, for $\tau_1>0$ and $\tau_2>-1$, a (conjugate) gamma prior on $\theta=\ell^{-m}$  of the form $\pi_0(\theta)  \propto e^{-\eta(\theta)}$, with $\eta(\theta) =  \tau_1\theta 
 -\tau_2 \log \theta$ (for \(\tau_1=0\) and \(\tau_2=-1\) one obtains the improper Jeffreys prior for which our bounds also hold for $n\geq 2$). Therefore $   \lambda(\theta) = \tau_1\theta 
 -(\tau_2 + n ) \log \theta$ and we take $\Theta = (0, \infty)$.
By the same reasoning as for the previous example, Assumption A\ref{item:A2} holds if $n+\tau_2>0$. 
 Simple computations yield  $  \tilde \theta = ({1 + \tau_2/n})/{(\overline{x^m} + \tau_1/n})$ and $ K= \sqrt{n} ({\overline{x^m} + \tau_1/n})/({\sqrt{1 + \tau_2/n})}$ (here $\overline{x^m} = n^{-1} \sum_{i=1}^n x_i^m$).  The situation is very similar to the 
  previous example and we immediately get 
    \begin{align*}  
     d_\mathrm{TV}( \theta^\star_{\mathrm{MAP}}, N)  
      \le  \sqrt{\frac{\pi}{8}}\frac{1}{  \sqrt{n+\tau_2} }, \quad d_\mathrm{Wass}( \theta^\star_{\mathrm{MAP}}, N)  
      =  \frac{1}{  \sqrt{n+\tau_2} }
                \end{align*}
                (equality in the Wasserstein distance follows again because the upper and lower bounds coincide; in particular we deduce that the distributions are stochastically ordered). 

   We can also work in the MLE standardisation, i.e.\ standardisation around the MLE. That is, we set $\theta^\star_{\mathrm{MLE}} = \widehat K (\theta - \hat \theta)$ with $\hat \theta = 1/\overline{x^m}$ the MLE and $\widehat K = ({n \beta''(\hat \theta)})^{1/2} = \sqrt{n} \, \overline{x^m}$. Then keeping $\eta$ unspecified we get 
   \begin{align*}
      & \mathfrak{D}_1(\lambda, \hat \theta, \widehat K) \le |\eta'(\hat \theta)|,\quad \mathfrak{D}_2(\lambda, \hat \theta, \widehat K) \le |\eta''(\hat \theta)| \mathbb{E}[|\theta - \hat \theta|],\\
      & \mathfrak{D}_3(\lambda, \hat \theta, \widehat K) \le\frac{n}{2 \hat \theta^2} \mathbb{E}\bigg[\frac{(\theta - \hat \theta)^2}{|\theta|}\bigg]+ \mathbb{E}[|\Upsilon^3(\eta, \hat \theta, \theta^\star_{\mathrm{MLE}})| (\theta - \hat \theta)^2].
   \end{align*}
   Also, $\| F\|_\mathrm{diag} = 1/{(\sqrt{n}\, \overline{x^m}})$. 
   Taking $\eta$ the uniform (improper) prior, the first two terms disappear and the  bounds become 
   \begin{equation*}
       d_\mathrm{TV}( \theta^\star_{\mathrm{MLE}}, N)  
      \le  \sqrt{\frac{\pi}{8}}\frac{1}{\sqrt{n} } , \quad d_\mathrm{Wass}( \theta^\star_{\mathrm{MLE}}, N)  
      =  \frac{1}{  \sqrt{n} }
   \end{equation*}
   which do not depend on the sample. Finally, in case  a conjugate prior is chosen we set 
   \begin{align*}
       \Delta & = \frac{| \tau_1 - \tau_2 \overline{x^m}|}{\sqrt{n}\overline{x^m}} + \frac{\tau_2|\tau_1 - \overline{x^m}(1+ \tau_2)|}{\sqrt{n}(n \overline{x^m} + \tau_1)}  + 2 \tau_2 \frac{\Gamma(n+ \tau_2+1, n + \tau_1/ \overline{x^m})}{\sqrt n\Gamma(n+ \tau_2+1)}  \\
       & \quad  +  2 \overline{x^m}  \tau_2 \frac{n + \tau_2 +1}{n \overline{x^m} + \tau_1}\frac{\Gamma(n+ \tau_2+1, n + \tau_1/ \overline{x^m})}{\sqrt n\Gamma(n+ \tau_2+1)} + \frac{(n+\tau_2) (\overline{x^m})^2   + (\tau_1 - \tau_2 \overline{x^m})^2}{\sqrt{n} 2  \overline{x^m}  (n \overline{x^m} + \tau_1)}
   \end{align*} 
   and get
      \begin{align*}
        d_\mathrm{TV}( \theta^\star_{\mathrm{MAP}}, N)  
      \le \sqrt\frac{\pi}{8}\Delta , \quad  \sqrt{n} \frac{|\tau_1-(1+\tau_2) \overline{x^m}|}{n \overline{x^m} + \tau_1} \le d_\mathrm{Wass}( \theta^\star_{\mathrm{MAP}}, N)  
      \le \Delta,
      \end{align*}
 and thus bounds in total variation and Wasserstein distances, which, contrarily to the MAP standardisation, now depend on all parameters from the model and the data $\mathbf{x}$ in an intricate way.  

  \subsection{Multinomial data
  }
  \label{ex:multi}

Consider i.i.d.\ multinomial data with parameter $p=(p_1,\ldots,p_k)^\intercal$ and individual  likelihood 
 $   p_1(x \vert p)
 = \prod_{j=1}^kp_j^{x_j}
 $,
where $0<p_1,\ldots,p_k<1$ with $\sum_{j=1}^kp_j=1$, and $x_1,\ldots,x_k\in\{0,1\}$ with $\sum_{j=1}^kx_j=1$. 
We express the likelihood 
as an exponential family by writing 
\begin{align*}
    p_1(x \vert p) &=\mathrm{exp}\bigg\{ \sum_{j=1}^{k-1} x_j\log(p_j)+\bigg(1-\sum_{j=1}^{k-1}x_j\bigg)\log(p_k)\bigg\}=\mathrm{exp}\bigg\{ \sum_{j=1}^{k-1} x_j\log\left(\frac{p_j}{p_k}\right)+\log(p_k)\bigg\},
\end{align*}
so $h(x_j) = x_j$, $g(x) = 0$ and the natural parameter is
$\theta=(\theta_1,\ldots,\theta_{k-1})^\intercal$, where
$\theta_j = \log (p_j/p_k)$, $j=1,\ldots,k-1$, and
$\Theta = \mathbb{R}^{k-1}$ (thus $d = k-1$ in the notations of the
previous section). Expressing the density in terms of $\theta$, we get
$p_j(\theta) = p_k e^{\theta_j}$. Also, since $\sum_{j=1}^kp_j=1$, we
obtain that $p_k+p_k\sum_{j=1}^{k-1}e^{\theta_j}=1$, so that
$p_k=(1+\sum_{j=1}^{k-1}e^{\theta_j})^{-1}$.  Thus,
$ \beta(\theta) = -\log(p_k(\theta)) =  \log(1+ \sum_{j=1}^{k-1}e^{\theta_j}).
$
We fix a {conjugate Dirichlet prior} on
$p=(p_1,\ldots,p_{k-1})^\intercal$ of the form
$\pi_0(p) \propto (1 - \sum_{j=1}^{k-1} p_j)^{\tau_k}\prod_{j=1}^{k-1}
p_j^{\tau_j} $ for some parameters $\tau_1,\ldots,\tau_k>-1$.  Let
$\tilde \tau = \sum_{j=1}^k \tau_j$.  Expressing everything in terms
of $\theta$ we get
$\pi_0(\theta) \propto \mathrm{exp}({-\eta(\theta)})$, where
$\eta(\theta) = -\sum_{j=1}^{k-1}\tau_j\theta_j + \tilde \tau \log
\big(1+ \sum_{j=1}^{k-1}e^{\theta_j} \big).$ Therefore
$\lambda(\theta) = -\sum_{j=1}^{k-1}\tau_j \theta_j + (n +\tilde \tau)
\log \big( 1+ \sum_{j=1}^{k-1}e^{\theta_j} \big)$, and, letting
$\pi_u = \tau_u + n \bar x_u$ with  $\barx_u =  n^{-1}\sum_{i=1}^n
x_{i, u}$, $u = 1, \ldots, k$, the posterior is
easily seen to be of the form
\begin{align}\label{eq:p2dir}
  p_2(\theta \, | \, \bm\pi) = \frac{\Gamma(n+\tilde\tau)}{ \prod_{j=1}^k
  \Gamma(\pi_j)} \frac{\mathrm{exp}( \sum_{j=1}^{k-1} \pi_j
  \theta_j)}{(1 + \sum_{j=1}^{k-1} e^{\theta_j})^{n+\tilde\tau}}
\end{align}
with $\bm\pi = (\pi_1, \ldots, \pi_{k-1})^{\intercal}$ (the
normalising constant is known because this is, after a change of
variables, a Dirichlet distribution). Note how
$\sum_{u=1}^{k} \pi_u = n + \tilde \tau$.

We  standardise  around the posterior mode i.e. $\tilde \theta$ that
solves the system of equations 
$\partial_{u}\lambda(\theta)- n \bar x_u=0$, for
$u=1,\ldots,k-1$.  
   A simple calculation reveals that
      $\tilde\theta=(\tilde\theta_1,\ldots,\tilde\theta_{k-1})^\intercal$
      is given by
\begin{equation*}
  \tilde\theta_u = \log \bigg(\frac{\tau_u+ n  \barx_{u}}{{\tau}_k +
    n \barx_k}\bigg) =  \log \Big( \frac{\pi_u}{\pi_k}\Big), \quad u=1,\ldots,k-1.
\end{equation*}
The second order partial derivatives of $\lambda$ are given by
\begin{align*}
    \partial_{ u, u}\lambda(\theta)=(n+\tilde \tau)\frac{e^{\theta_u}(1+\sum_{\ell=1}^{k-1}e^{\theta_\ell}-e^{\theta_u})}{(1+\sum_{\ell=1}^{k-1}e^{\theta_\ell})^2},\quad
    \partial_{u, v}\lambda(\theta)=-(n+\tilde \tau)\frac{e^{\theta_u+\theta_v}}{(1+\sum_{\ell=1}^{k-1}e^{\theta_\ell})^2}, \quad u\not=v,
\end{align*}
from which it follows that 
\begin{align}
   \label{hes12}  \partial_{u, u}\lambda(\tilde\theta) = \pi_u - \frac{1}{n+\tilde \tau} \pi_u^2, \quad
 \partial_{u,
  v}\lambda(\tilde\theta)=-\frac{1}{n+\tilde \tau}\pi_u \pi_v, \quad u\not=v. 
\end{align}
The entries in \eqref{hes12}
correspond to the
elements in the $(k-1)\times (k-1)$ Hessian matrix
$\widetilde{\matr H}:= \matr{H}_\lambda(\tilde\theta)$, which can therefore be rewritten
as
\begin{equation*}
\widetilde{\matr H} = \mathrm{diag}(\bm{\pi}) -
(n+\tilde \tau)^{-1} \bm\pi \bm \pi^\intercal, 
\end{equation*}
where $ \bm \pi$ is the $(k-1) \times 1$ vector with entries
$\pi_u = \tau_u+  n \barx_u$, $u = 1, \ldots, k-1$. 
Standard arguments (the Sherman-Morrisson formula)   
give 
\begin{align*}
 \widetilde{\matr F} :=  \mathbf{F}_\lambda(\tilde\theta) & =  \mathrm{diag}(\frac{1}{\bm{\pi}}) +
                                     \frac{ 1}{\pi_k}  \mathbf 1
                                     \mathbf
                                     1^\intercal,
\end{align*}
where $ \mathrm{diag}({1}/{\bm{\pi}})$ is the diagonal matrix with
diagonal entries $1/\pi_u$, $u = 1, \ldots, k-1$ and $\mathbf 1$ is the
$(k-1) \times 1$ vector of 1's.  Also, the matrices
$\widetilde{\matr H}$ and $\widetilde{\matr F}$ admit unique
symmetric square roots
$\widetilde{\matr K}$ and $\widetilde{\matr J}$,
respectively, which are, moreover, computable explicitly 
(diagonal-plus-rank-one matrices are well-studied, see e.g.\
\cite{bunch1978rank, stor2015forward}). 

With these preliminaries, we are ready to tackle the bound from
Theorem~\ref{thm:abstract-results}. 
From Example~\ref{sec:bounds-map} we know that all we need is to
control two terms. For the first, we  
have by definition, 
  \begin{align}\label{eq:goodboundonJgrad}
    \| \widetilde{\matr F}\|_{\mathrm{diag}} = \sqrt{\max_{1 \le u \le k-1}\big|  \widetilde{\matr F}_{uu} \big|} \le  \max_{1 \le u \le k-1}  {\sqrt{\frac{1}{\pi_u} + \frac{1}{\pi_k} }}
    \le  \sqrt{\frac{2}{\min \bm \pi}}.
  \end{align}
Next, for $\mathfrak D_3(\lambda, \tilde \theta, \widetilde{\matr K})$
the term depending on the third derivatives of $\lambda$, we start
by noting that
\begin{align*}
\partial_{u, u, u}\lambda(t)&=
                       (n+ \tilde \tau)     \frac{e^{t_u}(1+\sum_{\ell=1}^{k-1}e^{t_\ell}-e^{t_u})
                            (1+\sum_{\ell=1}^{k-1}e^{t_\ell}-2e^{t_u}
                            )}{(1+\sum_{\ell=1}^{k-1}e^{t_\ell})^3}
                            \mbox{ for all } 1 \le u \le k-1,
                          \\
\partial_{u, u, v}\lambda(t)&= -                        (n+ \tilde \tau)\frac{e^{t_u+t_v} (
                            1+\sum_{\ell=1}^{k-1}e^{t_\ell} - 2e^{t_u}
                            )}{(1+\sum_{\ell=1}^{k-1}e^{t_\ell})^3}
\mbox{ for all } 1 \le u\not=v \le k-1, \\
\partial_{u, v, w}\lambda(t)&=                        (n+ \tilde \tau) \frac{2e^{t_u+t_v+t_w}
                            }{(1+\sum_{\ell=1}^{k-1}e^{t_\ell})^3} 
\mbox{ for all } 1 \le u\not=v\not=w \le 1.
\end{align*}
Let
\begin{equation}
  \label{eq:8}
  p_u(t) =  \frac{e^{t_u}}{1+ \sum_{\ell = 1}^{k-1} e^{t_{\ell}}}.
\end{equation}
Then, since $0 \le p_u(t) \le 1$  and $0 \le \sum_{u=1}^{k-1} p_u(t) \le 1$ for all $t$,
clearly
\begin{align*}
&   |\partial_{u, u, u}\lambda(t)|  \le (n+\tilde\tau) p_u(t),\ \quad |\partial_{u,
  u, v}\lambda(t)| \le (n+\tilde\tau)  p_u(t) p_v(t) \\
  & \mbox{ and } |\partial_{u, v,
  w}\lambda(t)|  \le  2 (n+\tilde\tau) p_u(t)p_v(t)p_w(t)
\end{align*}
for all $1 \le u \neq v \neq w \le k-1$.  Assumption
\hyperref[assumptionP2]{P1} is therefore satisfied on the entire
parameter space with
$C_1 = C_2 = (n+ \tilde \tau)$ and $C_3 = 2(n+ \tilde \tau)$, as well
as $C_{\infty}(t) \le 4(n+ \tilde \tau)$ (for all $t$) and we can
apply Lemma \ref{sec:main-results} to get
\begin{align}\label{eq:y2sumputh}
 \Delta_{3}(\lambda,
    \tilde \theta, \widetilde{\matr K} ) =|
  \mathfrak{D}_{3}(\lambda, \tilde \theta,\widetilde{\matr K} )  | \le 
                                               4(n+ \tilde \tau)    \mathbb{E} \bigg[
                                                   Y_1
                                                   \sum_{u=1}^{k-1} p_u(R_2) (\theta_u - \tilde\theta_u)^2                                                    \bigg]   
\end{align}
(as before $R_2 = \tilde \theta + Y_1Y_2(\theta - \tilde \theta)$). 
A proof 
of the following lemma is given in
Section \ref{sec:about-dichl-distr} of the Appendix. 
\begin{lemma}\label{sec:multnomoments}
  Under the previous notations, and taking without loss of generality that $\pi_k=\max\bm\pi$, we have, for $u=1,\ldots,k-1$,
\begin{align*}
\mathbb{E}\big[(\theta_u-\tilde\theta_u)^2\big]&\leq \frac{2}{\pi_u}\bigg(1+\frac{3}{2\pi_u}\bigg),\quad 
\mathbb{E}\big[(\theta_u-\tilde\theta_u)^4\big]\leq\frac{12}{\pi_u^2}\bigg(1+\frac{11}{3\pi_u}+\frac{15}{4\pi_u^2}\bigg).
\end{align*}  
\end{lemma}
                                                 
Explicit bounds on the total variation and Wasserstein
distances now follow easily. We  apply \eqref{eq:boundwass} and
\eqref{eq:bounTV} with the simplification 
$p_u(R_2) \in (0, 1)$ to get 
\begin{align*}
  \Delta_{3}(\lambda, \tilde \theta, \widetilde{\matr K}) &  \le 2(n+ \tilde \tau)  \sum_{u=1}^{k-1}   \mathbb{E} \left[
                                                    (\theta_u - \tilde
                                                              \theta_u)^2
                                                \right]  \le  2(n+ \tilde
                                                  \tau)\frac{k}{\min\bm\pi}\bigg(1+\frac{3}{2\min\bm\pi}\bigg)
                                                  .
\end{align*}
Hence 
\begin{align*}
   d_{\mathrm{Wass}}(\theta^{\star}_{\mathrm{MAP}}, N)  &\le 4\sqrt{2}(n+\tilde\tau)\bigg(\frac{k}{\min\bm\pi}\bigg)^{3/2}\bigg(1+\frac{3}{2\min\bm\pi}\bigg),\\
   d_{\mathrm{TV}}(\theta^{\star}_{\mathrm{MAP}}, N) & \le 4\sqrt{\pi}(n+\tilde\tau)\frac{k}{(\min\bm\pi)^{3/2}}\bigg(1+\frac{3}{2\min\bm\pi}\bigg).
\end{align*}
The Wasserstein and total variation bounds are of orders
$O((k^6/n)^{1/2})$ and $O((k^5/n)^{1/2})$, respectively, if $\min \bm \pi \approx n/k$ (
consistency is proved if $k^6 \log k / n  \to 0$ in \cite{ghosal2000asymptotic}, and if $k^{4}/n \to 0$ in
\cite{belloni2014posterior}). 

The dependence on the dimension can be reduced in the Wasserstein and total variation distances at
the cost of some more work. Here, for the sake of simplicity, we work only up to
an absolute constant, but in Appendix \ref{sec:about-dichl-distr} we derive bounds with the constants as stated in Example \ref{ex:multiintro}. First, we note that, by definition,
\begin{align}\label{max}
  p_u (\tilde \theta) = \frac{\pi_u}{n+ \tilde \tau}, \quad u = 1, \ldots, k-1.
\end{align}
Observe that (\ref{max}) is of order $O(1/k)$ if $\max\bm\pi\approx n/k$, which is exactly what is needed to reduce the
dimensional dependence of our bounds.  Following \cite[Equation
(3.4)]{ghosal2000asymptotic}, we thus expand $p_u(t)$ around the
posterior mode $\tilde \theta$ to get, for all $u = 1, \ldots, k$,
\begin{align*}
   p_u(R_2)  &\le p_u (\tilde
             \theta)  + | p_u(R_2) -  p_u(\tilde \theta)|  \nonumber \\
           & \le \frac{\pi_u}{n + \tilde \tau} + \frac{| e^{R_{2, u}} -
             e^{\tilde \theta_u}|}{1+ \sum_{\ell=1}^{k-1}e^{R_{2,\ell}}} +
 \frac{\pi_u}{\pi_k}   \sum_{v=1}^{k-1} \frac{ |e^{R_{2, v}} -
             e^{\tilde \theta_v}|}{(1+ \sum_{\ell=1}^{k-1}e^{R_{2,
             \ell}}) (1 + \sum_{\ell = 1}^{k-1} e^{\tilde
             \theta_{\ell}})} \nonumber\\
  & \le \frac{\pi_u}{n + \tilde \tau} +  \frac{\pi_uM(\theta)}{n+\tilde\tau} {|\theta_u - \tilde \theta_u |}
  +  \frac{\pi_u M(\theta)}{(n+\tilde\tau)^2} \sum_{v=1}^{k-1} \pi_v|\theta_v - \tilde \theta_v |,
  \nonumber 
\end{align*}
where $M(\theta)=\max_{1\leq \ell\leq k-1}e^{2|\theta_\ell-\tilde\theta_\ell|}$. Here we used the basic inequality $|e^x-e^y|\leq |x-y|e^ye^{|x-y|}$, from which we obtain that $|e^{R_{2,u}}-e^{\tilde\theta_u}|\leq |\theta_u-\tilde\theta_u|e^{\tilde \theta_u}e^{|\theta_u-\tilde\theta_u|}=(\pi_u/\pi_k)|\theta_u-\tilde\theta_u|e^{|\theta_u-\tilde\theta_u|}$ We also used that $1/(1+\sum_{\ell=1}^{k-1}e^{\tilde\theta_\ell})\leq (\pi_k/(n+\tilde\tau))\max_{1\leq \ell\leq k}e^{|\theta_\ell-\tilde\theta_\ell|}$, which follows from similar considerations along with the fact that $\sum_{\ell=1}^k\pi_\ell=n+\tilde\tau$.
From         \eqref{eq:y2sumputh} we get 
\begin{align*}
 \Delta_3(\lambda, \tilde \theta, \widetilde{\matr K})  &
                                                                    \le  4 (n+ \tilde \tau)  \sum_{u=1}^{k-1} \mathbb{E} \left[ Y_1 p_u(R_2)  (\theta_u-
    \tilde
    \theta_u)^2   \right] \le I_1 + I_2 + I_3\end{align*}
with 
\begin{align*}
  I_1&    
= 2
 \sum_{u=1}^{k-1}  \pi_u    \mathbb{E} \left[
       (\theta_u-
       \tilde
       \theta_u)^2  \right] ,  \quad
  I_{2}  
=     2 
                     \sum_{u=1}^{k-1}  
         \pi_u \mathbb{E} \left[ |\theta_u - \tilde \theta_u|^3 M(\theta)\right],\\
  I_{3}       &  
 =  \frac{2}{n+ \tilde \tau} 
    \sum_{u=1}^{k-1}    \sum_{v=1}^{k-1}
    \pi_u\pi_v
     \mathbb{E} \left[ |\theta_v - \tilde \theta_v||\theta_u - \tilde
    \theta_u|^2 M(\theta)\right]. 
\end{align*}
With H\"older's inequality we obtain
\begin{align}
  &         \mathbb{E} \big[
    |\theta_u-
    \tilde
    \theta_u|^3M(\theta) \big]   \le \big(\mathbb{E}\big[(\theta_u-\tilde\theta_u)^4\big]\big)^{3/4}\big(\mathbb{E}\big[(M(\theta))^4\big]\big)^{1/4},\label{cs1} \\
  &   
   \mathbb{E} \big[ |\theta_v - \tilde \theta_v||\theta_u - \tilde
    \theta_u|^2M(\theta)\big]  \le  \big(\mathbb{E}\big[(\theta_v-\tilde\theta_v)^4\big]\mathbb{E}\big[(\theta_u-\tilde\theta_u)^4\big]\big)^{3/4}\big(\mathbb{E}\big[(M(\theta))^4\big]\big)^{1/4}. \label{cs2} 
\end{align}
Note that $\mathbb{E}[(M(\theta))^4]<\infty$ if $\pi_k>8$ (this condition is easily seen from formula (\ref{eq:p2dir}) for the posterior density). Recalling that $e^x=1+O(x)$ as $x\downarrow0$ and using the crude bound that $\max_{1\leq i\leq k-1}y_i\leq\sum_{i=1}^{k-1}|y_i|$ we also have that 
$\mathbb{E}[(M(\theta))^4]=1+O(k(\min \bm\pi)^{-1/2})$. Therefore applying the approximations of Lemma \ref{sec:multnomoments} gives that
\begin{align*}
  I_1 \lesssim 
k, \;  I_2 \lesssim   \frac{k}{\sqrt{\min \bm \pi}}
\mbox{
  and }  I_3 \lesssim  \frac{k^2}{n+\tilde\tau}\sqrt{\max \bm\pi}. \end{align*}
From the basic inequalities $\max\bm\pi\leq\sum_{u=1}^k\pi_u= n+\tilde\tau$ and $n+\tilde\tau=\sum_{u=1}^k\pi_u\geq k\min\bm\pi$ we obtain the bound $k\sqrt{\max\bm\pi}/(n+\tilde\tau)\leq \sqrt{k/\min\bm\pi}$. Thus, if we assume that $\min\bm\pi\geq ck$ for some constant $c>0$, then we get from \eqref{eq:boundwass}
and \eqref{eq:bounTV} that there exist positive constants $\mathcal{C}_1$ and $\mathcal{C}_2$ such that  
\begin{align}
d_{\mathrm{TV}}(\theta^{\star}_{\mathrm{MAP}}, N) \leq \frac{\mathcal{C}_1k}{\sqrt{\min \bm \pi}}, \quad d_{\mathrm{Wass}}(\theta^{\star}_{\mathrm{MAP}}, N) \leq\frac{\mathcal{C}_2k^{3/2}}{\sqrt{\min \bm \pi}}.
\label{our_bounds}
\end{align}
 
\begin{remark} In Appendix \ref{sec:about-dichl-distr}, we derive an estimate for the constant $\mathcal{C}_2$ under the assumption that $\sqrt{\min\bm\pi}\geq 7.40k/\sqrt{2}$. (We can derive an estimate for $\mathcal{C}_1$ with no need for assumptions on $\min\bm\pi$ thanks to the trivial bound $d_{\mathrm{TV}}(\theta^{\star}_{\mathrm{MAP}}, N) \leq1$.) This assumption is mild (the bound will not converge to zero if this condition does not hold), but somewhat arbitrary. However, in this remark, we provide some motivation for our choice. If this assumption does not hold, then we have that $7.40k^{3/2}/\sqrt{\min\bm\pi}\geq \sqrt{2k}$, in which case we can obtain the following alternative estimate.   Recalling that $
  d_{\mathrm{Wass}}(\theta_{\mathrm{MAP}}^*, N) = \sup_{\phi \in \mathfrak{F}_{\mathrm{Wass}}} | \mathbb{E} [\phi(\theta_{\mathrm{MAP}}^*)]
  - \mathbb{E} [\phi(N)]| $, where $\mathfrak{F}_{\mathrm{Wass}}$ is the class of Lipschitz functions on $\mathbb{R}^d$ such that $\|\phi\|_{\mathrm{Lip}}:=\sup_{s\not=t}|\phi(s)-\phi(t)|/\|s-t\|_2 \le 1$, we obtain that
\begin{align*}
d_{\mathrm{Wass}}(\theta_{\mathrm{MAP}}^*, N)&\leq\mathbb{E}[\|\theta_{\mathrm{MAP}}^*-N\|_2]\leq\sqrt{\sum_{u=1}^{k-1}\mathbb{E}[(\theta_u^*-N_u)^2]}=\sqrt{\sum_{u=1}^{k-1}(\mathbb{E}[(\theta_u^*)^2]+1)}=\sqrt{2(k-1)}+r_n,  
\end{align*} 
where $r_n$ is a remainder term of of order $o(1)$ as $n\rightarrow\infty$ and $\theta_u^*=(\theta_{\mathrm{MAP}}^*)_u$, $1\leq u\leq k-1$. That $r_n=o(1)$  as $n\rightarrow\infty$ follows since $\mathbb{E}[(\theta_u^*)^2]=\mathbb{E}[N_u^2]+o(1)=1+o(1)$, where $N_u\sim N(0,1)$.  
So also in this case we obtain a bound that is $O(k^{1/2})$.
\end{remark}

Finally, we note that the entire procedure can be applied to the MLE standardisation, irrespective of the choice of prior. We have done the computations but the bounds take a rather unsavory form and do not lead to any further improved  intuition beyond  the fact that, under mild conditions on the prior, $d_{\mathrm{TV}}(\theta_{\mathrm{MLE}}^*, N)$ and $d_{\mathrm{Wass}}(\theta_{\mathrm{MLE}}^*, N)$ converge to zero if $k/\sqrt{n\min\bar{\bm x}}\rightarrow0$ and $k^{3/2}/\sqrt{n\min\bar{\bm x}}\rightarrow0$, respectively, where $\min\bar{\bm x}=\min_{1\leq u\leq k}\bar{x}_u$. 
We dispense with the details here and leave the derivation of the bounds to the motivated reader.

\begin{remark}[About the literature]
  As mentioned in the introduction, BvM theorems for multinomial data have been extensively studied, under different parameterisations. As far as we are aware, the closest results to ours are to be found in \cite{ouimet2022multivariate}  and \cite{katsevitchimprovedevenbetter}. The result from   \cite[Theorem 2]{ouimet2022multivariate} reads,  in our notation, as 
\begin{align}
  \label{eq:3}
  d_{\mathrm{TV}}( X_n, N_n) \lesssim \frac{k-1}{\sqrt{n + \tilde \tau +  1}} \sqrt{\frac{\max \bm \pi}{\min \bm\pi}} ,
\end{align}
where $X_n \sim \mathrm{Dir}(\bm\pi)$ and $N_n \sim \mathcal{N}(\mu_n, \bm \Sigma_n)$ with
$\mu = {\bm\pi}/(n+\tilde \tau)$ and
$\bm \Sigma_n = (n+ \tilde \tau+1)^{-1}(n + \tilde \tau)^{-1} \left(
  \mathrm{diag}(\bm \pi) -  \bm \pi \bm
  \pi^{\intercal} \right). $ 
  \cite{katsevitchimprovedevenbetter} considers the exact same setting (with a flat prior $\bm\tau=  \bm1$) and  obtains under some mild conditions 
\begin{equation} \label{eq:33}
   \frac{1}{18} \sum_{j=1}^k \Big| \frac{1}{k}- {\bar x_j}\Big|  \frac{(k-1)}{\sqrt{n}} \le   d_{\mathrm{TV}}(X_n,  N_n)   \lesssim \sqrt{\frac{k}{\min \bm{\bm \pi}}} + \frac{k^2}{\min{\bm \pi}} + r_n 
\end{equation}
($r_n$ is a negligible remainder term).  Our total variation bound in \eqref{our_bounds} requires the same condition to go to zero as  \eqref{eq:3} and  \eqref{eq:33}, and the lower bound in \eqref{eq:33} indicates that the condition $k /\sqrt n \to 0$ is necessary in the setting of \cite{katsevitchimprovedevenbetter,ouimet2022multivariate}.  It may be worth noting that \cite{katsevitchimprovedevenbetter} also upper bound the difference of means in this case, and obtain an extra $\sqrt k$ in their bound exactly as we do in our Wasserstein distance bounds.  To the best of our understanding, the conclusions from \eqref{eq:3} and  
 \eqref{eq:33} do not carry through to our setting (we are considering normal approximation for a non-trivial transformation of a Dirichlet distribution).  
\end{remark}

 \subsection{Univariate normal with unknown mean and
  variance}\label{ex:normal2}

Consider i.i.d.\  normal data with mean $\mu$ and variance $\sigma^2$.
Here $\beta(\theta) = 1/2 \log (2\pi /\theta_2 )+\theta_1^2/(2\theta_2)$, $g(x) = 0$, $h_1(x)=x$, $h_2(x)=-x^2/2$ and the natural parameter is   $\theta=(\theta_1,\theta_2)^\intercal=(\mu/\sigma^2,1/\sigma^2)^\intercal$, which lives on $\mathbb{R} \times (0,\infty)$.

We consider {normal approximation of the posterior standardised
  around the posterior mode}.  Fix, for
$\tau_1, \tau_2, \tau_3, \tau_4 \in \mathbb{R}$ with
$\tau_4>0, \tau_1<3/2, 4\tau_4 \tau_2 > \tau_3^2$, a (conjugate) prior
on $\theta$ of the form $\pi_0(\theta) \propto \exp({-\eta(\theta)})$,
where
$\eta(\theta) = {\log(2\pi/\theta_2) \tau_1 + \theta_1^2
  \tau_2/\theta_2 - \theta_1 \tau_3 + \theta_2 \tau_4}$.  Therefore
\begin{equation*}
        \lambda(\theta)  =  \frac{(n/2+\tau_2)\theta_1^2}{\theta_2} - \theta_1 \tau_3 + \theta_2 \tau_4 
 + \bigg( \frac{n}{2} + \tau_1 \bigg) \log (2\pi /\theta_2 ).
\end{equation*}
Letting $\tilde \tau_1 = 1+2\tau_1/n$, $\tilde \tau_2 = 1+2\tau_2/n$,
$\tilde \tau_3 = \tau_3/n$ and $\tilde \tau_4 = 2\tau_4/n$ it follows from a simple calculation
that the posterior mode $\tilde \theta$ is given by
\begin{align*}
\tilde{\theta}_1= \frac{\tilde \tau_1 (\bar x +
                 \tilde \tau_3)}{\tilde \tau_1 (\overline{x^2}+\tilde
                 \tau_4)-(\bar x +\tilde \tau_3)^2} = \frac{\bar x}{s^2_x}+O(n^{-1}),
                 \quad \tilde{\theta}_2=
                 \frac{\tilde \tau_1 \tilde \tau_2}{\tilde \tau_2(
                 \overline{x^2}+\tilde\tau_4)-( \bar x+\tilde
                 \tau_3)^2}= \frac{1}{s^2_x}+O(n^{-1}),
\end{align*}
with $s^2_x = \overline{x^{2}} - (\bar x)^2$.
We also have that
\begin{align*}
  \widetilde{\matr H} :=  \mathbf{H}_\lambda(\theta)=n \begin{pmatrix}
      \displaystyle\frac{\tilde \tau_2}{\theta_2} &\displaystyle
      -\tilde \tau_2  \frac{\theta_1}{\theta_2^2} \\[4mm]\displaystyle
      -\tilde \tau_2 \frac{\theta_1}{\theta_2^2}  & \displaystyle
      \tilde \tau_2  \frac{\theta_1^2}{\theta_2^3} + \frac{\tilde
        \tau_1}{2} \frac{1}{\theta_2^2} \end{pmatrix}
\mbox{ and } \widetilde{\matr F} := \mathbf{F}_{\lambda}( \theta)     = n^{-1}
                                                    \begin{pmatrix}
                                                     \displaystyle 2 \frac{
                                                        \theta_1^2}{\tilde\tau_1}
                                                      +
                                                      \frac{\theta_2}{\tilde\tau_2}
                                                      & \displaystyle
                                                      \frac{2}{\tilde\tau_1}
                                                      \theta_1
                                                      \theta_2\\[4mm]
\displaystyle \frac{2}{\tilde\tau_1}
                                                      \theta_1
                                                      \theta_2 & \displaystyle
                                                      \frac{2}{\tilde
                                                        \tau_1} \theta_2^2
                                                    \end{pmatrix}.
\end{align*}
The expressions for
$\widetilde{\matr K} := \mathbf{K}_\lambda(\tilde \theta)$ and
$\widetilde{\matr J}:= \mathbf{J}_\lambda(\tilde \theta)$ are quite
inelegant but perfectly explicit and can easily be
computed. 
In particular
\begin{equation*}
  \| \widetilde{\matr F}\|_{\mathrm{diag}}^2  =
\frac{1}{n}\max \bigg( 2 \frac{\tilde \theta_1^2}{\tilde \tau_1} +
  \frac{\tilde \theta_2}{\tilde \tau_2}, 2\frac{\tilde
    \theta_2^2}{\tilde \tau_1} \bigg)
\end{equation*}
which, as we shall soon see, provides the rate of convergence.

As in the previous example,
the main work is to bound 
\begin{align*}
  \mathfrak{D}_{3}(\lambda, \tilde \theta, \tilde K \, | \, \Theta_0^{\star})               &
                                                                                 = \sum_{u=1}^2
                                                                                 \sum_{v=1}^2
                                                                                 \sum_{w=1}^2
                                                                                 \mathbb{E}
                                                                                 \left[
                                                                                 Y_1
                                                                                 |\partial^3_{u, v, w} \lambda(
                                                                                 R_2)|
                                                                      \, 
                                                                                 |\theta_v
                                                                                 -
                                                                                 \tilde\theta_v|\,
                                                                                 |\theta_w
                                                                                 -
                                                                                 \tilde\theta_w| \,
                                                                                 \mathbb{I}_{\Theta_0}(\theta)
                                                                \right]
\end{align*}
(as usual $R_2 := \tilde \theta + Y_1 Y_2 (\theta - \tilde \theta)$ and
$\Theta_0 = \tilde \theta  + \widetilde{\matr J}\Theta_0^{\star}$).
Simple computations yield
\begin{align*}
  \partial_{1,1,1}\lambda(w)  & = 0, \quad
                                     \partial_{2,2,2}\lambda(w)
                                     =
                                     -\frac{3(n+2\tau_2)w_1^2}{w_2^4}-\frac{n+2\tau_1}{w_2^3},
  \\
\partial_{1, 1, 2}\lambda(w) & = \partial_{1, 2, 1}\lambda(w)
                                   = \partial_{2, 1, 1}\lambda(w)= -\frac{n+2\tau_2}{w_2^2}, \\ 
\partial_{1,2,2}\lambda(w)  & = \partial_{2,1,2}\lambda(w) =\partial_{2,2,1}\lambda(w) = \frac{2(n+2\tau_2)w_1}{w_2^3} 
. 
\end{align*}
As we cannot uniformly bound these third order derivatives,  and thus 
cannot
use \eqref{eq:boundwass} and \eqref{eq:bounTV}, we 
use
\eqref{eq:TVbounnnnthetastarcond}. An application of the basic
inequalities $|w_1|\leq(1+w_1^2)/2$ and
$|w_2|^{-3}\leq(w_2^{-2}+w_2^{-4})/2$ shows that the absolute values
of each of the third order partial derivatives are bounded above by
\begin{align*}
(n+2(\tau_1+\tau_2))(1+3w_1^2)\bigg(\frac{1}{w_2^2}+\frac{1}{w_2^4}\bigg)
  =: (n+2(\tau_1+\tau_2))T(w).
\end{align*}
We apply
Lemma~\ref{sec:main-resultsv2} (with $S(w)=2(n+2(\tau_1+\tau_2))T(w)$) on $\Theta_0^{\star}$ where
$\Theta_0 = \{ (\theta_1- \tilde \theta_1)^2 + (\theta_2- \tilde
  \theta_2)^2 \le 1\}$, leading to
\begin{align*}
  \mathfrak{D}_{3}(\lambda, \tilde \theta, \tilde K \, | \, \Theta_0^{\star})     
  & \le   4(n+2(\tau_1 + \tau_2))
    \mathbb{E}\bigg[                                   T(R_2)  \sum_{u=1}^{2}
    (\theta_u- \tilde \theta_u)^2 \mathbb{I}_{\Theta_0}(\theta)
    \bigg].
\end{align*}
Since $\tilde \theta_1>0$ and $\tilde \theta_2>0$, it holds that, for
$\theta\in \Theta_0 = \{ (\theta_1- \tilde \theta_1)^2 +
  (\theta_2- \tilde \theta_2)^2 \le 1\}$, we have
\begin{align*}
  T(\tilde\theta + Y_1Y_2(\theta - \tilde \theta)) \le \left( 1 + 6
  \left( \tilde \theta_1^2 + 1 \right) \right)
  \bigg(\frac{1}{\tilde \theta_2^2}+\frac{1}{\tilde \theta_2^4}\bigg)
  =: T^{\star}(\tilde \theta),
\end{align*}
so that 
\begin{align*}
 \mathfrak{D}_{3}(\lambda, \tilde \theta, \tilde K \, | \, \Theta_0^{\star})  \le 4(n + 2
  (\tau_1+ \tau_2)) T^{\star}(\tilde \theta)   \sum_{u=1}^{2}   \mathbb{E}[    
  (\theta_u- \tilde \theta_u)^2                            
  ]. 
\end{align*}
Also,
$p_2(\theta \, | \, \mathbf{x})$ is proportional to
\begin{equation*}   \mathrm{exp} \left( (n \bar x + \tau_3) \theta_1 - (n
                 \overline{x^2}/2 + \tau_4) \theta_2  -
                 \frac{(n/2+\tau_2)\theta_1^2}{\theta_2}  
 + \bigg( \frac{n}{2} + \tau_1 \bigg) \log (\theta_2 ) \right).
\end{equation*}
Hence the term $r^{\star}(\Theta_0^{\star})$ is negligible. Finally,
$ \mathbb{P}[\theta \notin \Theta_0^{\star}] =
\mathbb{P}[\sum_{u=1}^2(\theta_u - \tilde \theta_u)^2\ge 1]
\le\sum_{u=1}^2 \mathbb{E}[(\theta_u - \tilde \theta_u)^2]$ so that it
only remains to compute
$\sum_{u=1}^2\mathbb{E}[(\theta_u - \tilde \theta_u)^2] =
\sum_{u=1}^2 \{\mathrm{Var}(\theta_u) + (\mathbb{E}[\theta_u] - \tilde
\theta_u)^2\}.$ We can either use properties of the exponential family
or direct computations from the above to reap explicit expressions for
all moments involved: 
\begin{align*}
  \mathbb{E}[\theta_1] & = -\frac{(n+2 \tau_1+3) (n \bar x + \tau_3)}{\bar{x}^2 n^2+2
                         \bar{x} n \tau_3-(n+2
                         \tau_2) (n \overline{x^2}+2
                         \tau_4)+\tau_3^2}=\frac{\bar x}{s_x^2}+O(n^{-1}), \\
  \mathbb{E}[\theta_2] & = -\frac{(n+2 \tau_1+3) (n+2
                         \tau_2)}{\bar{x}^2 n^2+2
                         \bar{x} n \tau_3-(n+2
                         \tau_2) (n \overline{x^2}+2
                         \tau_4)+\tau_3^2}=\frac{1}{s_x^2}+O(n^{-1}), 
\end{align*}
and similarly 
\begin{align*}
    \mathrm{Var}(\theta_1) & 
                         =\frac1n\frac{\bar{x}^2 + \overline{x^2}}{(s^2_x)^2} + O(n^{-2}), \quad 
\mathrm{Var}(\theta_2) {=\frac1n\frac{2}{(s^2_x)^2}+O(n^{-2})}.
\end{align*}
Wrapping up, on ignoring 
the lower order terms
which we collect into a remainder term $r_n = O(n^{-3/2})$ 
we get 
\begin{align}\nonumber 
   d_{\mathrm{TV}}(\theta^{\star}_{\mathrm{MAP}}, N) &\leq  \frac{28 \sqrt{\pi}}{\sqrt{n} } \frac{\sqrt{1+(\bar x)^2}}{s^2_x}  \left( 1+ \frac{(\bar x)^2 }{ (s^2_x)^2}\right)  (1+(s^2_x)^2)  ( \bar x^2 + \overline{x^2} +2)   +\frac{\bar x^2 + \overline{x^2} +2}{n(s_x^2)^2} +  r_n.
\end{align}
To leading order the bound only depends on the sample size and the data; the prior hyperparameters only appear in the lower order remainder term $r_n$.

  
\subsection{Linear regression with unknown
  variance}\label{ex:linreg}
We conclude the paper with a non-identically distributed example.
Consider $X_1, \ldots, X_n $ from a standard linear regression model
with unknown variance, i.e.\
\begin{align*}
    X=M\xi+ \epsilon,
\end{align*}
where $M \in \mathbb{R}^{n \times k}$ is the design matrix, $\xi \in \mathbb{R}^k$ and $\epsilon \sim N(0,\sigma^2 I_n)$. We will make a standard assumption that $M^\intercal M$ is invertible. 
We have that $X_i\sim N(\mathbf{m}_i \xi,\sigma^2)$,
where $\mathbf{m}_i = (m_{i1}, \ldots, m_{ik})$ denotes the $i$th row
of the matrix $M$, $i = 1, \ldots, n$. The natural parameter is
$\theta=(\xi_1/\sigma^2,\ldots,
\xi_k/\sigma^2,1/\sigma^2)^{\intercal} =: (\vec \theta, \theta_{k+1})^{\intercal} \in \mathbb{R}^{k} \times
(0,\infty)$. Moreover, we have for $i=1,\ldots,n$ the functions
$h_{i,j}(x)=m_{ij}x$, $j=1,\ldots,k$, $h_{i,k+1}(x)=-x^2/2$ and
$\beta_i(\theta)=-\frac{1}{2}
\log(\theta_{k+1})+\frac{1}{2\theta_{k+1}}(\mathbf{m}_i
\vec{\theta})^2$.
We     need the Hessian of $$\overline{\beta}(\theta)=\frac{1}{n}\sum_{i=1}^n
\beta_i(\theta) = -\frac{1}{2} \log(\theta_{k+1}) +
\frac{1}{2n\theta_{k+1}} \vec{\theta}^{\,\intercal}M^\intercal M
\vec{\theta}.$$
Direct computations give  
\begin{align*}
                                                  \nabla_{\vec{\theta}}\overline{\beta}(\theta)=
                                               \frac{1}{2n\theta_{k+1}^2}M^\intercal
  M \vec{\theta}, \mbox{ and } \partial_{\theta_{k+1}} \overline{\beta}(\theta)=&
                                               -\frac{1}{2\theta_{k+1}}
                                               -
                                               \frac{1}{2n\theta_{k+1}^2}
                                               \vec{\theta}^{\,\intercal}M^\intercal
                                               M\vec{\theta}, 
\end{align*}
and the second order derivatives 
\begin{align*}  
   &  \nabla_{\vec{\theta},\theta_{k+1}}  \overline{\beta}(\theta)=
                                                             -\frac{1}{2
  n\theta_{k+1}^2}M^\intercal M \vec{\theta},  \quad 
    \nabla_{\vec{\theta},\vec{\theta}}  \overline{\beta}(\theta)=
  \frac{1}{2n\theta_{k+1}}M^\intercal M,  \\
  & \mbox{ and }
     \nabla_{\theta_{k+1},\theta_{k+1}}  \overline{\beta}(\theta)=
  \frac{1}{2\theta_{k+1}^2}+\frac{1}{n\theta_{k+1}^3}\vec{\theta}^{\,\intercal}M^\intercal
  M \vec{\theta}.
\end{align*}
{We
standardise around the MLE}
\begin{align*}
    \hat{\theta}=  \frac{n}{{r_-}} ((M^\intercal M)^{-1}M^\intercal
  X,1)^{\intercal} \in \mathbb{R}^{k+1},
\end{align*}
for which 
\begin{align*}
 \mathbf{H}_{\overline{\beta}}(\hat \theta)&  =  \frac{r_{-}}{2n^2}
  \begin{pmatrix}
    M^{\intercal}M & - M^{\intercal} X \\
  -  X^{\intercal} M & r_{+}
\end{pmatrix}, \\
                     \mathbf{F}_{\overline{\beta}}(\hat \theta) & =  \frac{2n^2}{r_-}
  \frac{1}{X^{\intercal} X}
  \begin{pmatrix}
   (M^{\intercal}M)^{-1}X^{\intercal} X  +  (M^{\intercal}M)^{-1} M^\intercal
   X X^{\intercal} M  (M^{\intercal}M)^{-1} &  (M^{\intercal}M)^{-1}
   M^{\intercal} X \\
X^{\intercal} M  (M^{\intercal}M)^{-1} &  1
  \end{pmatrix},
\end{align*}
with  $r_{\pm} = X^{\intercal} X \pm X ^{\intercal} M
(M^{\intercal} M)^{-1} M^{\intercal} X 
$. Recall that
$
   \widehat{\matr K} = \sqrt n  {\matr K}_{\overline{\beta}}(\hat \theta)
$
 and define all other matrices accordingly.

We now choose a conjugate prior of the form
\begin{align*}
  \eta(\theta)= \tau_0 \log(\theta_{k+1}) + \vec{\tau}^{\intercal} \vec{\theta}
  +\tau_{k+1} \theta_{k+1} + \frac{1}{2\theta_{k+1}}
  {\theta}^{\,\intercal} \bm \sigma {\theta},
\end{align*}
where $\tau_0 >-1$, $\tau_{k+1}>0$,
$\bm\sigma=(\sigma_{i,j})_{1\leq i,j \leq k}$ is a symmetric positive
definite matrix and
$\vec{\tau}=(\tau_1,\ldots,\tau_k)^{\intercal} \in \mathbb{R}^k$ with
$\tau_{k+1}-\vec{\tau}^{\, \intercal} \bm\sigma^{-1}\vec{\tau}/2 <0$. 
We compute
\begin{align*}
  \nabla_{\theta_{k+1}} \eta(\theta)= \frac{\tau_0}{\theta_{k+1}}
  +\tau_{k+1} - \frac{1}{2\theta_{k+1}^2}
  \vec{\theta}^{\,\intercal}\bm \sigma\vec{\theta}, \quad
  \nabla_{\vec{\theta}}\eta(\theta)=& \vec{\tau} +
                                        \frac{1}{2\theta_{k+1}^2}\bm
                                        \sigma \vec{\theta}, 
\end{align*}
as well as
\begin{align*}
   \nabla_{\vec{\theta},\theta_{k+1}}  \eta(\theta)=
                                                        -\frac{1}{\theta_{k+1}^2}\bm
                                                        \sigma
                                                        \vec{\theta},
\quad 
   \nabla_{\vec{\theta},\vec{\theta}}  \eta(\theta)=
                                                        \frac{1}{2\theta_{k+1}}\bm
                                                        \sigma,
                                                        \mbox{ and }
   \nabla_{\theta_{k+1},\theta_{k+1}}  \eta(\theta)= \frac{-\tau_0}{\theta_{k+1}^2}+\frac{1}{\theta_{k+1}^3}\vec{\theta}^{\,\intercal}\bm \sigma \vec{\theta}.
\end{align*}
      The third order
      partial derivatives of {$\lambda=n\overline{\beta}+\eta$} that are nonzero are given
      by
\begin{align*}
    \nabla_{\vec{\theta},\vec{\theta},\theta_{k+1}} \lambda(\theta)=& -\frac{1}{\theta_{k+1}^2}(M^\intercal M+\bm\sigma), \quad
    \nabla_{\vec{\theta},\theta_{k+1},\theta_{k+1}} \lambda(\theta)= -\frac{2}{\theta_{k+1}^3}(M^\intercal M+\bm\sigma)\vec{\theta}, \\
    \nabla_{\theta_{k+1},\theta_{k+1},\theta_{k+1}} \lambda(\theta)=& \frac{2\tau_0-n}{\theta_{k+1}^3}-\frac{3}{\theta_{k+1}^4}\vec{\theta}^{\,\intercal}(M^\intercal M+\bm\sigma)\vec{\theta}. 
\end{align*}
Let $\kappa=\Vert M+\bm\sigma \Vert_{\infty}$. With the inequalities $\vert x\vert\leq (1+x^2)/2$ and $\vert x\vert^{-3}\leq (x^{-2}+x^{-4})/2$ we know that the absolute value of each of the partial derivatives is bounded by
\begin{align*}
  T(\theta)= \bigg(\frac{1}{\theta_{k+1}^2} + \frac{1}{\theta_{k+1}^4}
  \bigg) \bigg( \frac{n+2\tau_0}{2} + 3\kappa \|\vec{\theta}\|_\infty^2 + \kappa \big(1+\|\vec \theta\|_\infty^2/2\big)  \bigg).
\end{align*}
Now let $\Theta_0= B( \hat \theta, \epsilon)$ be a ball centred
around $\hat \theta$ with radius $\epsilon=\hat{\theta}_{k+1}/2>0$. Integrating out the uniform component with $\theta$ fixed we get (computations done
with \texttt{Mathematica}) 
\begin{align*} 
  &     
    \mathbb{E}\big[                                         
    T\big(\hat \theta + Y_1Y_2(\theta - \hat \theta)\big)  \, | \, \theta \big]\\
  & \le    \frac{1}{3 \hat \theta_{k+1}^3}  \bigg( \frac{n+2\tau_0}{2}+ \frac{7}{2} \kappa \epsilon^2 + \kappa \bigg) 
    \bigg( \frac{3 (\hat \theta_{k+1}-\epsilon) + \hat \theta_{k+1}}{2
    (\hat \theta_{k+1}-\epsilon)^2} + (1+3 \hat \theta_{k+1}^2)
    \frac{ \log(\hat
    \theta_{k+1}/(\hat \theta_{k+1}-\epsilon))}{\epsilon} \bigg)\\
    &=n\frac{1}{3 \hat \theta_{k+1}^3}  \bigg( \frac{1}{2}+\frac{\tau_0}{n} + \frac{\kappa}{n}\bigg(\frac{7}{8}  \hat{\theta}_{k+1}^2+1\bigg) \bigg) 
    \bigg( \frac{5+2\log(2)}{\hat{\theta}_{k+1}} +6\log(2) \hat \theta_{k+1} \bigg)\\
  & = : nT^{\star}(\hat \theta_{k+1}).
\end{align*}
Since $\nabla_{\vec{\theta},\vec{\theta},\vec{\theta}}\lambda(\theta)=0$, we can apply the bound (\ref{lembd}) to improve the dependence on the dimension $d=k+1$ over directly applying the bound of Lemma \ref{sec:main-resultsv2}.
Finally,
$ \mathbb{P}[\theta \notin \Theta_0^{\star}] =
\mathbb{P}[\sum_{r=1}^{k+1}(\theta_r - \hat \theta_r)^2\ge \hat \theta_{k+1}^2/4]
\le(4/\hat\theta_{k+1}^2)\sum_{r=1}^{k+1} \mathbb{E}[(\theta_r - \tilde \theta_r)^2]$.
We conclude 
that
\begin{align}
&d_{\mathrm{TV}}(\theta^{\star}_{\mathrm{MLE}}, N)
\leq \sqrt{\frac{\pi }{2n}}\|\mathbf{F}_{\overline{\beta}}(\hat \theta)\|_{\mathrm{diag}}\bigg(n(k+1)
    T^{\star}(\hat
    \theta_{k+1})\sum_{r=1}^{k+1}  \sqrt{\mathbb{E}[    
    (\theta_{k+1}- \hat
    \theta_{k+1})^2 ]\mathbb{E}[    
    (\theta_r- \hat
    \theta_r)^2 ]}  \nonumber\\
    &\quad+\sum_{u=1}^{k+1}|\partial_u\eta(\hat\theta)|+(k+1)\|\matr{H}_{\eta}(\hat
    \theta)\|_\infty\sum_{v=1}^{k+1}\sqrt{\mathbb{E}[    
    (\theta_r- \hat
    \theta_r)^2 ]}\bigg)+\frac{4}{\hat\theta_{k+1}^2}\sum_{r=1}^{k+1}  \mathbb{E}[    
    (\theta_r- \hat
    \theta_r)^2 ]+r^{\star}(\Theta_0^{\star}) , \label{regbd}
\end{align}
where, as in the previous examples, the remainder $r^{\star}(\Theta_0^{\star})$ is negligible.

In terms of the dependence on the sample size and the dimension, the bound (\ref{regbd}) can be seen to be $O(k^2/\sqrt{n})$ provided (i) $X^\intercal X=\Theta(n)$, (ii) $r_-=\Theta(n)$, (iii) $\|(M^\intercal M)^{-1}M^\intercal X\|_\infty=\Theta(1)$ and (iv) $\|(M^\intercal M)^{-1}\|_\infty=\Theta(n^{-1})$; in each of parts (i)--(iv) we understand the order notation to additionally imply that $X^\intercal X$, $r_-$, $\|(M^\intercal M)^{-1}M^\intercal X\|_\infty$ and $\|(M^\intercal M)^{-1}\|_\infty$ do not grow with the dimension $k$. Under these conditions, $\|\mathbf{F}_{\overline{\beta}}(\hat \theta)\|_{\mathrm{diag}}=\Theta(1)$, and we also have that  $\mathbb{E}[(\theta_r- \hat\theta_r)^2 ]=O(n^{-1})$ for $r=1,\ldots,k+1$, from which it is readily seen that the bound (\ref{regbd}) is $O(k^2/\sqrt{n})$.

We conclude by arguing that conditions (i)--(iv) are natural.
Condition (i) is clearly natural. Conditions (ii) and (iii) are also natural since the MLE $\hat\theta$ is a consistent estimator. Finally, we consider condition (iv). For $i,j=1,\ldots,k$, $(M^\intercal M)_{i,j}=\sum_{\ell=1}^n m_{\ell,i}m_{\ell,j}$, where $m_{i,j}$ is the $(i,j)$-th element of the matrix $M$. Thus, for a typical design matrix $M$, the elements of $M^\intercal M$ are of order $n$, and so $\mathrm{det}(M^\intercal M)=\Theta(n)$, from which it would follow that $\mathrm{det}((M^\intercal M)^{-1})=\Theta(n^{-1})$, which would in turn imply that $\|(M^\intercal M)^{-1}\|_\infty=\Theta(n^{-1})$.

\section*{Acknowledgements}

AF and YS were funded in part by ARC Consolidator grant from ULB, as well as  FNRS Grant CDR/J.0200.24. RG was funded in part by EPSRC grants EP/Y008650/1 and UKRI068. GR was supported in part by EPSRC grants EP/T018445/1, EP/R018472/1, and EP/X0021951.  AF was also supported in part by EPSRC grants EP/T018445/1 and and EP/X0021951.

\small

\bibliographystyle{plain}

\normalsize
 
\appendix

\section{Elements of multivariable calculus} \label{sec:mult-calc}
Although the arguments are elementary, we rely on some cumbersome
facts about multivariate derivatives and multivariate Taylor
expansions with integral remainder. With $Y_1, Y_2, \ldots$ i.i.d.\
uniform on $(0,1)$, %
we set
\begin{align*}
  U_k = \prod_{j=1}^{k-1} Y_j^j \mbox{ and  }V_k = \prod_{j=1}^k Y_j.
\end{align*}
First recall that if $\phi : \mathbb{R}^d \to \mathbb{R}$ is
sufficiently differentiable then the first and second order Taylor
expansions of $\phi$ around $0$ with integral remainder read as
\begin{align*}
  \phi(t)          & = \phi(0) +
            \sum_{\ell_1=1}^d  
            \mathbb{E}\left[ \partial_{\ell_1} \phi(Y_1t
                     )\right]t_{\ell_1}  = \phi({0}) + \left\langle \mathbb{E}[ \mathcal{D}\phi(Y_1t
            )], t \right\rangle, \\
  \phi(t)
          & = \phi(0) +\sum_{\ell_1=1}^d\partial_{\ell_1}\phi(0) t_{\ell_1} +
            \sum_{\ell_1=1}^d  \sum_{\ell_2=1}^d 
            \mathbb{E}\left[ Y_1\partial^{2}_{\ell_1, \ell_2} \phi(
            Y_1 Y_2t )\right]t_{\ell_1} t_{\ell_2}\\
                      & = \phi(0) + \left\langle \mathcal{D}\phi(0), t \right\rangle + \left\langle \mathbb{E}[Y_1
            \mathcal{D}^2 \phi(Y_1 Y_2 t)], t^{\otimes 2}
            \right\rangle 
\end{align*}
(compare with \cite[p.20]{ferguson1996large}).  To generalise at
arbitrary order for  $r \in \left\{ 1, \ldots, k \right\}$ and $\bm
\ell^r = (\ell_1, \ldots, \ell_r)$ we set 
$t_{\bm \ell^r} =   \prod_{i=1}^r t_{\ell_i}$, $\partial^r_{\bm\ell^{r}} =  \partial^r_{\ell_{1}, \ldots,
    \ell_r}$,   $\sum_{\bm\ell^r
            = 1}^d = \sum_{\ell_1=1}^d \dots \sum_{\ell_r=1}^d$. Then 
\begin{align*}
  \phi(t) & =  \phi(0) +  \sum_{r=1}^k \frac{1}{r!}\sum_{\bm\ell^r
            = 1}^d  t_{\bm \ell^r}
            \partial^r_{\bm\ell^r}
            \phi(0)   + \sum_{\bm\ell^{k+1} = 1}^d 
            t_{\bm \ell^{k+1}} \mathbb{E} \left[ U_k 
            \partial^{k+1}_{\bm\ell^{k+1}}
            \phi(V_k t)  \right]\\
  & =  \phi(0) +  \sum_{r=1}^k \frac{1}{r!} \left\langle
    \mathcal{D}^r\phi(0), t^{\otimes r} \right\rangle + \left\langle 
    \mathbb{E}[U_k \mathcal{D}^{k+1} \phi(V_k t)],
    t^{\otimes k+1} \right\rangle
\end{align*}
(the tensor notation is much more efficient).

Consequently, if $\lambda : \mathbb{R}^d \to \mathbb{R}$ is
sufficiently differentiable and $\matr{J}$ is a symmetric $d\times d$
matrix then the chain rule combined with the above formulas yields,
after a change in the order of summation, for all $u = 1, \ldots, d$,
\begin{align}
  \partial_u\lambda( t_0 +  \matr{J}\, t )  &  =
                                              \partial_u\lambda(t_{0}) +
                                              \sum_{v=1}^d \mathbb{E}   [
                                              \partial^2_{u,
                                              v}\lambda(t_0+Y_1 \matr{J}\, 
                                              t) ]
                                              ( \matr J \, t)_v, \label{eq:tm1} \\
  \partial_u\lambda( t_0 + \matr{J}\, t) 
                                            & =
                                              \partial_u\lambda(t_{0})
                                              + \sum_{v=1}^d 
                                              \partial^2_{u,
                                              v}\lambda(t_0)  (\matr
                                              J \, t) _v 
 + 
                                              \sum_{v=1}^d  \sum_{w =1}^d 
                                              \mathbb{E}\left[ Y_1
                                              \partial^{3}_{u, v, w} \lambda(t_0 +
                                              Y_1 Y_2 \matr J t
                                              )\right] (\matr J t)_v
                                              (\matr J t)_v,\label{eq:tm2}
\end{align}
and, more generally, 
\begin{align}
   \partial_u\lambda( t_0 + \matr{J}\, t) 
                                            & =
                                              \partial_u\lambda(t_{0})
                                              + \sum_{r=1}^k \frac{1}{r!}\sum_{\mathbf{v}^r=1}^d 
                                              \partial^{r+1}_{u,
                                              \mathbf{v}^r}\lambda(t_0)  (\matr
                                              J \, t)_{\mathbf{v}^r} 
 + 
                                              \sum_{\mathbf{v}^r=1}^d 
                                              \mathbb{E}\left[ U_k
                                              \partial^{k+1}_{u, \mathbf{v}^k} \lambda(t_0 +
                                              V_k \matr J t
                                              )\right] (\matr J t)_{\mathbf{v}^{k}}.\label{eq:tmk}
\end{align}

\section{Properties  of exponential families in canonical form}
\label{sec:moments-expon-famil}

\subsection{Generalities}
\label{sec:generalities}
Let $X_{\theta}\in \mathbb{R}^{r}$ be a random variable drawn from an exponential family
in canonical form with density (with respect to some dominating
measure $\mu$)
\begin{align}\label{eq:canonexpfam}
  p (  {{x}} \, | \, \theta) 
  =   \exp \big(  \left\langle  
  h({x}), \theta \right\rangle  - \alpha(x) - \beta(\theta) 
  \big)
\end{align}
with $\theta = (\theta_1, \ldots, \theta_d)^{\intercal}$,
$h({x}) \in \mathbb{R}^d$, $\alpha(x), \beta(\theta) \in \mathbb{R}$;
we suppose that this density is well-defined, that the parameter
$\theta$ is identifiable, and that $p (\cdot \, | \, \theta)$ has
support $\mathcal{S}\subset \mathbb{R}^r$ some open set (not depending
on $\theta$). This is exactly the setting considered in
\eqref{eq:likili}.  By definition the log normalising constant
(a.k.a.\ cumulant generating function) is
\begin{equation*}
\beta(\theta) = \log \int_{\mathcal{S}}  \exp \big(   \left\langle  
{h}({{x}}), \theta \right\rangle   
-\alpha({x})\big) \mu(dx).
\end{equation*}
It then holds that we can take derivatives under the integral sign to
get (see \cite[Theorem 2.2]{brown})
\begin{align*}
  \mathcal{D}\beta(\theta) & = \mathbb{E}[h(X_{\theta})], \quad
                             \mathcal{D}^2\beta(\theta)  =
                             \mathrm{Var}[h(X_{\theta})], \quad
                             \mathcal{D}^3\beta(\theta)  =
                             \mathbb{E}[(h(X_{\theta}) -
                             \mathbb{E}[h(X_{\theta})])^{\otimes 3}] 
\end{align*}
(as claimed in Remark \ref{rk:main-results-1}) and, more generally,
\begin{align*}
    \mathcal{D}^k\beta(\theta)  & = \kappa_k[h(X_{\theta})],
\end{align*}
where $\kappa_k[h(X_{\theta})]$ denotes the 
$k$-th cumulant of $h(X_{\theta})$
for all $k \ge 1$.

\subsection{Taylor expansion of log likelihoods in standardised
  exponential families}\label{sec4.2}

\begin{lemma} \label{lma:fact2a2b}
  Let  $ \theta^{\star} \sim
  p^{\star}(\cdot \, | \, \mathbf{x})$ as in (\ref{eq:thetastargen}). 
Then   \begin{align}
\nabla \log p^{\star}(t \, | \, \mathbf{x})  & =  
       \matr{J}  
                                                 \left( n \bar
                                               h(\mathbf{x}) - \nabla
                                                 \lambda(\theta_0 +  \matr{J}  t
                                                )\right).   \label{eq:tstvcoo} 
                            \end{align}
\end{lemma}
\begin{proof}
  From (\ref{eq:pstarlambda}) we see that ($c$ here is an irrelevant
  constant that will disappear with the derivatives)
\begin{align*}
  \log   p^{\star}(t \, | \, \mathbf{x}) & =
                                           c
                                           +  n\left\langle 
                                          \overline{h}(\mathbf{x}),
                                           \theta_0 + \matr{J}\,
                                           t \right\rangle  
                                           - \lambda \big(\theta_0 +
                                           \matr{J}\, t\big)    =     c + n \sum_{v=1}^d
                                          \overline{h}_v(\mathbf{x}) (\matr{J}\,
                                           t)_v 
                                           - \lambda \big(\theta_0 +
                                           \matr{J}\, t\big)\\
                                         & =     c + n \sum_{u=1}^d \sum_{v=1}^d
                                          \overline{h}_v(\mathbf{x})  \matr{J}_{vu} t_u 
                                           - \lambda \big(\theta_0 +
                                           \matr{J}\, t\big). 
\end{align*}
From the chain rule we immediately obtain that the Stein operator
$T_{\theta^{\star}}$ for
$ \theta^{\star} \sim p^{\star}(\cdot \, | \, \mathbf{x})$ is given by
$T_{\theta^{\star}}f(t) = \Delta f(t) +\left\langle L(t), \nabla f(t)
\right\rangle $, where $L(t)$ has components
  \begin{align*}
L_u(t) =  \partial_u \log
  p^{\star}(t \, | \, \mathbf{x}) =                                                  \sum_{v=1}^d \matr J_{uv} 
                                                 \left( n \bar h_{v}(\mathbf{x}) - \partial_v
                                                 \lambda(                                                \theta_0 +  \matr{J}  t
)\right) 
  \end{align*}
  (we also use the symmetry of $\matr J$) which gives
  \eqref{eq:tstvcoo}.
\end{proof}

\begin{lemma}[Taylor expansions of the log likelihood
  ratio] \label{lem:diffexp} Let $Y_1, Y_2, \ldots$ be i.i.d.\ uniform on
  $(0,1)$, independent of all other randomness involved.  Let
  $ p^{\star}(\cdot \, | \, \mathbf{x})$ be given by
  \eqref{eq:pstarlambda}. Finally, let
  $t^{\star} = \matr K(t - \theta_0)$ with $t \in \mathbb{R}^d$.  The
  following hold, for all $\xi \in \mathbb{R}^d$.
    \begin{itemize}
    \item If $\lambda$ is twice differentiable around $\theta_0$, then
      for all $t$ sufficiently close to $\theta_0$,
\begin{align}
& \left\langle  \nabla \log p^{\star}(t^{\star} \, | \, \mathbf{x}) +t^{\star}, \, \xi
  \right\rangle  = \left\langle  \nabla \lambda(\theta_0) - n \bar
                                              h(\mathbf{x}) 
                                              ,       \,         
                                                \matr{J}  \xi \right\rangle   +  
                                              \mathbb{E} \left[ \left\langle
                                              \matr{H}_\lambda( r_1(t^{\star})
                                              )
                                              -  \matr{H} , \,  \left(
    \matr{J}  \xi\right) \otimes (\matr J t^{\star})  \right\rangle \right]  \label{eq:2dorder}
\end{align}
with $r_1(u) =       \theta_0 + Y_1 \matr J u$. 

\item 
If $\lambda$ is three times differentiable around $\theta_0$,  then
Lemma \ref{lma:3rdordertayl} holds.
\item If $\lambda$ is $k+1$ times differentiable (for
  $k \in \mathbb{N}_0$) around $\theta_0$, then for all $t$
  sufficiently close to $\theta_0$, 
    \begin{align}
      &  \left\langle  \nabla \log p^{\star}(t^{\star} \, | \, \mathbf{x}) +t^{\star}, \, \xi
        \right\rangle \nonumber {=}
        \left\langle  \nabla \lambda(\theta_0) - n \bar
                h(\mathbf{x}), \,                                         
                \matr{J} \xi \right\rangle  +                                                                       
\left\langle                                                                        
                                      \matr{H}_{\lambda}(                                         \theta_0) 
         - \matr{H}, \,   \left( \matr{J} \xi \right) 
                                                                        \otimes
        (\matr J t^{\star})
        \right\rangle \nonumber \\
      & \qquad + \frac{1}{2} \left\langle \mathcal{D}^3 \lambda(
 \theta_0),  \left( \matr{J} \xi
                                                                               \right)
                                                                               \otimes
                                                                              (\matr
        J t^{\star})^{\otimes  2}
        \right\rangle + 
        \cdots + \frac{1}{k!} 
        \left\langle \mathcal{D}^k \lambda(
 \theta_0),    \left( \matr{J} \xi \right) \otimes (\matr
        J t^{\star})^{\otimes (k-1)}
\right\rangle \nonumber \\
      &\qquad
   -   \left\langle  \mathbb{E} \left[ U_k \mathcal{D}^{k+1} \lambda(
r_k(t^{\star}))\right], 
\left( \matr{J} \xi \right) \otimes (\matr
        J t^{\star})^{\otimes (k+1)}
                                                                       \right\rangle
                                                                          \label{eq:kdorder},
   \end{align}
   where $U_k = \prod_{j=1}^{k-1} Y_j^j$ and
   $r_k(u) = \theta_0+V_k \matr J u$ with
   $V_k = \prod_{j=1}^k Y_j$.
       \end{itemize}
\end{lemma}

\begin{proof}
  Let, as in the previous proof,
  $L(t) = \nabla \log p^{\star}(t \, | \, \mathbf{x})$.  Recall the
  notation $\matr{H} = \matr{K}^{2}$ and $\matr{K} =
  \matr{J}^{-1}$. We control the expression
  \begin{align*}
      L(t) + t      &      =  \matr{J} \left(
                          n\overline{h}(\mathbf{x})  - \nabla
                          \lambda (\theta_0 +  \matr{J}\, t)  
 \right)+ t = \matr{J} \left(
                          n\overline{h}(\mathbf{x})  - \nabla
                          \lambda ( 
           \theta_0 +  \matr{J}\, t)  
+\matr K t \right)
  \end{align*}
  by expanding $\nabla \lambda ( \theta_0+\matr{J}\, t)$ around
  $t = 0$ using the material from Section
  \ref{sec:mult-calc}. Expanding to the first order through
  \eqref{eq:tm1} gives
\begin{align*}
  L(t) + t      &      =     \matr{J}
                 \left( n \overline{h}(\mathbf{x})- \nabla  \lambda(
                 \theta_0) \right)   - \matr{J}  \left(   \mathbb{E}   [
                                               H_{\lambda}(t_0+Y_1 \matr{J}\, 
                                              t) ]
                                              ( \matr J \, t) -
                  (\matr{K} t) \right)\\
   &  =     \matr{J}
                 \left( n \overline{h}(\mathbf{x})- \nabla  \lambda(
                 \theta_0) \right)   - \matr{J}  \left(   \mathbb{E}   [
                                               H_{\lambda}(t_0+Y_1 \matr{J}\, 
                                              t) ]
                                               -
                      \matr{K}^2 \right)( \matr J \, t).
\end{align*}
Replacing $t$ by $t^{\star} = \matr K(t - \theta_0)$ and introducing
the notation $R_1=\theta_0+Y_1 (t-\theta_0)$ yields
\begin{align*}
  L(t^{\star}) + t^{\star}
  &=    \matr{J}
                 \left( n \overline{h}(\mathbf{x})- \nabla  \lambda(
                 \theta_0) \right)   - \matr{J}  \left(   \mathbb{E}   [
                                               \matr H_{\lambda}(R_1) ]
                                               -
                      \matr{H} \right)(t - \theta_0).
\end{align*}
It follows, thanks to the symmetry of all matrices involved, that
\begin{align*}
  \left\langle  L(t^{\star}) + t^{\star}, \xi \right\rangle & = \left\langle \matr{J}
                 \left( n \overline{h}(\mathbf{x})- \nabla  \lambda(
                 \theta_0) \right) , \xi \right\rangle - \left\langle  \matr{J}  \left(   \mathbb{E}   [
                                               \matr H_{\lambda}(R_1) ]
                                               -
                                                              \matr{H} \right)(t - \theta_0), \xi \right\rangle \\
  & = \left\langle 
                n \overline{h}(\mathbf{x})- \nabla  \lambda(
                 \theta_0), \matr{J} \xi \right\rangle - \left\langle   \mathbb{E}   [
                                               \matr H_{\lambda}(R_1) ]
                                               -
                                                              \matr{H}
,(\matr{J}   \xi) \otimes(t - \theta_0)  \right\rangle. 
\end{align*}
This proves  \eqref{eq:2dorder}. 
Expanding to the next order through  \eqref{eq:tm2} we get 
\begin{align*}
  L_u(t) + t_u             
& =       ( \matr{J}
                 \big( n \overline{h}(\mathbf{x})- \nabla  \lambda(\tilde
                 \theta) \big) )_u  -
                                   (      \matr{J}
\big(  \matr{H}_{\lambda}(\theta_0)    - \matr{H} 
                                                                       \big)
                                                                              \matr{J}
                          t)_u\\
  & \qquad -  \sum_{v=1}^d \sum_{s_1=1}^d
\sum_{s_2=1}^d \matr J_{uv}   \mathbb{E} \left[ Y_1 \partial^3_{s_1,   s_2, v}
 \lambda(\theta_0+
    Y_1Y_2 \matr{J} t )  (\matr{J}t)_{s_1}(\matr{J}t)_{s_2}
                                                                       \right].
                                                                          \end{align*}
Replacing $t$ with $t^{\star}= \matr K (t - \theta_0)$ and introducing $R_2 =
\theta_0 + Y_1Y_2 (t - \theta_0)$ leads to
\begin{align*}
  L_u(t^{\star}) + t^{\star}_u             
& =       ( \matr{J}
                 \big( n \overline{h}(\mathbf{x})- \nabla  \lambda(\tilde
                 \theta) \big) )_u  -
                                   (      \matr{J}
\big(  \matr{H}_{\lambda}(\theta_0)    - \matr{H} 
                                                                       \big)
                                                                              (t
                                              - \theta_0))_u\\
  & \qquad -  \sum_{v=1}^d \sum_{s_1=1}^d
\sum_{s_2=1}^d \matr J_{uv}   \mathbb{E} \left[ Y_1 \partial^3_{s_1,   s_2, v}
 \lambda(R_2)  (t-\theta_0)_{s_1}(t- \theta_0)_{s_2}
                                                                       \right],
                                                                          \end{align*}
                                                                          and thus
\begin{align*}
  \left\langle  L(t^{\star}) + t^{\star}, \xi \right\rangle & 
= \left\langle 
                n \overline{h}(\mathbf{x})- \nabla  \lambda(
                 \theta_0), \matr{J} \xi \right\rangle - \left\langle   
                                               \matr H_{\lambda}(\theta_0) 
                                               -
                                                              \matr{H}
                                                              ,(\matr{J}
                                                              \xi)
                                                              \otimes(t
                                                              -
                                                              \theta_0)
                                                              \right\rangle
  \\
  & \qquad - \left\langle \mathbb{E} \left[ Y_1 \mathcal{D}^3
 \lambda(R_2) \right] ,( \matr J \xi)\otimes (t- \theta_0)^{\otimes 2} \right\rangle,
\end{align*}
which is \eqref{eq:3dorder}. 
The arbitrary order expansion \eqref{eq:kdorder}
 follows along the exact same lines from \eqref{eq:tmk}.
\end{proof} 
\section{Further proofs from Section \ref{ex:multi}
}
\label{sec:about-dichl-distr}


\subsection{Proof of Lemma \ref{sec:multnomoments}}\label{appendixd1}

\begin{proof}[Proof of Lemma \ref{sec:multnomoments}]
The cumulant generating function  for \eqref{eq:p2dir}
is
\begin{align*}
  \kappa( \bm \pi) = \sum_{j=1}^{k-1} \log \Gamma(\pi_j) + \log
  \Gamma \bigg(n+ \tilde \tau - \sum_{j=1}^{k-1} \pi_j\bigg) - \log \Gamma(n+\tilde\tau).
\end{align*}
With $\psi$ the digamma function, direct and straightforward computations yield
\begin{align*}
  &   \partial_i \kappa(\bm\pi) = \psi(\pi_i) -  \psi(\pi_k) , \quad   \partial_{ii} \kappa(\bm\pi) = \psi'(\pi_i) +  \psi'(\pi_k), \quad  \partial_{ij} \kappa(\bm\pi) = \psi'(\pi_k), \quad i\not=j, \\
   & \partial_i^3 \kappa(\bm \pi) = \psi''(\pi_i)  - \psi''(\pi_k), \quad   \partial_i^4 \kappa(\bm \pi) = \psi'''(\pi_i)  + \psi'''(\pi_k),
\end{align*}
for all $i, j = 1, \ldots, k-1$.
It follows that
\begin{align*}
  & \mu_i(\bm \pi) := \mathbb{E}[\theta_i] = \psi(\pi_i) - \psi(\pi_k), \quad \sigma^2_{i}(\bm \pi) := \mathrm{Var}(\theta_i) = \psi'(\pi_i) + \psi'(\pi_k),\\
  & \mu_{3, i}(\bm \pi) :=  \mathbb{E}[(\theta_i -\mu_i(\bm \pi))^3] =  \psi''(\pi_i)
    - \psi''(\pi_k), \\
      & \mu_{4, i}(\bm \pi) :=  \mathbb{E}[(\theta_i - \mu_i(\bm \pi))^4] =
        \psi'''(\pi_i)  + \psi'''(\pi_k)  + 3 \left( \psi'(\pi_i) + \psi'(\pi_k) \right)^2.
\end{align*}

The moments needed in the statement of the lemma are centred around the posterior
mode $\tilde \theta_i = \log (\pi_i) - \log (\pi_k)$. With the
notation $\rho(x) = \psi(x) - \log(x)$, we get
\begin{equation}\label{222}
 \mathbb{E}[(\theta_i - \tilde \theta_i)^2] =
\psi'(\pi_i) + \psi'(\pi_k) + (\rho(\pi_i) - \rho(\pi_{k}))^2=\sigma^2_i(\bm \pi) + \mathfrak b_i(\bm \pi)^2
\end{equation}
with
$\mathfrak b_i(\bm \pi) = \mu_i(\bm \pi) - \tilde \theta_i
=\rho(\pi_i) - \rho(\pi_{k})$ the bias. Similarly, we have that
\begin{align}\label{444}
\mathbb{E}[(\theta_i - \tilde \theta_i)^4]= \mu_{4, i}(\bm \pi)+4 \mathfrak b_i(\bm \pi)\mu_{3, i}(\bm \pi)+6\mathfrak b_i(\bm \pi)^2\sigma^2_i(\bm \pi)+\mathfrak b_i(\bm \pi)^4.  
\end{align}
In order to bound $\mathbb{E}[(\theta_i - \tilde \theta_i)^2]$ and $\mathbb{E}[(\theta_i - \tilde \theta_i)^4]$ we will make use of the following bounds for the digamma function (see \cite[Lemma 1]{gq13}). Let $k\in\mathbb{N}$. Then for $x>0$,
\begin{equation*}
-1/x<\psi(x)-\log(x)<-1/(2x)   
\end{equation*}
and
\begin{equation*}
\frac{(k-1)!}{x^k}+\frac{k!}{2x^{k+1}}<(-1)^{k+1}\psi^{(k)}(x)< \frac{(k-1)!}{x^k}+\frac{k!}{x^{k+1}}.   
\end{equation*}
Applying these inequalities and using that without loss of generality $\pi_k=\max\bm\pi$ we obtain the bounds
\begin{align*}
|\mathfrak b_i(\bm \pi)|&\leq\frac{1}{\pi_i}, \quad
\sigma^2_i(\bm \pi)\leq\frac{2}{\pi_i}+\frac{2}{\pi_i^2}, \\
|\mu_{3, i}(\bm \pi)|&\leq\frac{1}{\pi_i^2}+\frac{2}{\pi_i^3}, \quad \mu_{4, i}(\bm \pi)\leq\frac{12}{\pi_i^2}+\frac{28}{\pi_i^3}+\frac{24}{\pi_i^4}.
\end{align*}
Applying these bounds to the formulas (\ref{222}) and (\ref{444}) now yields the desired bounds for $\mathbb{E}[(\theta_i - \tilde \theta_i)^2]$ and $\mathbb{E}[(\theta_i - \tilde \theta_i)^4]$.
\end{proof}

\subsection{Derivation of explicit constants
}

\noindent{\bf{Part 1: Initial bounds
}}

\vspace{2mm}

\noindent To derive Wasserstein and total variation distance bounds with explicit constants, we require an explicit bound for $\mathbb{E}[(M(\tilde\theta))^4]$. In the following lemma, we derive a crude bound, which will suffice for our purposes.

\begin{lemma}\label{glem}Let $G(x)=2\Gamma(x+8)\Gamma(x-8)/(\Gamma(x))^2$.
Then, $G$ is a decreasing function of $x$ on $(8,\infty)$. Moreover, provided $\min\bm\pi>8$, we have that
\begin{align}\label{ggg0}
\mathbb{E}\big[(M(\theta))^4\big]\leq (k-1)G(\min\bm\pi).    
\end{align}  
\end{lemma}

\begin{proof}
Using the crude inequalities $\max_{1\leq i\leq k-1}x_i\leq\sum_{i=1}^{k-1}|x_i|$ and $e^{|x|}\leq e^x+e^{-x}$ we have that
\begin{align}\label{ggg1}
\mathbb{E}\big[(M(\theta))^4\big]\leq\sum_{u=1}^{k-1}\big(\mathbb{E}\big[e^{8(\theta_u-\tilde\theta_u)}\big]+\mathbb{E}\big[e^{-8(\theta_u-\tilde\theta_u)}\big]\big).    
\end{align}
Using formula for the the posterior as given in equation (\ref{eq:p2dir})
we have that, for $\pi_k>8$ (which is ensured since we assume that $\min\bm\pi>8$),
\begin{align*}
\mathbb{E}\big[e^{8\theta_u}\big]&=\int_{\mathbb{R}^{k-1}}e^{8x_u}\frac{\Gamma(n+\tilde\tau)}{ \prod_{j=1}^k
  \Gamma(\pi_j)} \frac{\mathrm{exp}( \sum_{j=1}^{k-1} \pi_j
  x_j)}{(1 + \sum_{j=1}^{k-1} e^{x_j})^{n+\tilde\tau}}\,dx_1\cdots dx_{k-1}\\
  &=\frac{\Gamma(\pi_u+8)}{\Gamma(\pi_u)}\frac{\Gamma(\pi_k-8)}{\Gamma(\pi_k)}\int_{\mathbb{R}^{k-1}}\frac{\Gamma(n+\tilde\tau)}{ \prod_{j=1}^k
  \Gamma(\tilde{\pi}_j)} \frac{\mathrm{exp}( \sum_{j=1}^{k-1} \tilde{\pi}_j
  x_j)}{(1 + \sum_{j=1}^{k-1} e^{x_j})^{n+\tilde\tau}}\,dx_1\cdots dx_{k-1}\\
  &=\frac{\Gamma(\pi_u+8)}{\Gamma(\pi_u)}\frac{\Gamma(\pi_k-8)}{\Gamma(\pi_k)},
\end{align*}
where $\tilde\pi_u=\pi_u+8$, $\tilde\pi_k=\pi_k-8$ and $\tilde\pi_\ell=\pi_\ell$ otherwise. Since $e^{-8\tilde\theta_u}=(\pi_k/\pi_u)^8$, we thus obtain that
\begin{align*}
\mathbb{E}\big[e^{8(\theta_u-\tilde\theta_u)}\big]=\frac{\Gamma(\pi_u+8)}{\pi_u^8\Gamma(\pi_u)}\frac{\pi_k^8\Gamma(\pi_k-8)}{\Gamma(\pi_k)}.
\end{align*}
Now, $\Gamma(x+8)/(x^8\Gamma(x))=\prod_{j=1}^7(1+j/x)$ is clearly a decreasing function of $x$ on $(8,\infty)$, and similarly $x^8\Gamma(x-8)/\Gamma(x)=\prod_{j=1}^8(1-j/x)^{-1}$ is clearly a decreasing function of $x$ on $(8,\infty)$. These monotonicity results allow us to infer that $G(x)$ is decreasing on the interval $(8,\infty)$, and moreover that
\begin{align}\label{ggg2}
\mathbb{E}\big[e^{8(\theta_u-\tilde\theta_u)}\big]\leq\frac{\Gamma(\min\bm\pi+8)}{(\min\bm\pi)^8\Gamma(\min\bm\pi)}\frac{(\min\bm\pi)^8\Gamma(\min\bm\pi-8)}{\Gamma(\min\bm\pi)}=\frac{1}{2}G(\min\bm\pi).
\end{align}
By a similar argument, one can show that $\mathbb{E}\big[e^{-8(\theta_u-\tilde\theta_u)}\big]\leq G(\min\bm\pi)/2$ for $\min\pi>8$.
Combining this inequality with inequalities (\ref{ggg1}) and (\ref{ggg2}) 
now yields the desired bound (\ref{ggg0}).
\end{proof}

On combining the bounds (\ref{cs1}) and (\ref{cs2}) together with  the bounds of Lemmas \ref{sec:multnomoments} and \ref{glem} we can now obtain explicit bounds on the quantities $I_1$, $I_2$ and $I_3$. Provided $\min\bm\pi>8$,
\begin{align*}
I_1&\leq 4(k-1)\bigg(1+\frac{3}{2\min\bm\pi}\bigg)=:(k-1)A(\min\bm\pi),\\ 
I_2&\leq 2\cdot12^{3/4}(k-1)\bigg(1+\frac{11}{3\min\bm\pi}+\frac{15}{4(\min\bm\pi)^2}\bigg)^{3/4}\frac{(k-1)^{1/4}H(\min\bm\pi)}{\sqrt{\min\bm\pi}}=:(k-1)B_k(\min\bm\pi),\\
I_3&\leq2\cdot12^{3/4}(k-1)\bigg(1+\frac{11}{3\min\bm\pi}+\frac{15}{4(\min\bm\pi)^2}\bigg)^{3/4}\frac{(k-1)^{3/4}H(\min\bm\pi)}{\sqrt{\min\bm\pi}}=:(k-1)C_k(\min\bm\pi),
\end{align*}
where $H(x):=(G(x))^{1/4}$. In deriving the bound for $I_3$, we used the basic inequalities $\max\bm\pi\leq\sum_{u=1}^k\pi_u= n+\tilde\tau$ and $n+\tilde\tau=\sum_{u=1}^k\pi_u\geq k\min\bm\pi$ to obtain the bound $\sqrt{\max\bm\pi}/(n+\tilde\tau)\leq 1/\sqrt{k\,\min\bm\pi}$.  We thus obtain the bounds
\begin{align*}
d_{\mathrm{Wass}}(\theta^{\star}_{\mathrm{MAP}}, N)\leq\frac{\sqrt{2}(k-1)^{3/2}}{\sqrt{\min\bm\pi}}D_k(\min\bm\pi), \quad d_{\mathrm{TV}}(\theta^{\star}_{\mathrm{MAP}}, N)\leq\frac{\sqrt{\pi}(k-1)}{\sqrt{\min\bm\pi}}D_k(\min\bm\pi),    
\end{align*}
where $D_k(\min\bm\pi)=A(\min\bm\pi)+B_k(\min\bm\pi)+C_k(\min\bm\pi)$. These bounds are valid for $\min\bm\pi>8$. Finally, we put these rather unpalatable bounds into a compact form. 

\vspace{3mm}

\noindent{\bf{Part 2: Compact total variation distance bounds}}

\vspace{2mm}

\noindent We now complete the derivation of the total variation distance bound given in Example \ref{ex:multiintro} with constant $\mathcal{C}_1=8.46$; that is, we will prove that
\begin{equation}
\label{constant1} d_{\mathrm{TV}}(\theta^{\star}_{\mathrm{MAP}}, N)\leq \frac{8.46k}{\sqrt{\min\bm\pi}}.    
\end{equation}
We note that we may assume that $\sqrt{\min\bm\pi}\geq 8.46k$, as otherwise the bound is trivial. Observe that this assumption also ensures that we satisfy the requirement that $\min\bm\pi>8$. We begin by verifying that inequality (\ref{constant1}) holds for $k=2$. In this case, $8.46k=16.92$, and we compute $D_2(16.92^2)=5.95640$. We also observe that since, for fixed $k$, $D_k(x)$ is a decreasing function of $x$, it follows that $D_2(\min\bm\pi)\leq D_2(16.92^2)=5.95640$. Therefore, for $k=2$, we get that 
\begin{align*}d_{\mathrm{TV}}(\theta^{\star}_{\mathrm{MAP}}, N)\leq \frac{\sqrt{\pi}(2-1)}{\sqrt{\min\bm\pi}}\times5.95640\leq \frac{10.5574}{\sqrt{\min\bm\pi}}= \frac{5.2787\times 2}{\sqrt{\min\bm\pi}}< \frac{8.46\times 2}{\sqrt{\min\bm\pi}},
\end{align*}
and so inequality (\ref{constant1}) holds for $k=2$. Similarly, we can verify that inequality (\ref{constant1}) holds for $3\leq k\leq51$. We provide here the details for $k=16$; we do so because, for $3\leq k\leq 51$, numerical computations confirm that $((k-1)/k)D_k((8.46k)^2)\leq (15/16)D_{16}((8.46\times 16)^2)$. In this case, $D_{16}((8.46\times 16)^2)=5.08787$ so that  
\begin{equation*}d_{\mathrm{TV}}(\theta^{\star}_{\mathrm{MAP}}, N)\leq \frac{\sqrt{\pi}(16-1)}{\sqrt{\min\bm\pi}}\times5.08787=\frac{8.45438\times 16}{\sqrt{\min\bm\pi}}< \frac{8.46\times 16}{\sqrt{\min\bm\pi}},
\end{equation*}
as required. Finally, we consider the case $k\geq52$. In this case, we may assume that $\sqrt{\min\bm\pi}\geq 8.46k\geq 8.46\times 52=439.92$, and we get that
\begin{align*}
D_k(\min\bm\pi)&\leq D_k(439.92^2)\\&\leq4\bigg(1+\frac{3}{2\times439.92^2}\bigg)+2\cdot12^{3/4}\bigg(1+\frac{11}{3\times439.92^2}+\frac{15}{4\times439.92^4}\bigg)^{3/4}\times\\
&\quad\times\bigg(\frac{H(439.92^2)}{8.46\times 52^{3/4}}+\frac{H(439.92^2)}{8.46\times 52^{1/4}}\bigg)=4.76870.
\end{align*}
Thus, for $k\geq52$,
\begin{equation*}d_{\mathrm{TV}}(\theta^{\star}_{\mathrm{MAP}}, N)\leq \frac{\sqrt{\pi}k}{\sqrt{\min\bm\pi}}\times4.76870=\frac{8.45231k}{\sqrt{\min\bm\pi}}< \frac{8.46k}{\sqrt{\min\bm\pi}},
\end{equation*}
as required. We have therefore proved that inequality (\ref{constant1}) holds for all $k\geq2$. 


\vspace{3mm}

\noindent{\bf{Part 3: Compact Wasserstein distance bounds}}

\vspace{2mm}

\noindent We now derive the compact bound 
\begin{equation}
\label{constant2} d_{\mathrm{Wass}}(\theta^{\star}_{\mathrm{MAP}}, N)\leq \frac{7.40k^{3/2}}{\sqrt{\min\bm\pi}}    
\end{equation}
that was stated in Example \ref{ex:multiintro}. We shall derive this bound under the assumption that $\sqrt{\min\bm\pi}\geq7.40/\sqrt{2}$. We proceed as in part 2 of the proof, but with a more efficient exposition. The verification that inequality (\ref{constant2}) holds for $2\leq k\leq18$ follows by a similar approach as in part 2 of the proof. For $2\leq k\leq 18$,  numerical computations show that $((k-1)/k)^{3/2}D_k((7.40k/\sqrt{2})^2)\leq (16/17)^{3/2}D_{17}((7.40\times 17/\sqrt{2})^2)$.
In this case, $D_{17}((7.40\times 17/\sqrt{2})^2)=5.72873$ so that  
\begin{equation*}d_{\mathrm{Wass}}(\theta^{\star}_{\mathrm{MAP}}, N)\leq \frac{\sqrt{2}(17-1)^{3/2}}{\sqrt{\min\bm\pi}}\times5.72873=\frac{7.39741\times 17^{3/2}}{\sqrt{\min\bm\pi}}< \frac{7.40\times 17^{3/2}}{\sqrt{\min\bm\pi}},
\end{equation*}
as required. 
We verify that inequality (\ref{constant2}) holds for $k\geq19$. In this case, we may assume that $\sqrt{\min\bm\pi}\geq 7.40k/\sqrt{2}\geq 7.40\times 19/\sqrt{2}=99.4192$, and we get that
\begin{align*}
D_k(\min\bm\pi)&\leq D_k(99.4192^2)\\&\leq4\bigg(1+\frac{3}{2\times99.4192^2}\bigg)+2\cdot12^{3/4}\bigg(1+\frac{11}{3\times99.4192^2}+\frac{15}{4\times99.4192^4}\bigg)^{3/4}\times\\
&\quad\times\bigg(\frac{H(99.4192^2)}{7.40\times 19^{3/4}}+\frac{H(99.4192^2)}{7.40\times 19^{1/4}}\bigg)=5.22329.
\end{align*}
Thus, for $k\geq19$,
\begin{equation*}d_{\mathrm{Wass}}(\theta^{\star}_{\mathrm{MAP}}, N)\leq \frac{\sqrt{2}k^{3/2}}{\sqrt{\min\bm\pi}}\times5.22329=\frac{7.38684k^{3/2}}{\sqrt{\min\bm\pi}}< \frac{7.40k^{3/2}}{\sqrt{\min\bm\pi}},
\end{equation*}
as required. This completes the proof. \hfill\qed

\end{document}